\documentclass[a4paper,12pt]{article}

\usepackage[latin1]{inputenc} 
\usepackage[T1]{fontenc}      

\usepackage{geometry}         
\usepackage[francais, english]{babel}  
                              
\usepackage{amsmath}
\usepackage{amsfonts}
\usepackage{amssymb}
\usepackage{fancyhdr}
\usepackage{url}
\usepackage{mathrsfs}

\usepackage[all]{xy}

\newtheorem{df}{Definition}[section]
\newtheorem{prop}[df]{Proposition}
\newtheorem{thm}[df]{Theorem}
\newtheorem{thrm}{Theorem}
\newtheorem{lem}[df]{Lemma}
\newtheorem{crl}[df]{Corollary}
\newcommand{\prf}{\textit{Proof}}
\newtheorem{rmk}[df]{Remark}

\newcommand{\mom}{\mathcal{PM}_{\Omega}}
\newcommand{\tmom}{\widetilde{\mathcal{PM}}_{\Omega}}

\newcommand{\dbar}{\overline{\partial}}
\newcommand{\ddbar}{\partial\overline{\partial}}

\newcommand{\jbar}{\bar{\jmath}}
\newcommand{\kbar}{\bar{k}}
\newcommand{\R}{\mathbb{R}}

\newcommand{\C}{\mathbb{C}}

\newcommand{\N}{\mathbb{N}}
\newcommand{\XD}{X\backslash D}
\newcommand{\UD}{U\backslash D}

\newcommand{\ric}{\operatorname{ric}}
\newcommand{\tr}{\operatorname{tr}}

\newcommand{\vol}{\operatorname{vol}}

\newcommand{\vareps}{\varepsilon}

\newcommand{\vl}{\operatorname{Vol}}

\newcommand{\pr}{\operatorname{pr}}

\newcommand{\scal}{\mathbf{s}}
\newcommand{\energy}{\operatorname{\textbf{E}}}
\newcommand{\riem}{\operatorname{Rm}}

\title{\textbf{The space of Poincaré type Kähler metrics on the complement of a divisor}}  
\author{\textsc{Hugues AUVRAY}}
\date{}

\begin{document}

\makeatletter
\renewcommand%
   {\section}%
   {%
   \@startsection{section}%
      {1}%
      {0mm}%
      {\baselineskip}%
      {0.5\baselineskip}%
      {\sc\large\centering}%
   }%
\makeatother

\makeatletter
\renewcommand%
   {\subsubsection}%
   {%
   \@startsection{subsubsection}%
      {1}%
      {0mm}%
      {1.25\baselineskip}%
      {0.25\baselineskip}%
      {\sf\normalsize}%
   }%
\makeatother

\maketitle

\renewcommand{\abstractname}{Abstract}
 \begin{abstract}
  Consider a divisor $D$ with simple normal crossings in a compact Kähler manifold $X$. It has been known since the work by G. Tian and S.T. Yau that if $K[D]$ is ample there exists on $\XD$ a unique Kähler-Einstein metric with cusp singularities along the divisor (implying completeness and finite volume). 
  We show in this article that a Kähler metric in an arbitrary class, with constant scalar curvature and singularities analogous to that constructed by Tian and Yau, is unique in this class when $K[D]$ is ample.
This we do by generalizing Chen's construction of approximate geodesics in the space of Kähler metrics, and proving an approximate version of the Calabi-Yau theorem, both independently of the ampleness of $K[D]$. 
 \end{abstract}

 \section*{Introduction}

In the setting of compact Kähler manifolds, the existence of smooth geodesics for the Mabuchi metric between any two metrics among a fixed Kähler class is strongly related to, and in particular implies, uniqueness of canonical metrics like extremal metrics or constant scalar curvature metrics, up to the action of automorphisms of the identity component. 

Whereas it is now known that such smooth geodesics do not exist in general, see \cite{lv}, it is possible to construct less regular paths verifying the same equation in some more formal sense, and the difficulties arising from the lack of regularity can be bypassed, see e.g. \cite{chen1} when the canonical line bundle is ample, and \cite{ch-ti} in the general case, to give the expected uniqueness results. 
Such formal geodesics have nonetheless some regularity properties up to some of their second order derivatives (and higher order outside "small" sets), they are the unique such paths verifying the equation in question \cite{ph-st}, and they can be approached in the appropriate topology by smooth paths verifying some perturbed equation.

In this direction, this article lies within the framework of generalizing the results of \cite{chen1} to the setting of Kähler metrics with cusp singularities along a divisor. 
Namely, let $(X, \omega_0, J)$ be a compact Kähler manifold of complex dimension $m$, in which we consider a divisor $D$ with \textit{simple normal crossings}; write its decomposition into smooth irreducible components as $D=\sum_{j=1}^N D_j$. 
Let us endow each line bundle $[D_j]$ with a smooth hermitian metric $|\cdot|_j$, and denote by $\sigma_j\in\mathcal{O}([D_j])$ a holomorphic section such that $D_j=\{\sigma_j=0\}$, $j=1,\dots,N$. Up to multiplying $|\cdot|_j$ by a positive constant or a smooth positive function for those $j$, we can assume that $|\sigma_j|_j\leq e^{-1}$ so that $\rho_j:=-\log(|\sigma_j|^2_j)\geq 1$ out of $D_j$; notice that $i\ddbar\rho_j$ extends to a \textit{smooth} real (1,1)-form on the whole $X$, the class of which is $2\pi c_1([D_j])$.
Let $\lambda$ be a nonnegative real parameter, and set $u_j:=\log(\lambda+\rho_j)=\log\big(\lambda-\log(|\sigma_j|^2_j)\big)$. 
Choose $A_1,\dots, A_N>0$, and increase $\lambda$ if necessary; then,  
 \begin{equation}  \label{def_omega}
  \omega:=\omega_0-i\ddbar\mathfrak{u}
         =\omega_0-\sum_{j=1}^N\Big(A_ji\ddbar u_j\Big), \text{ where }\mathfrak{u}=\sum_{j=1}^N A_j u_j,
 \end{equation}
defines a genuine Kähler form on $\XD$, that we will take as a reference metric in what follows. 
This because if $U$ is a polydisc of coordinates $(z_1,\dots, z_m)$ around some point of $D$ such that $U\cap D=\{z_1\cdots z_k=0\}$, then $\omega$ is mutually bounded near the divisor with $\sum_{j=1}^k \tfrac{idz_j\wedge d\overline{z_j}}{|z_k|^2\log^2(|z_k|^2)}+\sum_{j=k+1}^m idz_j\wedge d\overline{z_j}$, and moreover has bounded derivatives at any order with respect to some orthonormal frame for this local model metric. 
In the same way, we shall look at $\mathfrak{u}$ as a reference potential. 
Indeed, we state:
 \begin{df}  \label{df_poinc_met}
  Let $\varpi$ be a locally smooth closed real (1,1) form on $\XD$. 
We say that $\varpi$ is a \textit{Kähler metric of Poincaré type in the class $\Omega=[\omega_0]_{dR}$}, denoted by $\varpi\in\mom$, if:
  \begin{itemize}  \setlength{\itemsep}{0pt}
   \item[(1)] $\varpi$ is $C^{\infty}$-quasi-isometric to $\omega$, meaning that $c\omega\leq\omega'\leq c^{-1}\omega$ on $\XD$ for some $c>0$, and $|\nabla_{\omega}^j\varpi|_{\omega}$ is bounded for any $j\geq1$;
   \item[(2)] $\varpi=\omega_0+i\ddbar v$ for some $v$ locally smooth on $\XD$ (the potential) such that $v=O(\mathfrak{u})$ near $D$, and $|\nabla_{\omega}^j v|$ is bounded on $\XD$ for any $j\geq 1$.
  \end{itemize}
Similarly, we denote by $\tmom$ the space of potentials of such metrics. 
 \end{df}
One can extend the notion of Mabuchi metric, see section \ref{gmtry_sp_mtrcs}, and endow these spaces with a Riemannian structure the geodesics of which are of particular interest, as in the compact setting. 

The main result of this article is a partial resolution of the equation of such geodesics, and as an application of this partial resolution and of the construction of approximate geodesics, we give a uniqueness result for Kähler metrics of Poincaré type with constant scalar curvature, provided the line-bundle $K[D]$ is ample. 
Let us sum up these results in the following statements: 
 \begin{thrm}  \label{thm_intro_gdscs+unqnss}
  Let $X$ be a compact Kähler manifold and $D$ a divisor with simple normal crossings in $X$. 
Consider the space $\mom$ of Poincaré type Kähler metrics on $\XD$ relative to some Kähler class on $X$, endowed with its Mabuchi metric. 
Then any two potentials of metrics in this space can be joined by a continuous geodesic, which furthermore can be approached by $C^{\infty}$ deformations of the segment joining them. 
There exists some uniform control on these approximate geodesics, seen as paths between potentials: they and their first order derivatives (in space and time directions) are bounded, as well as their time-time, space-time and some of their space-space second order derivatives, namely their complex Hessians (Theorem \ref{thm_MA_hmg} and Corollary \ref{crl_eps_gdscs}).
 \end{thrm}
 \begin{thrm}
Under the same assumptions and if in addition $K_{X}[D]$ is ample on $X$, then any metric lying in the space considered in Theorem \ref{thm_intro_gdscs+unqnss} and with constant scalar curvature is unique (Theorem \ref{thm_unqnss_csck_met}).
 \end{thrm}

A key ingredient to proving the last point is the existence among the class of metrics $\mom$ we consider of a metric with negative (in some strong sense) Ricci form, which we state as: 
 \begin{thrm}
  Assume $K_X[D]$ is ample on $X$. Then there exists a metric $\varpi\in\mom$ such that $\varrho(\varpi)<-c\varpi$ for some $c>0$ (Theorem \ref{thm_negative_ricci}).
 \end{thrm}
Notice that such metrics lie in $\mom$ for any fixed Kähler class $\Omega$ on $X$, and not necessarily for $\Omega=k c_1(K_X[D])$ with $k>0$, as the Kähler-Einstein metric obtained by Tian and Yau in \cite{tian-yau1} would do. 
This construction typically requires an analogue of the celebrated Calabi-Yau theorem \cite{yau1, aub}, and ours can be stated as the resolution of some Monge-Ampère equation, \textit{independently of the ampleness of }$K_X[D]$. 
Even if we use the method of its proof rather than the theorem itself, let us quote it now:
 \begin{thrm}
  Let $\omega$ be a Kähler metric of Poincaré type, and $f\in C^{\infty}_{loc}(\XD)$ which is a $O(e^{-\nu\mathfrak{u}})$ at any order for some $\nu>0$, such that $\int_{\XD}(e^f-1)\vol^{\omega}=0$. 
Then there exists a function $\varphi\in C^{\infty}_{loc}(\XD)$ bounded at any order such that $\big(\omega+i\ddbar\varphi\big)^m=e^f\omega^m$ on $\XD$. 
 \end{thrm}

The class of metrics we consider calls for a few comments. 
A first interest of considering such Kähler metrics with cusp singularities along a divisor is that this kind of singularities, at least when looking at the local model, can be reasonably well understood, and can allow canonical metrics involving a contribution of the divisor, like the Tian-Yau's Kähler-Einstein metric of \cite{tian-yau1}.  
Moreover, in the very active area devoted to the existence of Kähler-Einstein metrics on Fano manifolds, or more generally to the Tian-Yau-Donalson conjecture relying the existence of constant scalar curvature metrics among some integer Kähler class and algebro-geometric properties of the underlying polarized manifold, see for example \cite{don1}, there has recently been a renewed interest for certain type of singular metrics, see e.g. \cite{jmr, don2, rt}.  
Those singular metrics are roughly speaking those with conical singularities along a divisor or orbifold metrics, and one of the aims of considering them is to analyze smooth metrics by letting the cone angle go to $2\pi$. Since symmetrically cups singularities can be thought of as limits of conical singularities with cone angle going to 0, the study of the metrics we consider here certainly takes part in this growing theory, as a counterpart of the study of smooth metrics.

Now, let us give a few precisions on conditions (1) and (2) of Definition \ref{df_poinc_met}. 
First, condition (1) can somehow appear loose, and certainly having restricted our attention to metrics with a sharper asymptotic behaviour, as in \cite[\S 3.2]{sze}, would have made more precise some analytic considerations. 
Nonetheless, constructing approximate geodesics as in Theorem \ref{thm_MA_hmg} in such a restricted class (in other words, getting precise asymptotics for approximate geodesics) seems delicate, whereas working in $\mom$ had no notable drawback for this construction. 
Moreover, as it is not known whether one can ask a precise asymptotic behaviour for constant scalar curvature Kähler metrics of Poincaré type, it might be useful to have a rather general uniqueness result at one's disposal. 

On the other hand, condition (2) can appear a bit artificial, and one could think about replacing it by the more natural
 \begin{itemize}
  \item[(2')] $\varpi=\omega+d\psi$ for some 1-form $\psi\in L^2(\XD, \omega)$,
 \end{itemize}
and then deducing (2) from (1) and (2'). 
This is actually manageable, at least when the divisor is smooth. 
The main issue here is the \textit{boundedness of differentials} for potentials, required in proving Theorem \ref{thm_MA_hmg}, see section \ref{cont_met}, so we assume it via condition (2), which still allows us a sufficient range of examples. 
However, the control $v=O(\mathfrak{u})$ can always be performed assuming only (1) and (2'), see paragraph \ref{ctrl_potl_grth}.

~

The article is organized as follows. 
The first part of this work mainly deals with examples and analytic preliminaries in the setting of what we called Kähler metrics of Poincaré type. 
After recalling a model (and investigating more precisely its behaviour near the divisor) for such metrics (section \ref{model}), we focus on extending quickly a few notions required when dealing with the geometry of a space of Kähler metrics and constant scalar curvature metrics, like Aubin-Yau functional, or Mabuchi's metric and $K$-energy. 
This leads us to the equation of geodesics in our spaces of metrics; we can also look at it on potentials, and it writes, given a path $(v_t)_{(t\in[0,1])}$ joining $v_0$ and $v_1$ in $\tmom$,
 \begin{equation*}
  \ddot{v_t}-|\partial \dot{v_t}|^2_{\omega_{v_t}}=0,
 \end{equation*}
as in the compact case. 
Let us precise that we conclude part \ref{sp_of_mtrc} with the complementary section \ref{twds_ddbar_lemma}, where are generalized results like Poincaré inequality and is discussed the control one can get when widening the definition of Poincaré type metrics.

We formally solve the equation of geodesics in part \ref{resolution_mah} (Theorem \ref{thm_MA_hmg}), mostly adapting Chen and B\l ocki's techniques \cite{chen1, blocki} to our framework, particularly using a continuity method for a homogeneous Monge-Ampère equation on $(\XD)\times[0,1]\times S^1$ derived from \eqref{geodesic_eqn} (section \ref{cont_met}) similar to Chen's.
At that point, the difficulty occurred by the boundary at infinity $D$ adds to that from the boundary at finite distance of $[0,1]\times S^1$. 
We nonetheless bypass this thanks to the Poincaré setting, e.g replacing balls of coordinates by balls of local covering coordinates ("quasi-coordinates") in the $(\XD)$ direction when reasoning locally, or using an adapted maximum principle (Lemma \ref{lem_strong_max_principle}), in which the boundary on which is usually required some nonpositivity is replaced by $(\XD)\times(\{0\}\sqcup\{1\})\times S^1$. 

In the two central parts \ref{approx_CY} and \ref{prf_approx_CY} we come back on $\XD$. In the former, we state our logarithmic Calabi-Yau theorem (section \ref{approx_CY_motiv}), which we use in the next two sections to construct, assuming the ampleness of $K_X[D]$, metrics of Poincaré type with negative Ricci forms.   
A weighted $\ddbar$-lemma is also required for this construction, and is proved in section \ref{wtd_ddbar_lemma}. 
As for part \ref{prf_approx_CY}, it is devoted to the proof of the logarithmic Calabi-Yau theorem, following a rather classical progression ($C^0$, second, third, and higher order estimates in section \ref{prf_approx_CY_unif_bnds}) concerning the uniform control, the decay properties of approximate solutions being dealt with in section \ref{prf_approx_CY_wtd_spaces}.

As an application of the constructions of approximate geodesics and metrics with negative Ricci forms, we give in the final part the uniqueness result for Kähler metrics of Poincaré type with constant scalar curvature, provided the line-bundle $K[D]$ is ample (Theorem \ref{thm_unqnss_csck_met}). 
This uniqueness does not need to be considered up to the action of some $D$-parallel automorphism group of $X$, since as explained in Lemma \ref{lem_holo_vf}, there is no non-trivial holomorphic vector field $L^2$ for a Poincaré type metric when $K[D]$ is ample. 
We conclude (section \ref{prf_unqnss_thm}) by a short outline of the proof of the uniqueness theorem, which follows closely Chen's \cite{chen1} of the compact case.

~

A question to be studied would be the generalization of the results of \cite{ch-ti} to get rid of the ampleness of $K[D]$ in Theorem \ref{thm_intro_gdscs+unqnss}, and to get the uniqueness of constant scalar curvature metrics of Poincaré type up to the action of automorphisms in the connected component of the identity and tangent to the divisor. 
Another question one can ask is that of the existence of such metrics, and of a definition of $K$-stability for the pair $(X,D)$, see the suggestions of \cite[\S 3.2]{sze}.

  \subsubsection*{Acknowledgements}

I am very grateful to my PhD advisor Olivier Biquard for his constant encouragement and his well-advised recommendations, as well as his reading over a preparatory version of this paper. I also thank Paul Gauduchon and Sébastien Boucksom for the numerous useful conversations I had with them.

\section{The space of metrics} \label{sp_of_mtrc}
 
The purpose of this section is to give, by reference to a simple model, the definition and a few basic properties of Kähler metrics of Poincaré type. 

 \subsection{Model metric} \label{model}
 
  \subsubsection{Construction}  \label{constr_omega}
 
Recall quickly the construction of the model metric $\omega$; let $(X, \omega_0, J)$ be a compact Kähler manifold of complex dimension $m$, in which we consider a divisor $D$ with \textit{simple normal crossings} with decomposition $D=\sum_{j=1}^N D_j$ into smooth irreducible components. 
Take $\sigma_j\in\big(\mathcal{O}([D_j]),|\cdot|_j \big)$ a holomorphic defining section for $D_j$, $j=1,\dots,N$. 
We can assume that $\rho_j:=-\log(|\sigma_j|^2_j)\geq 1$ out of $D_j$; notice that $i\ddbar\rho_j$ extends to a \textit{smooth} real (1,1)-form on the whole $X$, whose class is $2\pi c_1([D_j])$.
Now let $\lambda$ be a nonnegative real parameter. If we set $u_j:=\log(\lambda+\rho_j)=\log\big(\lambda-\log(|\sigma_j|^2_j)\big)$, one has:
 \begin{lem}
  Let $A>0$. For sufficiently big $\lambda$ (depending on $A$ and $\omega_0$), the (1,1)-form $\omega_0-Ai\ddbar u_j$ defines a Kähler form on $X\backslash D_j$.
 \end{lem}
\prf. This comes from a simple computation; indeed, 
 \begin{equation*}
  -Ai\ddbar u_j=\frac{Ai\partial\rho_j\wedge\dbar\rho_j}{(\lambda+\rho_j)^2}-\frac{Ai\ddbar\rho_j}{\lambda+\rho_j}.
 \end{equation*}
The first summand is a nonnegative (1,1)-form, whereas $\pm\tfrac{Ai\ddbar\rho_j}{\lambda+\rho_j}\leq \tfrac{CA}{\lambda+\rho_j}\omega_0$ in the sense of (1,1)-forms where $C$ is such that $\pm i\ddbar\rho_j \leq C\omega_0$ on $X$. Since $\rho_j$ goes to $+\infty$ near $D_j$, we have $\omega_0-\tfrac{Ai\ddbar\rho_j}{\lambda+\rho_j}>0$ on $X\backslash D_j$ when $\lambda$ is big enough. $\square$

~

Choosing $A_1,\dots, A_N>0$, replacing $\omega_0$ by $\tfrac{1}{N}\omega_0$ and increasing $\lambda$ if necessary, one has $\tfrac{1}{N}\omega_0-i\ddbar u_j>0$ on $X\backslash D_j$ for $j=1,\dots,N$, hence
$\omega:=\omega_0-i\ddbar\mathfrak{u}=\sum_{j=1}^N\Big(\frac{1}{N}\omega_0-A_ji\ddbar u_j\Big)$ (recall $\mathfrak{u}=\sum_{j=1}^N A_j u_j$) defines a genuine Kähler form on $\XD$. 
Now we have checked this point, we shall also see why $\omega$ can be compared to a product (cusp metrics across the divisor)$\times$(smooth metric on the divisor), and even get precise asymptotics near $D$.

  \subsubsection{Asymptotic behaviour}  \label{asymp_prop}

Before describing the asymptotics of $\omega$ near $D$, let us fix the setting briefly. 
Around a codimension $k$ crossing, $D_1\cap\cdots\cap D_k$ say, consider an open set $U$ of holomorphic coordinates $(z_1,\dots,z_k,z_{k+1},\dots, z_m)\in \Delta^m$ with $\Delta\subset \C$ the open unit disc. 
The simple normal crossing assumption allows us to write $D\cap U=(D_1\cup\cdots\cup D_k)\cap U=\{z_1=0\}\cup\cdots\cup\{z_k=0\}$, that is $U\backslash D=(\Delta^*)^k\times\Delta^{m-k}$, $z_j=0$ being the equation of $D_j$ in $U$. We then have:
 \begin{prop}  \label{prop_asymptotics_omega}
  Set 
   \begin{equation*}
    \omega_{U,A}=\frac{A_1 idz_1\wedge d\overline{z_1}}{|z_1|^2\log^2(|z_1|^2)}+\cdots+\frac{A_k idz_k\wedge d\overline{z_k}}{|z_k|^2\log^2(|z_k|^2)}
                 +\Big(\omega_0-\sum_{j=k+1}^NA_ji\ddbar u_j\Big)|_{D_1\cap\cdots\cap D_k},
   \end{equation*}
  where the last summand on the right hand side is the metric $\omega_0-\sum_{k=j+1}^NA_ji\ddbar u_j$ restricted to $\Lambda^{1,1}_{D_1\cap\cdots\cap D_k}$. Then $\|\nabla^p_{\omega_{U,A}}(\omega-\omega_{U,A})\|_{\omega_{U,A}} =O(\rho_1^{-1}+\cdots+\rho_k^{-1})$ for all $p\geq 0$.
 \end{prop}
\prf. Let us start by the  $k=1$ and $p=0$ case. Notice that $|\sigma_1|^2_1=e^f |z_1|^2$ with some smooth (through $D$) $f$, thus $\rho_1=f+\log(|z_1|^2)\sim \log(|z_1|^2)$, $\partial\rho_1=\tfrac{dz_1}{z_1}+\partial f$ and $i\ddbar \rho_1=i\ddbar f$. As a consequence,
 \begin{align*}
  -i\ddbar u_1=\frac{idz_1\wedge d\overline{z_1}+i(z_1dz_1\wedge \dbar f+\overline{z_1}\partial f\wedge d\overline{z_1})+|z_1|^2i\partial f\wedge\dbar f} 
                    {|z_1|^2\rho_1^2(1+(\lambda+f)/\rho_1)^2}
               -\frac{i\ddbar f}{\rho_1+\lambda}.
 \end{align*}
As $\omega_{U,A}$ dominates $\omega_0$ and $i\ddbar f$ is smooth, $i\ddbar f$ is bounded i.e. $-\frac{i\ddbar f}{\rho_1+\lambda}$ is a $O(\rho_1^{-1})$ for $\omega_{U,A}$. Similarly, $df$ is bounded for $\omega_0$ hence for $\omega_{U,A}$, hence $\big|i(z_1dz_1\wedge \dbar f+\overline{z_1}df\wedge d\overline{z_1})\big|\leq C|dz_1|_{\omega_{U,A}}=CA_1^{-1/2}|z_1|^2\big|\log(|z_1|^2)\big|$, which gives a $O(\rho_1^{-1})$ after dividing by $|z_1|^2\log^2(|z_1|^2)$. 
Again we have $|z_1|^2i\partial f\wedge\dbar f=O(|z_1|^2)$, which gives $O(\rho_1^{-2})$ after dividing by $|z_1|^2\log^2(|z_1|^2)$, hence with respect to $\omega_{U,A}$, $-i\ddbar u_1=\frac{A_1 idz_1\wedge d\overline{z_1}}{|z_1|^2\log^2(|z_1|^2)}$ up to some $O(\rho_1^{-1})$. 

On the other hand $\omega'=\omega_0-\sum_{j\geq2}A_ji\ddbar u_j$ is smooth on $U$, hence $(\omega'_{j\kbar}-(\omega'|_{D_1})_{j\kbar})idz_j\wedge d\overline{z_k}=O(z_1)$ for $j,k\geq2$, which is easily a $O(\rho_1^{-1})$, and finally $\omega'_{1\kbar}idz_1\wedge d\overline{z_k}$ is a $O(|dz_1|_{\omega_{U,A}})$ and $\omega'_{1\bar{1}}idz_1\wedge d\overline{z_1}$ is a $O(|dz_1|^2_{\omega_{U,A}})$, with $|dz_1|_{\omega_{U,A}}=A_1^{-1/2}|z_1|\big|\log(|z_1|^2)\big|$, thus both are easily again $O(\rho_1^{-1})$.

We use the same technique when $p\geq 1$, and simply add the developments when the number of $D_j$ increases. \hfill $\square$

~

In the $k=1$ case, i.e. away from the crossings, these asymptotics clearly give a notion of what is the metric \textit{induced} by $\omega$ on any $D_j$ away from $\bigcup_{j'\neq j} D_{j'}$, and even on $D_j\backslash \bigcup_{j'\neq j} D_{j'}$, looking closer and closer to but not through $\bigcup_{j'\neq j} D_{j'}$. 
Then we can actually give a sharper formulation if there do exist crossings, using this notion of \textit{induced} metrics. 
For instance in the $k=2$ case, with $D_1$ (resp. $D_2$) given by $z_1=0$ (resp. $z_2=0$) we can show that $\omega=(\omega|_{D_2\backslash D_1})_{1\bar{1}}+(\omega|_{D_1\backslash D_2})_{2\bar{2}}+\omega'|_{D_1\cap D_2}+O(\rho_1^{-1}\rho_2^{-1})$ with $\omega'=\omega_0-\sum_{j>2}A_j i\ddbar u_j$ (smooth through $D_1\cup D_2$), and $O(\rho_1^{-1}\rho_2^{-1})$ understood at any order, which we shall denote by $O^{\infty}(\rho_1^{-1}\rho_2^{-1})$.
This gives rise to a recursive notion of metrics induced on $k$ codimensional crossings away from crossings of higher codimensions.

In this way, near the divisor, our metric is asymptotically a product of Poincaré metrics, or \textit{cusp} metrics, on punctured discs with a smooth metric, and we get back properties such as: $\omega$ is complete, has finite volume (equal to the volume associated to the class of $\omega_0$, cf. section \ref{gmtry_sp_mtrcs}), and its injectivity radius goes to 0. 
Notice also that the metrics induced on the (successive crossings of the) divisor have a similar behaviour.

About the standard cusp metric, notice the following, which is due to the homogeneity of Poincaré's half plane: let $\omega_{cusp}=\tfrac{idz\wedge d\overline{z}}{|z|^2\log^2(|z|^2)}$ on the punctured disc, and for $\delta\in(0,1)$ set $\varphi_{\delta}:\tfrac{3}{4}\Delta\rightarrow \Delta^*$, $\zeta\mapsto \exp\big(-\tfrac{1+\delta}{1-\delta}\tfrac{1+\zeta}{1-\zeta}\big)$. Then for \textit{all} $\delta\in(0,1)$, 
 \begin{equation*}
  {\varphi_{\delta}}^*\omega_{cusp}=\frac{id\zeta\wedge d\overline{\zeta}}{(1-|\zeta|^2)^2},
 \end{equation*}
which \textit{does not depend on} $\delta$ and is $C^{\infty}$-quasi-isometric to (mutually bounded with and with bounded derivatives at any order w.r.t. orthogonal systems for) the euclidian metric. 
Moreover ${\varphi_{\delta}}^*\log(|z|^2)=2\tfrac{1+\delta}{1-\delta}\tfrac{|\zeta|^2-1}{|1-\zeta|^2}$, which has the size of $\tfrac{1}{1-\delta}$ (with fixed factors). Besides $c\Delta^*\subset\bigcup_{\delta\in(0,1)}\varphi_{\delta}\big(\tfrac{3}{4}\Delta\big)$ with $c>0$ small enough ($0<c\leq e^{-25/7}$). This tells us for instance that in the $k=1$ case of the latter proposition, on the considered neighbourhood $U$, that there exists for all $p\in\N$ a constant $C_p$ such that
 \begin{equation*}
  \sup_{\delta\in(0,1)} \frac{1}{1-\delta}\Big\|\nabla_{euc}^p\Big({\Phi_{\delta}}^*\omega- 
                                     A_1\omega_{cusp}-\big(\omega_0-\sum_{j=2}^NA_ji\ddbar u_j\big)|_{D_1\cap\cdots\cap D_k}\Big)\Big\|_{euc}\leq C_p
 \end{equation*} 
with $\Phi_{\delta}: \tfrac{3}{4}\Delta\times \Delta^{m-1}\rightarrow \Delta^*\times\Delta^{m-1}$, $(\zeta_1,z_2,\dots,z_m)\mapsto \big(\varphi_{\delta}(\zeta_1),z_2,\dots,z_m\big)$.

  \subsection{Metrics of Poincaré type}  \label{class_ric_form}

The viewpoint of "quasi-coordinates" (the $\Phi_{\delta}$ are generally only holomorphic immersions, not charts) --- see \cite{tian-yau1}, p.580 --- is useful to define likewise Hölder spaces, for which the usual way of defining them is quite inconvenient because of the injectivity radius going to 0. 
Thus, if $\Phi_{\delta}: \mathcal{P}_k:=\big(\tfrac{3}{4}\Delta\big)^{k}\times \Delta^{m-k}\rightarrow (\Delta^*)^k\times\Delta^{m-k}$, $(\zeta_1,\dots,\zeta_k, z_{k+1},\dots,z_m)\mapsto \big(\varphi_{\delta_1}(\zeta_1),\dots,\varphi_{\delta_k}(\zeta_k), z_{k+1},\dots,z_m\big)$ for $\delta\in(0,1)^k$ and if $U$ is a polydisc neighbourhood of a  crossing $D_1\cap\cdots\cap D_k$ with $U\cap D_j$ given by $\{z_j=0\}$, $j=1,\dots,k$, we define for $f\in C^{p,\alpha}_{loc}(U\backslash D)$, $(p,\alpha)\in\N\times[0,1[$, 
 \begin{equation*}
  \|f\|_{C^{p,\alpha}(U\backslash D)}=\sup_{\delta\in(0,1)^k}\big\|{\Phi_{\delta}}^*f\big\|_{_{C^{p,\alpha}(\mathcal{P}_k)}},
 \end{equation*}
assuming that $U\subset (c\Delta)^k\times\Delta^{m-k}$. Then given a finite number of such open sets $U\in\mathcal{U}$ covering $D$, $V$ such that $X=V\cup\bigcup_{U\in\mathcal{U}}U$ and a partition of unity $\{\chi_V,\chi_U, U\in\mathcal{U}\}$, we define the Hölder space
 \begin{equation*}
  C^{p,\alpha}(\XD)=\big\{f\in C^{p,\alpha}_{loc}(\XD)|\, \|\chi_V f\|_{C^{p,\alpha}(V)}
                                                     +\sup_{U\in\mathcal{U}}\|\chi_U f\|_{C^{p,\alpha}(U\backslash D)}<+\infty\big\}
 \end{equation*}
endowed with the obvious norm, for the same $(p,\alpha)$ as above. We set similarly
 \begin{equation*}
  C^{\infty}(\XD)=\bigcap_{(p,\alpha)\in\N\times[0,1)}C^{p,\alpha}(\XD).
 \end{equation*}
Those spaces do not depend on the covering. 
In order to avoid ambiguity, let us precise that the Hölder spaces defined above with $\alpha=0$ are the same than the $C^k$ spaces defined with the help of the Levi-Civita of $\omega$, which strictly contain the spaces of $C^k_{loc}$ functions on $\XD$ with bounded derivatives up to order $k$ with respect to a smooth metric on the whole $X$.

We proceed similarly for tensors; for instance for $\varpi\in \Gamma^{p,\alpha}_{loc}(\Lambda^{1,1},U\backslash D)$ with $U$ as above, we set 
 \begin{equation*}
  \|\varpi\|_{\Gamma^{p,\alpha}(\Lambda^{1,1},U\backslash D)}= \sup_{\delta\in(0,1)^k}\|{\Phi_{\delta}}^*{\varpi}\|_{\Gamma^{p,\alpha}(\Lambda^{1,1},\mathcal{P}_k)}
 \end{equation*}
with the norm of ${\Phi_{\delta}}^*{\varpi}$ computed with the standard euclidian metric, and then we define $\Gamma^{p,\alpha}(\Lambda^{1,1},\XD)$ and $\Gamma^{\infty}(\Lambda^{1,1},\XD)$ by means of a partition of unity. 
An immediate example is $\omega\in \Gamma^{\infty}(\Lambda^{1,1},\XD)$. Let us restate the definition of the class of metrics we are investigating in this paper: 
 \begin{df}
  We say that a locally smooth real closed (1,1)-form $\omega'$ is a Kähler metric of Poincaré type in the class $\Omega$ of $\omega$, denoted by $\omega'\in\mom$, if:
   \begin{enumerate}  \setlength{\itemsep}{0pt}
    \item[1.] $\omega'$ is $C^{\infty}$-quasi-isometric to $\omega$, meaning that $c\omega\leq\omega'\leq c^{-1}\omega$ on $\XD$ for some $c>0$, and $\omega'$ has bounded derivatives in the quasi-coordinates used above (i.e. $\omega'$ is quasi-isometric to $\omega$ and $\omega'\in \Gamma^{\infty}(\Lambda^{1,1},\XD)$);
    \item[2.] $\omega'=\omega_0+i\ddbar v$ for some $v$ locally smooth on $\XD$ such that $v=O(\mathfrak{u})$ near $D$, and $dv$ has bounded derivatives at any nonnegative order in the quasi-coordinates used above (i.e. $dv\in \Gamma^{\infty}(\Lambda^{1},\XD)$).
   
   \end{enumerate}
Similarly, we denote by $\tmom$ the space of potentials (computed with respect to some fixed $\omega_{bp}\in\mom$ chosen as a base-point) of such metrics.
 \end{df}

 \begin{rmk}
Assuming condition 1. above, one can consider a class of metrics which is a priori wider, by relaxing condition 2. into 2'. $\omega'=\omega+d\psi$ for some 1-form $\psi\in L^2(\XD,\omega)$; we discuss this point in section \ref{twds_ddbar_lemma} below.
 \end{rmk}

We will take again the viewpoint of quasi-coordinates below, when looking at \textit{weighted} Hölder spaces, see section \ref{approx_CY_motiv}. 

We have defined an infinite dimensional manifold $\tmom$, which is a \textit{convex} open subset of the Fréchet space $\mathcal{E}$ of $C^{\infty}_{loc}$ functions on $\XD$ which are $O(\mathfrak{u})$ and with differential in $\Gamma^{\infty}(\Lambda^1,\XD)$; hence the tangent space $\tmom$ at any point is $\mathcal{E}$ itself. 
Since potentials are unique up to constants (their growth authorizes integrations by parts without boundary terms, use e.g. Gaffney-Stokes' theorem \cite{gaf} for complete manifolds), we take the normalization $\int_{\XD}f\vol^{\omega'}=0$ to fix the tangent space $\mathcal{E}/\R$ to $\mom$ at a point $\omega'$. 

Before going deeper into the geometries of $\mom$ and $\tmom$, let us observe the following fact, which contrasts with the compact setting ($K$ stands for the canonical line bundle $\Lambda^{m,m}$ on $X$ and $K[D]=K\otimes[D_1]\otimes\cdots\otimes [D_N]$):
 \begin{prop} \label{prop_class_ric}
  Let $\omega'$ be a metric 
of Poincaré type. 
Then its Ricci form $\varrho_{\omega'}$ is in $\Gamma^{\infty}(\Lambda^{1,1},\XD)$, and is $L^2(\XD,\omega)$-cohomologous to some (any) smooth real (1,1)-form of $-2\pi c_1(K[D])$.
 \end{prop}
\prf. By definition and according to Proposition \ref{prop_asymptotics_omega}, we can write $(\omega')^m=\tfrac{f}{\prod_{j=1}^N |\sigma_j|_j^2\log^2(|\sigma_j|^2_j)} \omega_0^m$ with $f\in C^{\infty}(\XD)$ and bounded below by some constant $c>0$, that is $\log f\in C^{\infty}(\XD)$. Thus on $\XD$, 
 \begin{equation*}
  \varrho_{\omega'}=\varrho_{\omega_0}+\sum_{j=1}^N i\ddbar\log(|\sigma_j|_j^2)+2\sum_{j=1}^N i\ddbar\log\big(\log(|\sigma_j|_j^2)\big)-i\ddbar\log f.
 \end{equation*}
Now for $j=1,\dots,N$, $i\ddbar\log(|\sigma_j|_j^2)$ extends to a smooth form of class $-2\pi c_1([D_j])$, so that $\varrho_{\omega_0}+\sum_{j=1}^N i\ddbar\log(|\sigma_j|_j^2)$ extends to a smooth form in $-2\pi c_1(K[D])$. 
Finally $\tfrac{1}{2}d^c\big(2\sum_{j=1}^N \log\big(\log(|\sigma_j|_j^2)\big) +\log f\big)$ is bounded hence $L^2$ for $\omega$, and its differential is $\Gamma^{\infty}(\Lambda^{1,1},\XD)$, hence the result. \hfill $\square$

 \subsection{Geometry of the spaces of metrics and of potentials} \label{gmtry_sp_mtrcs}

Let us now give a brief list of objects defined on $\mom$ and $\tmom$ necessary to our study, who share similar properties to their homonyms defined in the compact case (see \cite[ch.4]{gau} for a review of this subject):
 \begin{itemize} \setlength{\itemsep}{0pt}
  \item[$\bullet$] \textit{The volume.} We already know that the metrics $\omega_v=\omega_{bp}+i\ddbar v\in\mom$ have finite volume; in fact, they all have the same volume, $\vl$ say, which is also that attached to $\Omega$. 
Indeed, taking $\omega_{bp}=\omega$ (potentials are only translated in $\mathcal{E}$), it is easy to see that for any $v\in\tmom$, $\omega^m=(\omega_v)^m+d\Theta_v$ with $\Theta_v$ a polynomial in $\omega_v$, $d^c v$ and $i\ddbar v$, which is therefore bounded (for $\omega$ say), as well as its differential, and hence $\int_{\XD}d\Theta_v=0$ by Gaffney-Stokes. 
Since the $\mathfrak{u}\in \mathcal{E}$, we get in the same way that $\vl=\tfrac{1}{m!}[\omega_0]^m=\tfrac{1}{m!}\Omega^m$.

  \item[$\bullet$] \textit{Aubin-Yau functional} $\mathscr{J}$. 
For the same reason linked to integrations by parts as for the volume, the 1-form $\widetilde{\vol}:v \mapsto \big\{ f \mapsto \frac{1}{\vl} \int_{\XD} f\vol^{\omega_v}\big\}$ defined on $\tmom$ is closed, hence gives rise to a functional called $\mathscr{J}$ on $\tmom$ we fix saying it vanishes at 0. 
Moreover $\mathscr{J}$ is $\R$-equivariant (compute $\mathscr{J}(v)$ using the path $(tv)_{t\in[0,1]}$ from $0$ to $v$), thus we can identify $\mom$ to $\mathscr{J}^{-1}(0)$ in $\tmom$.

  \item[$\bullet$] \textit{Mabuchi metric.} One can give $\tmom$ and $\mom$ Riemannian structures writing $\langle f_1, f_2\rangle_v=\tfrac{1}{\vl}\int_{\XD}f_1f_2\vol^{\omega_v}$, or $\langle f_1, f_2\rangle_{\omega'}=\tfrac{1}{\vl}\int_{\XD}f_1f_2\vol^{\omega'}$ when $\int_{\XD}f_j\vol^{\omega'}=0$, $j=1,2$. 
Those metrics give an isometry $\tmom\ni v \mapsto \big(v-\mathscr{J}(v),\mathscr{J}(v)\big)\in \mom \times \R$.
  
  \item[$\bullet$] \textit{The equation of geodesics.} Denote by $v\in\mathcal{E}\big((\XD)\times[0,1]\big)$ a function on $(\XD)\times[0,1]$ such that for all $t$, $\dot{v}_t$ et $\ddot{v}_t\in\mathcal{E}$, or even $\big|\tfrac{d^k v}{dt^k}\big|\leq C_k\mathfrak{u}$ and $\big|\tfrac{d^k \nabla^jv}{dt^k}\big|\leq C_{k,j}$ for all $k\geq0$, $j\geq1$ and $\nabla$ the Levi-Civita connection of $\omega$ --- a segment in $\tmom$ clearly verifies those conditions.
Then such a function can be seen as a \textit{path} in $\tmom$ if $\omega_0+i\ddbar v_t\in\mom$ for all $t\in[0,1]$, and is a geodesic for the Mabuchi metric \textit{iff}
  \begin{equation}   \label{geodesic_eqn}
   \ddot{v}_t-|\partial\dot{v}_t|^2_{\omega_{v_t}}=0.
  \end{equation}
  
  \item[$\bullet$] \textit{The mean scalar curvature.} Once again integrations by parts work as in the compact case and tell us that the mean scalar curvature $\tfrac{1}{\vl}\int_{\XD} \scal(\omega')\vol^{\omega'}$ is the same for all $\omega'\in\mom$. 
If this quantity is denoted by $\overline{\scal}$, due to Proposition \ref{prop_class_ric}, one has 
 \begin{equation*}
  \overline{\scal}=-4\pi m\frac{c_1(K[D])\cdot[\omega_0]^{m-1}}{[\omega_0]^m}.
 \end{equation*}

  \item[$\bullet$] \textit{Mabuchi $K$-energies.} Again, the 1-form $\tilde{\scal}:v \mapsto\big\{ \tilde{\scal}_{v}:f\mapsto \int_{\XD} f\big(\scal(\omega_v)-\overline{\scal}\big)\vol^{\omega_v}\big\}$ is closed on $\tmom$ hence gives rise to a functional $\widetilde{\energy}$ that descends to a functional $\energy$ on $\mom$. 
They are called $K$-energies, and one can fix them saying they vanish at the base-points considered above. Their critical points are (potentials of) constant scalar curvature metrics.
Moreover, if one considers a path $(v_t)_{t\in[0,1]}\in\mathcal{E}\big((\XD)\times[0,1]\big)$ of potentials and sets $E:t\mapsto\widetilde{\energy}(t)$, then one can show:
 \end{itemize}

  \begin{prop} \label{prop_second_derivative_energy}
  For all $t\in[0,1]$, 
   \begin{equation}
    \begin{aligned}
    \ddot{E}(t)=2\|\nabla_{v_t}^- d\dot{v}_t\|^2_{L^2_{\omega_{v_t}}}
           -\int_{\XD} \big(\ddot{v}_t-|\partial \dot{v}_t|^2_{v_t}\big)(\scal_{v_t}-\overline{\scal})\vol^{\omega_{v_t}}.
    \end{aligned}
   \end{equation}
   where $\nabla_{v_t}^-$ is the $J$-anti-invariant part of the Levi-Civita connection of $\omega_{v_t}$ acting on $1$-forms.
  \end{prop}

The latter formula illustrates the importance of geodesics, because along such paths the $K$-energy would be convex.
Now take a path $(v_t)\in \mathcal{E}\big((\XD)\times[0,1]\big)$, and look at it as an element of $\mathcal{E}\big((\XD)\times[0,1]\times S^1\big)$ (similar definition) independent of the last variable, $s$ say, and set $\Phi(z,t,s)=v_t(z)$ for all $(z,t,s)\in(\XD)\times[0,1]\times S^1$. 
Give $\Sigma:=[0,1]\times S^1$ its natural complex structure. Then an easy computation, see e.g. \cite{sem}, or \cite{chen1}, p.197, gives:
 \begin{prop} \label{prop_semmes}
The path $(v_t)_{t\in[0,1]}$ is a geodesic iff 
 \begin{equation}  \label{MA_hmg}
   \big({\pr_{\XD}}^*\omega_{bp}+i\ddbar\Phi\big)^{m+1}\equiv 0,
 \end{equation} where operators $\partial$ and $\dbar$ are those of $(\XD)\times \Sigma$. 
In other words, the datum of a geodesic on $\tmom$ with extremities $v_0$ and $v_1$ is equivalent to that of a function $\Phi$ in $\mathcal{E}\big((\XD) \times \Sigma\big)$ which is $S^1$-invariant, which verifies equation \eqref{MA_hmg} and boundary conditions $\Phi(\cdot,\tau,\cdot)=v_{\tau}$, $\tau=0,1$, and such that for all $(t,s)$, $\Phi(\cdot,t,s)\in\tmom$. 
 \end{prop}

The solutions we can get for equation \eqref{MA_hmg}, and hence the "geodesics" we can get on $\mom$, are the purpose of the next part, to which the reader can jump directly, next section being devoted to complementary considerations for a priori more general metrics.

 \subsection{Analytic complements in Poincaré type metric} \label{twds_ddbar_lemma}

With the auxiliary aim %
of proving that metrics which are $C^{\infty}$-quasi-isometric to $\omega$ and in the same $L^2$ cohomolgy class are actually precisely those of $\mom$ when $D$ is smooth,  
we develop here a few basic tools for analysis in Poincaré type metric, some of which are also used in part \ref{prf_approx_CY}. 
For $k\geq1$, $\alpha\in[0,1]$, let us define
 \begin{equation} \label{eq_df_ekalpha}
  \mathcal{E}^{k,\alpha} = \big\{v\in C^{k,\alpha}_{loc}|\,v=O(\mathfrak{u}),\, dv\in\Gamma^{k-1,\alpha}(\Lambda^1) \big\}.
 \end{equation} 
The result we get in this section states as:
 \begin{prop} \label{prop_roughpnc=pnc}
  Assume $D$ smooth.
  Let $\eta\in \Gamma^{k,\alpha}(\Lambda^{1,1})$, $(k,\alpha)\in \N\times(0,1)$, an exact $(1,1)$-form one can write as $d\psi$ with $\psi\in L^2(\XD,\omega)$.
  Then there exists $v\in \mathcal{E}^{k+2,\alpha}$ such that $\eta=i\ddbar v$.
 \end{prop}
As an immediate corollary we have:
 \begin{prop}  \label{prop_roughpnc=pnc2}
  Assume $D$ smooth. 
  If $\omega'$ is $C^{\infty}$-quasi-isometric to $\omega$ and in the same $L^2$ cohomolgy class, then $\omega'$ writes as $\omega+i\ddbar\varphi$, with $\varphi\in\bigcup_{k,\alpha}\mathcal{E}^{k,\alpha}$, that is: $\omega'\in\mom$.
 \end{prop}
We first solve the equation $\Delta_{\omega} v=f$ with $f\in L^2$ and $\int_{\XD}f \vol^{\omega}=0$, and for this, establish a Poincaré inequality for (metrics quasi-isometric to) $\omega$ (\S \ref{poinc_ineq}). 
Then we take $f=-2\tr^{\omega}(\eta)$, and show that $i\ddbar v=\eta$ ; we also get the control $v=O(\mathfrak{u})$, and from classical elliptic theory, $v\in \mathfrak{u}C^{k+2,\alpha}(\XD)$ follows (\S \ref{prgrph_ddbar_lemma} and \ref{ctrl_potl_grth}). 
So far we do not need $D$ to be smooth, but assuming this we can improve regularity to get $v\in \mathcal{E}^{k+2,\alpha}$ (\S \ref{prgrph_prf_prop_roughpnc=pnc}).

  \subsubsection{Poincaré inequality}  \label{poinc_ineq}
 
  We consider a metric $g$ quasi-isometric to the model $\omega$ of the section \ref{model}. In order to solve in $H^1=H^1(\XD,g)$ the equation $\Delta_g v=f$, where $f$ is $L^2$ and has zero mean, by the classical variational method (minimization of the functional $\big\{v\mapsto\tfrac{1}{2}\int_{\XD}|dv|_g^2 \vol^g-\int_{\XD}vf\vol^g\big\}$ on zero mean functions), we show for $g$ a \textit{Poincaré inequality}:
   \begin{lem} \label{lem_unwtd_poinc_ineq}
  Assume $\XD$ is equipped with a metric $g$ quasi-isometric to the metric $\omega$ defined by (\ref{def_omega}). Then there exists a constant $C_P>0$ such that for all $v\in H^1(\XD,g)$ verifying $\int_{\XD}v\vol^g=0$ we have
  \begin{equation} \label{unwtd_poinc_ineq} \tag{PI}
   \int_{\XD}|v|^2\vol^g\leq C_P \int_{\XD} |dv|^2_g\vol^g.
  \end{equation}
 \end{lem}
\prf. Start, for simplicity, by the case where $D$ is smooth. We cover it in $X$ with open sets of coordinates $U_j$, $j=1,\dots,M$, of the form $\{|z|<a\}\times \Delta^{m-1}$, so that $D\cap U_j=\{|z|=0\}$. Consider also a neighbourhood $U$ of $D$ such that $U\subset \bigcup_{j=1}^M U_j$. Let $v\in C^{\infty}_c\big(\overline{U}\backslash D\big)$ such that $v|_{\partial U}\equiv 0$. We are first seeing there exists $c>0$ such that for all $j$, 
$\int_{U_j\backslash D}|v|^2\vol^g\leq c \int_{U_j\backslash D}|dv|^2_g\vol^g$. 
We can assume, up to modifying $c$, that $g$ restricted to $U_j\backslash D$ writes $\tfrac{4|dz|^2}{|z|^2\log^2(|z|^2)}+ds^2$, with $ds^2$ the euclidian metric on $\Delta^{m-1}$. Now change the coordinates by setting $t=\log(\log^2(|z|^2))\in(A,\infty)$ and $\theta=\arg z\in S^1$; $g$ becomes $dt^2+e^{-2t}d\theta^2+ds^2$, with volume form $e^{-t}dt d\theta ds$. Thus $\int_{U_j}|v|^2\vol^g=\int_{S^1\times\Delta^{m-1}}d\theta ds\int_{A}^{+\infty}|v|^2e^{-t}dt$ (resp. $\int_{U_j}|dv|^2_g\vol^g=\int_{S^1\times\Delta^{m-1}}d\theta ds\int_{A}^{+\infty}|dv|^2_ge^{-t}dt$), and we just need an inequality $\int_{A}^{+\infty}|v|^2e^{-t}dt\leq c \int_{A}^{+\infty} |dv|_g^2e^{-t}dt$ for all $(\theta,s)$ to conclude. 
Moreover, since $|dv|_g^2=(\partial_t v)^2+e^{2t}(\partial_\theta v)^2+|d_{\Delta^{m-1}}v|_{ds^2}^2\geq (\partial_t v)^2$, an inequality $\int_{A}^{+\infty}v^2e^{-t}dt\leq c \int_{A}^{+\infty}(\partial_t v)^2 e^{-t}dt$ for all $(\theta,s)$ still suffices to conclude. 

Set $w(t)=e^{-t}$; if $'$ stands for $\partial_t$, wen the have $(v^2w)'=2vv'w+v^2w'=2vv'w-v^2w$, hence by integrating with fixed $\theta$ and $s$, $0=2\int_{A}^{+\infty}vv'e^{-t}dt-\int_{A}^{+\infty}v^2e^{-t}dt$ because $v\equiv0$ on $\{t=A\}$ and for $t$ big enough. We rewrite this as:
 \begin{align*}
  \int_{A}^{+\infty}v^2e^{-t}dt =2\int_{A}^{+\infty}vv'e^{-t}dt
                                \leq 2\Big(\int_{A}^{+\infty}v^2e^{-t}dt\Big)^{\tfrac{1}{2}} 
                                      \Big(\int_{A}^{+\infty}v'^2e^{-t}dt\Big)^{\tfrac{1}{2}}
 \end{align*}
by Cauchy-Schwarz, hence $\int_{A}^{+\infty}v^2e^{-t}dt\leq 4\int_{A}^{+\infty}v'^2e^{-t}dt$, which ends the first point of demonstration.
We then have:
 \begin{equation} \label{nearD_poinc_ineq}
  \int_{\UD}|v|^2\vol^g\leq\sum_{j=1}^M\int_{U_j\backslash D}|v|^2\vol^g\leq c\sum_{j=1}^M\int_{U_j\backslash D}|dv|^2_g\vol^g\leq 
  Mc\int_{U\backslash D}|dv|^2_g\vol^g,
 \end{equation}
as soon as $v\in C^{\infty}_c\big(\overline{U}\backslash D\big)$. 

Now seek a contradiction, and take a sequence of functions $f_j\in C^{\infty}_c(\XD)$ violating the theorem; we thus can consider that 
 \begin{itemize}\setlength{\itemsep}{0pt} 
  \item for all $j$, $\int_{\XD}f_j\vol^g=0$ and $\int_{\XD}f_j^2\vol^g=1$;
  \item $\lim_{j\rightarrow \infty}\int_{\XD}|df_j|^2_g\vol^g=0$.
 \end{itemize}
Observe that $(f_j)$ is bounded in $H^1(\XD,g)$, hence up to an extraction converges weakly in $H^1(\XD,g)$ to a function $f\in H^1(\XD,g)$. 
In particular $\|df\|_{L^2(\XD,g)}=0$, that is to say $f$ is constant, since the $df_j$ tend to 0 in $L^2$.
Now finally, by weak $L^2$ convergence, $\int_{\XD}f\vol^g=\lim_{j\rightarrow \infty}\int_{\XD}f_j\vol^g=0$, hence $f\equiv 0$.

Take $\vareps>0$ small, such that $3\vareps^2<(Mc)^{-1}$ say, and a domain $V\subset\subset\XD$ wide enough so that $U^c\subset\subset V$ and there exists a smooth cut-off function $\chi$ equal to $1$ on $U^c$, $0$ on $V^c$, and such that $0\leq\chi\leq1$ et $|d\chi|_g\leq\vareps$. 
For all $j$ set $u_j=(1-\chi)f_j$ and $v_j=\chi f_j$ so that $u_j\in C^{\infty}_c\big(\overline{U}\backslash D\big)$, $(u_j)_{|\partial U}\equiv 0$, $v_j\in C^{\infty}_c(V)$ and $f_j=u_j+v_j$. Thus for all $j$,
 \begin{align*}
  \int_{\XD}f_j^2\vol^g\leq 2\Big(\int_{\XD}u_j^2\vol^g+\int_{\XD}v_j^2\vol^g\Big)
                       =2\Big(\int_{\UD}u_j^2\vol^g+\int_{V}v_j^2\vol^g\Big).
 \end{align*}
Now on the one hand, $(v_j)$ converges weakly to $0$ in $H^1\big(\overline{V},g\big)$ --- just see that for all test function $\varphi$ (resp. test 1-form $\alpha$) on $V$, $\chi\varphi$ is again a test function (resp. $\chi\alpha$ a test 1-form and $(d\chi,\alpha)_g$ a test function) --- 
and since $\overline{V}$ is compact with boundary, we can assume (forgetting another extraction) that $(v_j)$ strongly converges to 0 in $L^2$, necessarily to $0$.
 
On the other hand, according to the beginning of this demonstration, for all $j$ we have
 \begin{align*}
  \int_{\UD}u_j^2&\vol^g\leq Mc\int_{\UD}|du_j|_g^2\vol^g\\
                 &= Mc \left(\int_{\UD} \chi^2|df_j|_g^2\vol^g+\int_{\UD} f_j^2|d\chi|_g^2\vol^g +
                                   2\int_{\UD} f_j\chi(df_j,d\chi)_g\vol^g \right).
 \end{align*}
In the latter line, the first integral is bounded above by $\int_{\XD}|df_j|^2_g\vol^g$ which tends to 0; the second one by $\vareps^2\int_{\XD}f_j^2\vol^g=\vareps^2$, and the third by the square root of the first two.
It thus follows that $\int_{\XD}f_j^2\vol^g\leq 2Mc\vareps^2<1$ when $j$ is big enough, a contradiction, hence the theorem for $C_c^{\infty}(\XD)$ functions, and then for $H^1(\XD,g)$ functions by density. 

Now let us consider the case where $D$ admits crossings. If we have an inequality for smooth functions with a compact support near $D$ like (\ref{nearD_poinc_ineq}), the end of the argument will apply unchanged. To get this inequality though, cover $D$ with polydiscs of coordinates $\mathcal{P}_k=\{|z|<a_k\}^k\times \Delta^{m-k}$ ($a_k<1$ to adjust) such that $D$ is given in those by $\{z_1\cdots z_k=0\}$. 
One point is that to get the desired inequality with $U$ an open set relatively compact in the union of our polydiscs, it is enough to show such an inequality for functions $v\in C^{\infty}_c\big(\overline{\mathcal{P}_k}\backslash D\big)$ with $v\equiv 0$ on $\{|z_1|=a_k\}\cap\cdots\cap\{|z_k|=a_k\}$. 
But this we can do assuming $g$ is the product metric $\tfrac{4|dz_1|^2}{|z_1|^2\log^2(|z_1|^2)}+\cdots+\tfrac{4|dz_k|^2}{|z_k|^2\log^2(|z_k|^2)}+ds^2$, i.e. 
$dt_1^2+\cdots+dt_k^2+e^{-2t_1}d\theta_1^2+\cdots+e^{-2t_k}d\theta_k^2+ds^2$ where $t_l=\log(\log^2(|z_l|^2))\in(A_k,\infty)$, $\theta_l=\arg z_l\in S^1$, $l=1,\dots,k$. 
Finally, express $(t_1,\dots,t_k)$ in polar coordinates $(r,\varphi_1,\dots,\varphi_{k-1})$, $\varphi_1$, ..., $\varphi_{k-1}\in(0,\pi/2)$, $r\in(r(\varphi_1,\dots,\varphi_{k-1}),\infty)$, and do the same integrations by parts as above with $'$ standing for $\partial_r$ in order to conclude. \hfill $\square$

~

\subsubsection{Resolving $\Delta v=f$; a $\ddbar$-lemma} \label{prgrph_ddbar_lemma}

Take a metric $g$ quasi-isometric to the model metric $\omega$ of (\ref{def_omega}). As a corollary of Lemma \ref{lem_unwtd_poinc_ineq} every $f\in L^2$ with zero mean for $\vol^g$ admits a $H^1$ function $v$ such that $\Delta_g v=f$, unique as soon as $\int_{\XD}v\vol^g=0$. 
Moreover, $v$ is $H^2_{loc}$ by local ellipticity of $\Delta_g$ if one assumes more regularity on $g$. 
Actually:
 \begin{lem}[Sobolev estimate on $\XD$]  \label{lem_sobolev_estimate}
  If $g$ is $C^{\infty}$-quasi-isometric to $\omega$ and $\Delta_g v=f$ with $v\in H^1$, $f\in H^k$, $k\geq 0$, $\int_{\XD}f\vol^g=0$, then $v$ is in $H^{k+2}(\XD,g)$. If $\int_{\XD}v\vol^g=0$ then $\|v\|_{H^{k+2}}\leq C_k\|f\|_{H^{k}}$ for some constant $C_k$ depending only on $g$ and $k$.
 \end{lem}
\prf. Even if the idea of the proof is rather simple, it is more complicated to write it down completely in a brief way. 
Let us nonetheless give a few indications to see how it goes. 
First, for $v$ and $f$ as in the statement, an integration by parts shows 
$\|dv\|_{L^2_g}^2=\int_{\XD} vf\vol^g\leq \|v\|_{L^2_g}\|f\|_{L^2_g}$, so that if $\int_{\XD}v\vol^g=0$, Poincaré inequality (\ref{unwtd_poinc_ineq}) for $g$ gives $\|v\|_{L^2_g}\leq C_P\|f\|_{L^2_g}$. Secondly, since $\Delta_g$ is elliptic on any relatively compact domain $V$ in $\XD$, standard Sobolev estimates on balls tell us that $v\in H^{k+2}(V,g)$ and that there exists some $C_{V,k}$ such that $\|v\|_{H^{k+2(V,g)}}\leq C_{V,k}\big(\|f\|_{H^k(\XD,g)}+\|v\|_{L^2_g}\big)$, which is less than $(C_{V,k}+C_P)\|f\|_{H^k(\XD,g)}$ when $v$ has zero mean.

So that there remains to estimate the $L^2$ norm of $v$ on a neighbourhood of $D$. Assume for simplicity that $D$ is smooth and $k=0$, and suppose it is as usual covered by polydiscs of coordinates $U=(c\Delta)\times\Delta^{m-1}$ with $c$ a small constant, with $D$ given by $z_1=0$. Since $g$ is $C^{\infty}$-quasi-isometric to $\omega$, we can replace it by $g_{cusp}+ds^2$ on $U\backslash D$. Now since the pull-backs by the $\Phi_{\delta}$ introduced in \S \ref{constr_omega} of this latter metric are all same, $g_0$ say, the game is to express the $H^l_g$ norms on the $U\backslash D$ with the help of $H^l$ norms on the pullbacks. 
Namely, it is possible to find a sequence $(\delta_l)$ increasing to 1 and two constants $c_1$, $c_2$ such that for any $H^2_{loc}$ function $w$ on the considered sets, if $\mathcal{P}$ denotes the polydisc $\frac{3}{4}\Delta\times\Delta^{m-1}$,
 \begin{equation*}
  \|\nabla^j_g w\|_{L^2_g(c'U\backslash D)}\leq 
                                    c_1\sum_{l=1}^{\infty}\frac{1}{2^l}\|{\Phi_{\delta_l}}^*(\nabla^j_g w)\|_{L^2_{g_0}(\tfrac{1}{2}\mathcal{P})}
                                  = c_1\sum_{l=1}^{\infty}\frac{1}{2^l}\|\nabla^j_{g_0}({\Phi_{\delta_l}}^* w)\|_{L^2_{g_0}(\tfrac{1}{2}\mathcal{P})},
 \end{equation*}
$j=0,1,2$, i.e. $\| w\|_{H^2_g(c'U\backslash D)}\leq c_1\sum_{l=1}^{\infty}\frac{1}{2^l}\|{\Phi_{\delta_l}}^* w\|_{H^2_{g_0}(\tfrac{1}{2}\mathcal{P})}$ with $c'>0$ small independent of the covering, and conversely 
 \begin{equation*}
  \| w\|_{L^2_g(U\backslash D)}\geq c_2\sum_{l=1}^{\infty}\frac{1}{2^l}\|{\Phi_{\delta_l}}^* w\|_{L^2_{g_0}(\mathcal{P})}.
 \end{equation*}
Now, the standard Sobolev estimate on $\mathcal{P}$ for $g_0$ says there exists some constant $C>0$ such that for every $l$, $\|{\Phi_{\delta_l}}^* v\|_{H^2_{g_0}(\tfrac{1}{2}\mathcal{P})}\leq C\big(\|\Delta_{g_0}({\Phi_{\delta_l}}^* v)\|_{L^2_{g_0}(\mathcal{P})}+\|{\Phi_{\delta_l}}^* v\|_{L^2_{g_0}(\mathcal{P})}\big)$ $=C\big(\|{\Phi_{\delta_l}}^* f\|_{L^2_{g_0}(\mathcal{P})}+\|{\Phi_{\delta_l}}^* v\|_{L^2_{g_0}(\mathcal{P})}\big)$. 
Then take the weighted sum over $l$ with weights $\frac{1}{2^l}$ to get $\|v\|_{H^2_g(c'U\backslash D)}\leq c_1c_2^{-1}C\big(\|f\|_{L^2_g(U\backslash D)}+\| v\|_{L^2_g(U\backslash D)}\big)$. To conclude take enough of those $U$ so that $D$ is covered by the $c'U$, take $V$ wide enough and collect the inequalities. \hfill $\square$

~

We are now able to state a $\ddbar$-\textit{lemma} adapted to metrics "roughly" of Poincaré type:
 \begin{prop}[$\ddbar$-lemma  on $\XD$]  \label{lem_ddbar}
Any real square integrable exact $(1,1)$-form $\eta$ such that $\eta=d\psi$ with $\psi$ a $C^{\infty}_{loc}$ square integrable $1$-form writes $i\ddbar v$ with $v$ in $H^2\cap C^{\infty}_{loc}$, unique up to a constant.
 \end{prop}
\prf. This is classical. First, take $v$ as the only possible candidate (with zero mean), that is the solution of $\Delta_{\omega}v=-2\tr^{\omega}(d\psi)$. 
Then consider the 1-form $\xi:=\tfrac{1}{2}d^c v-\psi$. 
Since by construction, $\tr^{\omega}(d\xi)=0$, at every point one has the identity $d\xi\wedge d\xi\wedge\omega^{m-2}=-\tfrac{|d\xi|_{\omega}^2}{m(m-1)}\omega^m$. 
But the left hand side term can also be written $d(\xi\wedge d\xi\wedge\omega^{m-2})$, so by Gaffney-Stokes' theorem \cite{gaf} ($\xi\wedge d\xi\wedge\omega^{m-2}$ and $d(\xi\wedge d\xi\wedge\omega^{m-2})$ are $L^1$ since $v$ is $H^2$ for $\omega$ according to Lemma \ref{lem_sobolev_estimate}), its integral over $\XD$ is zero, hence $d\xi\equiv 0$ i.e. $i\ddbar v=\eta$. 
The only point to be verified is that $\int_{\XD}\tr^{\omega}(d\psi)\vol^{\omega}=0$, but this is guaranteed by the formula $\tr^{\omega}(d\psi)\vol^{\omega}=d\psi\wedge\tfrac{\omega^{m-1}}{(m-1)!}$ and one more use of Gaffney-Stokes' theorem. The local smoothness of $v$ is due to local ellipticity of $\Delta_{\omega}$, and actually this is a standard fact that for every $(p,\alpha)\in \N\times(0,1)$ and relatively compact domain $V\subset\subset W$, there exists $C=C(p,\alpha,V,W)$ such that $\|v\|_{C^{p+2,\alpha}(V)}\leq C(\|\eta\|_{L^2}+\|\eta\|_{C^{p,\alpha}(W)})$.   \hfill $\square$

 \subsubsection{Control on the potentials growth}  \label{ctrl_potl_grth}

Our $\ddbar$-lemma provides potentials for Kähler metrics of Poincaré type, in $H^2\cap C^{\infty}_{loc}$. 
Of course such potentials are not bounded in general (for example with $\alpha$ small enough in absolute value, $\omega+\alpha i\ddbar \mathfrak{u}$ is of Poincaré type whereas $\alpha \mathfrak{u}$ is not bounded --- recall that $\mathfrak{u}$ is defined by formula \eqref{def_omega}), we can still get some control on their growth near the divisor. 
 \begin{lem} \label{lem_control_of_potentials}
  Let $f\in C^{\infty}(\XD)$ have zero mean against $\vol^{\omega}$ --- for instance, $f=-2\tr^{\omega}(\omega'-\omega)$ with $\omega'$ a metric roughly of Poincaré type in class of $\omega$. 
Then if $v\in C^{\infty}_{loc}\cap H^2$ is a solution of $\Delta_{\omega}v=f$ --- in the example, $i\ddbar v=\omega'-\omega$ ---, there exists $C$ such that $|v|\leq C\mathfrak{u}$. 
Moreover, if $v$ also has zero mean, then one can take $C=C'\|f\|_{C^0(\XD)}$ with $C'$ depending only on $\omega$.
 \end{lem}
\prf. There is no loss in generality in assuming that $A_1=\cdots=A_N=2$ in defining formula (\ref{def_omega}). Now for $j=1,\dots,N$, take $\tilde{\lambda}\geq0$ and set $\tilde{u}_j=\log(\tilde{\lambda}+\rho_j)$ so that
 \begin{equation*}
  i\ddbar \tilde{u}_j=-\frac{i\partial\rho_j\wedge\dbar\rho_j}{(\tilde{\lambda}+\rho_j)^2}+\frac{i\ddbar\rho_j}{\tilde{\lambda}+\rho_j}.
 \end{equation*}
In view of Proposition \ref{prop_asymptotics_omega} and since $i\ddbar\rho_j$ is smooth through $D_j$, it is clear that given $\vareps>0$, when $\lambda$ is big enough then $\Delta_{\omega}\tilde{u}_j\geq-\vareps$ on $X\backslash D_j$, and $\Delta_{\omega}\tilde{u}_j=1+O(\rho_j^{-1})$ near $D_j$. 
So taking $\vareps$ small enough and $\tilde{\lambda}$ big enough ensures that there exist a neighbourhood $U$ of $D$ in $X$ and some constant $c>0$ (which we can take arbitrarily close to 1 after adjusting $\vareps$, $\tilde{\lambda}$ and $U$) such that $\Delta_{\omega}\tilde{\mathfrak{u}}\geq c$ on $U$, where $\tilde{\mathfrak{u}}=\sum_{j=1}^N\tilde{u}_j$. Notice that $\tilde{\mathfrak{u}}$ and $\mathfrak{u}$ are equivalent near $D$.

Now write $V_0=X\backslash U$, take domains $V_p$, $p\geq1$, such  that $(V_p)_{p\geq0}$ is an increasing exhaustive sequence of compact domains of $\XD$, and set finally $U_p=U\cap \mathring{V_p}$ for all $p\geq 0$. On the other hand, set $\varphi:=\pm v-C\tilde{\mathfrak{u}}-A$, where $C$ is chosen so that $\Delta_{\omega}\varphi=\pm f-C\Delta_{\omega}\tilde{\mathfrak{u}}\leq0$ on $U\backslash D$ (so $C$ depends only on $\|f\|_{C^0}$) and $A$ is chosen so that $\varphi\leq 0$ on $\partial U$ (so $A$ depends only on $\|f\|_{C^0}$ and on $\|v\|_{C^0(U)}$, which is controlled by $\|f\|_{C^0}$ provided $v$ has zero mean).

Consider for $p\geq 0$ the solution $\varphi_p$ of the Dirichlet problem 
 \begin{equation*}
 \left\{
  \begin{aligned}
   \Delta_{\omega}\varphi_p& =\Delta_{\omega}\varphi &\text{ on }& U_p\\
   \varphi_p               &=\varphi                 &\text{ on }& \partial U \\
   \varphi_p               &=0                       &\text{ on }& \partial V_p. \\
  \end{aligned}
 \right.
 \end{equation*}
By the usual maximum principle those $\varphi_p$ are nonpositive on their domains $U_p$. Suppose (some subsequence of) $(\varphi_p)_{p\geq0}$ converges almost everywhere to $\varphi$ ; then $\varphi\leq0$, i.e. $\pm v\leq C\tilde{\mathfrak{u}}+A$ and we are done. So we want to control the $\varphi_p$ in some Sobolev space in order to get some convergence in a smaller space.

Set $\theta_p=\varphi$ on $\partial U$ and $0$ on $\partial V_p$. The techniques used to show Lemma \ref{lem_sobolev_estimate} generalize to show that $\varphi_p$ is $H^2$ and there exists a constant $C$ \textit{independent of} $p$ such that 
 \begin{equation*}
  \|\varphi_p\|_{H^2(U_p)}\leq C\left(\|\Delta_{\omega}\varphi_p\|_{L^2(U_p)}+\|\varphi_p\|_{L^2(U_p)}+\|\theta_p\|_{L^2(\partial U_p)}\right).
 \end{equation*}
Now $\|\Delta_{\omega}\varphi_p\|_{L^2(U_p)}=\|f\|_{L^2(U_p)}\leq \|f\|_{L^2(\XD)}$, and $\|\theta_p\|_{L^2(\partial U_p)}=\|\varphi\|_{L^2(\partial U)}$, which do not depend on $p$ (and are controlled by $\|f\|_{C^0}$). It remains to estimate $\|\varphi_p\|_{L^2(U_p)}$. 
Decompose $\varphi_p$ into $\psi_p+\chi_p$ where $\psi_p\equiv0$ on $\partial U_p$ and $\chi_p$ is harmonic on $U_p$. Then $\chi_p$ is nonpositive and reaches its infimum on $\partial U_p$, so that $\|\chi_p\|_{L^2(U_p)}\leq \left|\inf_{\partial U}\varphi\right|\cdot\vl(U)^{1/2}$. 

Finally, $\int_{U_p} |d\psi_p|^2_{\omega}\vol^{\omega}=\int_{U_p} \psi_p\Delta_{\omega}\psi_p\vol^{\omega}=\int_{U_p} \psi_pf\vol^{\omega}$. But $\psi_p$ extends to an $H^1$ function on $\XD$ declaring it is 0 on $(\XD)\backslash U_p$, so that if $a_p$ is its mean on $\XD$,
 \begin{equation*}
  \int_{\XD}(\psi_p-a_p)^2\vol^{\omega}\leq C_P\int_{\XD}|d\psi_p|^2_{\omega}\vol^{\omega}=C_P\int_{U_p}|d\psi_p|^2_{\omega}\vol^{\omega}.
 \end{equation*}
As $\int_{\XD}(\psi_p-a_p)^2\vol^{\omega}=\int_{U_p}\psi_p^2\vol^{\omega}-a_p^2\vl(\XD)$ and $|a_p|\leq\tfrac{\vl(U_p)^{1/2}}{\vl(\XD)}\|\psi_p\|_{L^2(U_p)}$, we get, going back up those inequalities that
 \begin{equation*}
  \|\psi_p\|_{L^2(U_p)}\leq C_P\Big(1-\frac{\vl(U_p)}{\vl(\XD)}\Big)^{-1}\|f\|_{L^2(U_p)}\leq C_P\Big(1-\frac{\vl(U)}{\vl(\XD)}\Big)^{-1}\|f\|_{L^2(\XD)},
 \end{equation*}
which does not depend on $p$. So the $\varphi_p=\psi_p+\chi_p$ are $H^2$-bounded in their domains, and the bound, $C$ say, does not depend on $p$.

A diagonal extraction gives us the weak convergence in all the $H^2(U_p)$ and strong convergence in the $H^1(U_p)$ of a subsequence of $(\varphi_p)$ to some $\varphi'$ lying in $\bigcap_{p\geq 0}H^2(U_p)$. Moreover $\|\varphi'\|_{H^2(U\backslash D)}=\sup_p \|\varphi'\|_{H^2(U_p)}$, and each $\|\varphi'\|_{H^2(U_p)}$ is less or equal than the $\liminf$ of the $\|\varphi_q\|_{H^2(U_p)}$ when $q$ goes to $\infty$, quantity bounded by $C$, so: $\|\varphi'\|_{H^2(U)}\leq C<+\infty$. 
It is not hard to see that $\varphi'|_{\partial U}=\varphi|_{\partial U}$ and $\Delta_{\omega}\varphi'=\Delta_{\omega}\varphi$ on $U$ because the equality $\varphi_p|_{\partial U}=\varphi|_{\partial U}$ (resp. $\Delta_{\omega}\varphi_p=\Delta_{\omega}\varphi$ on $U_q$) holds for every $p$ (resp. every $p\geq q$). So $\varphi$ and $\varphi'$ are two $H^2(U)$ functions satisfying the same Dirichlet problem on $U$, so by $H^1(U)$ uniqueness, $\varphi'=\varphi$, that is:
$\varphi$ is (up to an extraction) the $L^2$-limit of $(\varphi_p)$ on any $U_q$, so (up to another extraction) $\varphi$ is almost everywhere in $U$ the limit of this sequence of nonpositive $(\varphi_p)$. 
\hfill $\square$

 \subsubsection{The smooth divisor case: proof of Proposition \ref{prop_roughpnc=pnc}} \label{prgrph_prf_prop_roughpnc=pnc}
 
We assume now that $D$ is smooth, and reduced to one component for sake of simplicity (what follows easily generalizes to the case when $D$ has several disjoint components). 
We start from the following fibration:
 \begin{equation} \label{eqfibr}
  \xymatrix{
    S^1 \ar[r] & \mathcal{N}_A\backslash D \ar[d]^{q=(t,p)} \\
               & [A,+\infty[\times D
  }
 \end{equation}
Let us explain it briefly. 
The tubular neighbourhood $\mathcal{N}_A$ of $D$, with projection $p$, is obtained from the exponential map of a smooth metric on $X$, e.g. $\omega_0$. 
The $S^1$ action comes from the identification of $\mathcal{N}_A$ with a neighbourhood $\mathcal{V}$ of the null section of the holomorphic tangent bundle $N_D=\tfrac{T^{1,0}X|_D}{T^{1,0}D}$, and leaves $p:\mathcal{N}_A\simeq\mathcal{V}\subset N_D\to D$ invariant. 
The part $t$ of the projection $q$ in \eqref{eqfibr} is obtained from $\mathfrak{u}=\log\big(-\log(|\sigma|^2)\big)$ we make $S^1$-invariant (we take for example the mean of $|\sigma|$ under the $S^1$ action) near the divisor and extended smoothly away; it is easy to see that $t=\mathfrak{u}$ up to a perturbation which is $O(e^{-t})$ as well as its derivatives at any order (for $\omega$). 
Finally, $A$ et $\mathcal{N}_A$ are adjusted so that $\mathcal{N}_A\backslash D=\{t\geq A\}\subset\XD$.

One associates to the circle action on $\mathcal{N}_A$ a connection 1-form $\eta$, as follow: if $g$ the metric associated to $\omega$ and $T$ the infinitesimal generator of the action, of flow $\Phi_s$, we set at any point $x$ of $\mathcal{N}_A$
 \begin{equation*}
  \hat{\eta}_x= \int_0^{2\pi}\Phi_s^*\Big(\frac{g_x(\cdot,T)}{g_x(T,T)}\Big) ds \quad \text{ et } \eta_x=2\pi\Big(\int_{S^1}\hat{\eta}\Big)^{-1}\hat{\eta}_x,
 \end{equation*}
where $S^1$ in the last integral is the fiber $q^{-1}(x)$. 
In this way, for all $x\in\mathcal{N}_A$, $\int_{q^{-1}(x)}\eta=2\pi$.

Moreover, if one considers around a point of $D$ a neighourhood of holomorphic coordinates $(z_1,\dots,z_m)$ such that $D$ is given $z_1=0$, one ha $\eta=d\theta$ up to a term which is $O(1)$ at any order for $\omega$. 
We then have:
 \begin{equation}  \label{eq_asmpt_g}
  g=dt^2+e^{-2t}\eta^2+p^{*}g_D+O\big(e^{-t}\big)
 \end{equation}
with $g_D$ the metric associated to $\omega_0|_D$, and the perturbation $O\big(e^{-t}\big)$ is understood at any order for $\omega$. 
This means for example that $Jdt=2e^{-t}\eta+O\big(e^{-t}\big)$, the $O\big(e^{-t}\big)$ understood as well. 

One can use furthermore the fibration \eqref{eqfibr} as follows. 
Let $f\in C^{k,\alpha}\big(\XD\big)$; we write the decompositions
 \begin{equation}  \label{eqdecomp1}
  f=(\Pi_0 f)(t,z)+\Pi_{\perp}f = f_0(t)+f_1(t,z)+\Pi_{\perp}f
 \end{equation}
where $z=p(x)$, with :
 \begin{equation*}
  (\Pi_0 f)(t,z)=\frac{1}{2\pi}\int_{q^{-1}(x)}f\,\eta \quad \text{et} \quad f_0(t)=\frac{1}{\vl(D)}\int_D f(t,z)\,\vol^{g_D},
 \end{equation*}
and $\vl(D)$ computed with $g_D$, hence equal to $\tfrac{[\omega_0|_D]^{m-1}}{(m-1)!}$, or $\tfrac{c_1([D])\cdot[\omega_0|_D]^{m-1}}{(m-1)!}$

Using \eqref{eqdecomp1} and the definition of $C^{k,\alpha}\big(\XD\big)$, since the fibers $S^1$ are of length equivalent to $e^{-t}$ for $g$, it is easy to see that on an open set of coordinates as above and for all $j\leq k$,
 \begin{equation*}
  \mathcal{D}_{\ell,j-\ell}\big(\Pi_{\perp}f\big) =O(e^{-(k-\ell+\alpha)t})
 \end{equation*}
as soon as $\mathcal{D}_{\ell,j-\ell}$ denotes a product $(j-\ell)$ factors of which are equal to $e^t\partial_{\theta}$, and $\ell$ factors are among $\{r|\log r|\partial_r,\partial_{z_\beta},\dbar_{z_\beta}, \,\beta\geq2\}$, where $r=|z_1|$.

Having said this, we come to the promised proof. 
Now according to \S \ref{poinc_ineq} to \ref{ctrl_potl_grth}, we know that $v\in H^{k+2}$ and $v\in tC^{k,\alpha}$ (we know that $v=O(t)$; we get that $v\in tC^{k,\alpha}$ by Schauder estimates in a system of quasi-coordinates).

To see that $v\in \mathcal{E}^{k+2,\alpha}$, we consider the Dirichlet problem:
 \begin{equation*}
   \left\{\begin{aligned}
            \Delta_{\omega} w  &= g \quad \text{in }\mathcal{N}_A\backslash D, \\
                            w  &= 0 \quad \text{on }\partial\mathcal{N}_A=\{t=A\}
          \end{aligned}
   \right.
 \end{equation*}
with $g\in C^{k,\alpha}(\mathcal{N}_A\backslash D)$ and $w\in tC^{k+2,\alpha}(\mathcal{N}_A\backslash D)$ (obtained by exhaustion). 
Indeed if $\gamma$ is a smooth cut-off function equal to $1$ on $\{t\leq A\}$ and vanishing on $\mathcal{N}_{A+1}$, we get $v$ as $v_{int}+v_{ext}$, with 
 \begin{equation*}
  \left\{\begin{aligned}
          \Delta_{\omega}v_{int}     &= (\Delta_{\omega}\gamma)v+\gamma f-2(d\gamma,dv)_{\omega} \quad \text{in }X\backslash\mathcal{N}_{A+1}, \\
                         v_{int}    &= 0 \quad \text{on }\partial\mathcal{N}_{A+1}=\{t=A+1\}
         \end{aligned}
  \right. 
 \end{equation*}
and
 \begin{equation*}
  \left\{\begin{aligned}
          \Delta_{\omega}v_{ext}     &= (-\Delta_{\omega}\gamma)v+(1-\gamma) f+2(d\gamma,dv)_{\omega} \quad \text{in }\mathcal{N}_A\backslash D, \\
                         v_{ext}    &= 0 \quad \text{on }\partial\mathcal{N}_A=\{t=A\}
         \end{aligned}
  \right. 
 \end{equation*}
(the right-hand-side members being controlled by $\|f\|_{C^{k,\alpha}}$, since for all $C$ there exists $K=K(C)$ such that $\|v\|_{C^{k+2,\alpha}(\mathcal{N}_C)}\leq K\|f\|_{C^{k,\alpha}}$), and $v_{int}$, $v_{ext}$ extended by 0. 
The role of $v_{ext}$ will be played $w$, whereas $g$ will play that of $(-\Delta_{\omega}\gamma)v+(1-\gamma) f+2(d\gamma,dv)_{\omega}$.
 
We introduce the subspace 
 \begin{align*}
  \mathcal{F}^{k+2,\alpha}=\big\{v\in C^{k,\alpha}_{loc}(\mathcal{N}_A\backslash D)|\, v_0=O(t),&\,\partial_t v_0\in C^{k+1,\alpha}(\mathcal{N}_A\backslash D)\,; \\
                                                           &\,v_1,\,\Pi_{\perp}v\in C^{k+2,\alpha}(\mathcal{N}_A\backslash D)\,;\,v|_{t=A}\equiv0 \big\},
 \end{align*}
of $tC^{k,\alpha}$ endowed with the obvious norm, and we assume $A$ big enough so that $\Delta_h-\Delta_{\omega}: \mathcal{F}^{k+2,\alpha}\to C^{k,\alpha}(\mathcal{N}_A\backslash D)$ has a sufficiently small norm, 
where $h$ is the metric 
 \begin{equation*} \label{eq_formule_h}
  h=dt^2+e^{-2t}\eta^2+p^*g_D.
 \end{equation*}
(compare with the asymptotics \eqref{eq_asmpt_g}). 
If one shows that $\Delta_h : \mathcal{F}^{k+2,\alpha}\to C^{k,\alpha}(\mathcal{N}_A\backslash D)$ is invertible of inverse $G_h$, with $\|G_h\|$ remaining bounded if $A$ increases, a perturbation argument will tell us that $\Delta_{\omega}$ is also invertible ; 
one writes $\Delta_{\omega}=\Delta_h\big(1-G_h(\Delta_h-\Delta_{\omega})\big)$. 
In the final analysis, there remains to see that the solution of 
 \begin{equation*}
   \left\{\begin{aligned}
            \Delta_{h} w      &= g \quad \text{in }\mathcal{N}_{A}\backslash D, \\
                           w  &= 0 \quad \text{on }\partial\mathcal{N}_A=\{t=A\}
          \end{aligned}
   \right.
 \end{equation*}
which is in $tC^{k+2,\alpha}(\mathcal{N}_A\backslash D)$ is in $\mathcal{F}^{k+2,\alpha}$. 
Observe that $\Delta_h$ respects decomposition \eqref{eqdecomp1}, hence $\Delta_h w_0=g_0$, $\Delta_h w_1=g_1$ and $\Delta_h (\Pi_{\perp}w_1)=\Pi_{\perp}g$. 
We also show that the component $w_0$ has bounded derivative, and that the other two are bounded and with bounded derivatives:
 \begin{itemize}
  \item[$\bullet$]$w_0$: the condition $v_0(A)=0$, as well as the identity $\Delta_h w_0=-(\partial_t^2-\partial_t)w_0$ give 
   \begin{equation*}
    \partial_tw_0(t)= e^t\int_t^{+\infty} e^{-s}g_0(s)\,ds= O(1) 
   \end{equation*}
and
 \begin{equation*}
    w_0(t)= \int_A^t e^sds\int_s^{+\infty} e^{-u}g_0(u)\,du =O(t).
   \end{equation*}
Those formulas clearly give the norms of $w_0$ and $\partial_t w_0$ are controlled by $\|g\|_{C^{k,\alpha}}$, independently of $A$.
  \item[$\bullet$]$w_1$: set $a(t)=\int_D w_1(t,\cdot)^2\vol^{g_D}$; if one shows that $a(t)$ is bounded, then the classical theory will tell use that  $w_1$ is bound (with an effective bound coming from that of $a(t)$). 
Now $\partial_ta(t)=2\int_D w_1(t,\cdot)\partial_tw_1(t,\cdot)\vol^{g_D}$ and $\partial_t^2a(t)=2\big(\int_D w_1(t,\cdot)\partial_t^2w_1(t,\cdot)\vol^{g_D}+\int_D \big(\partial_t^2w_1(t,\cdot)\big)^2\vol^{g_D}\big)$. 
In this way:
 \begin{align*}
  (\partial_t^2-\partial_t)a(t) & \geq 2\int_D w_1(t,\cdot)\big(\partial_t^2w_1(t,\cdot)-\partial_tw_1(t,\cdot)\big)\vol^{g_D} \\
                                &  =   2\int_D w_1(t,\cdot)\big(\Delta_Dw_1(t,\cdot)-g_1(t,\cdot)\big)\vol^{g_D}               \\
                                &  =   2\int_D \big|d_Dw_1(t,\cdot)\big|^2_{g_D}\vol^{g_D}- 2\int_D w_1(t,\cdot)g_1(t,\cdot)\vol^{g_D} \\
                                & \geq c a(t)-C(g)a(t)^{1/2},
 \end{align*} 
where we go from the first to the second line by noticing that $\Delta_hv_1=-(\partial_t^2-\partial_t)w_1+\Delta_{g_D}w_1$, with $c$ coming from Poincaré inequality for $g_D$ (one has $\int_D w_1(t,\cdot)\vol^{g_D}=0$), and $C(g)$ is the supremum of $\big(\int_D g_1(t,\cdot)^2\vol^{g_D}\big)^{1/2}$. 
According to Lemma \ref{lem_appdx} following this proof, this inequality forces $a$ to be bounded, and $a(t)\leq (\tfrac{C(g)}{c})^2$. 
In other terms, the $L^2$ norm of $w_1$ on each $\{t\}\times D$ remains bounded, and is smaller than $\tfrac{C(g)}{c}\leq C' \|f\|_{C^{k,\alpha/2}}$, $C'$ independent of $A$, hence an analogous estimation on $\|w_1\|_{C^1}$. 
 \item[$\bullet$]$\Pi_{\perp}w$ : as $\Pi_{\perp}g\in e^{-(\alpha/2)t}C^{k,\alpha/2}$, according to the weighted analysis in \cite{biq}, $\Pi_{\perp}v$ is in $e^{-\beta t}C^{k+2,\alpha/2}$ for some $\beta>0$, and in particular is bounded, as well as its differential. 
We also have that $C^1$ is controlled by $\|g\|_{C^{k,\alpha}}$, independently of $A$.
 \end{itemize}

The classical elliptic theory gives us that $w\in \mathcal{F}^{k+2,\alpha}$, with $\|w\|_{\mathcal{F}^{k+2,\alpha}}\leq C\|g\|_{C^{k,\alpha}(\mathcal{N}_A\backslash D)}$, $C$ independent of $A$ 
(one applies Schauder estimates on balls $B$ of quasi-coordinates to $w$ to which is subtracted its mean on $B$; this gives a family uniformly bounded in $C^0$, since $w$ has bounded derivatives), which ends the proof. \hfill $\square$

~

We close this part with the statement and the proof of the lemma used in the previous proof:
 \begin{lem} \label{lem_appdx}
  Let $b$ a nonnegative $C^2_{loc}$ function on $[A,+\infty[$, vanishing at $A$. 
We assume that $b = O(t^{\beta})$ for some $\beta > 0$, that $b$, $\partial_tb$ and $\partial^2_tb$ are $L^1$ for $e^{-t}dt$, and that 
 \begin{equation} \label{ineq_lem_appdx}
  (\partial_t^2-\partial_t-c)b\geq -Cb^{1/2},
 \end{equation}
with $c>0$, $C\geq0$. 
Then $b$ is bounded above, and $\sup b\leq \big(\tfrac{C}{c}\big)^2$. 
 \end{lem}
\prf. Assume that $b$ is not identically 0, and that $\beta<1$, so that $b=o(t)$. 
Then $b_{\vareps}:t\mapsto b-\vareps(t-A)$ goes to $-\infty$ after reaching its upper bound at a point $t_{\vareps}\in]A,+\infty[$, and this for all $\vareps>0$. 
At such a point, $\partial_t^2b(t_{\vareps})=\partial_t^2b_{\vareps}(t_{\vareps})\leq 0$ and $\partial_t^2b(t_{\vareps})=\partial_t^2b_{\vareps}(t_{\vareps})+\vareps= \vareps$. 
From \eqref{ineq_lem_appdx}, we hence have that $cb(t_{\vareps})\leq Cb(t_{\vareps})^{1/2}-\vareps\leq Cb(t_{\vareps})^{1/2}$, that is $b(t_{\vareps})\leq (\tfrac{C}{c})^2$

Now, at fixed $t$, $b(t)=\lim_{\vareps\to 0}b_{\vareps}(t)$, and for all $\vareps>0$, $b_{\vareps}(t)\leq b_{\vareps}(t_{\vareps})\leq b(t_{\vareps})\leq (\tfrac{C}{c})^2$, d'où $b(t)\leq (\tfrac{C}{c})^2$. 
This holding for all $t$, we have that $b$ is bounded above, with the announced bound. 

There remains to see that we can take $\beta<1$. 
Set $B(t)= -(\partial_t^2-\partial_t-c)b(t)$ on $[A,+\infty[$. 
This can be integrated into
 \begin{equation*}
  b= e^{\nu t} \int_A^t   e^{(\mu-\nu) s}ds \int_s^{+\infty} e^{-\mu u} B(u)\, du,
 \end{equation*}
with $\mu>\nu$ the roots of $X^2-X-c$ ($\mu>1$, $\nu<0$). 
Now by \eqref{ineq_lem_appdx}, $B(t)\leq C't^{\beta/2}$, hence $\int_t^{+\infty} e^{-\mu u} B(u)\leq C't^{\beta/2}e^{-\mu t}$, \textit{etc.}, hence $b=O(t^{\beta/2})$ (since $b\geq0$). 
We concludes by an immediate induction.  \hfill $\square$

~

The remaining question is the following:

~

\noindent\textit{Question.} 
Does Proposition \ref{prop_roughpnc=pnc2} hold when $D$ has simple normal crossings ?

~

Actually, using integral formulas like (\ref{integral_formula_1}) below and the fact that components orthogonal to the constants on the $S^1$-fibers around the divisor have a harmless behaviour, there is not much difficulty seeing that the differential of such a potential has its component in the normal directions to the divisor bounded. 
However it seems delicate to adapt our proof of Proposition \ref{prop_roughpnc=pnc} in the normal crossing case.

 \section[Resolution of homogeneous Monge-Ampère equation]{Resolution of the homogeneous Monge-Ampère equation on the product $(\XD)\times\Sigma$} \label{resolution_mah}

  \subsection{The theorem and its interpretation in terms of geodesics}

The result we get in the present part is:
 \begin{thm}  \label{thm_MA_hmg}
  Equation \eqref{MA_hmg} with boundary conditions and $S^1$ invariance admits a solution in the sense of currents. 
More precisely, this solution is the increasing limit of $C^{\infty}\big((\XD)\times\Sigma\big)$ and $S^1$-invariant deformations $\Phi_r$ of the segment $\Xi:=\big((1-t)v_0+tv_1\big)_{t\in[0,1]}$, satisfying the equations 
  \begin{equation}  \label{MA_htrg}
   \big({\pr_{\XD}}^*\omega_{bp}+i\ddbar\Phi_r\big)^{m+1}\equiv cr\frac{i}{2}dw\wedge d\overline{w}\wedge\big({\pr_{\XD}}^*\omega_{bp})^m
  \end{equation}
for arbitrarily small $r>0$, where $dw=dt+ids$ and $c>0$ is a positive constant, and with ${\pr_{\XD}}^*\omega+i\ddbar\Phi_r$ positive and $C^{\infty}$-quasi-isometric to $\tfrac{i}{2}dw\wedge d\overline{w}+{\pr_{\XD}}^*\omega$. 
Finally, there exist uniform $C^0$ and $C^1$ bounds on $\Phi_r-\Xi$, as well as uniform bounds on $i\ddbar(\Phi_r-\Xi)$.
 \end{thm}

The proof, which follows Chen's \cite{chen1} of the compact case, itself in the line of works like \cite{ckns, guan}, consists in a \textit{continuity method} which requires several steps. The method is explained in next section, the estimates we need to achieve it are obtained in sections \ref{uniq_sol}, \ref{higher_order_est} and \ref{C2_and_C2alpha_est}, and proof is completed in section \ref{ccl_cont_meth}.

For now, we shall translate Theorem \ref{thm_MA_hmg} into the language of paths in $\tmom$ between $v_0$ and $v_1$, since this is what we need to show Theorem \ref{thm_unqnss_csck_met} of part \ref{unqnss_cscK}:
 \begin{crl}   \label{crl_eps_gdscs}
  For any $v_0$, $v_1\in\tmom$ and any small enough $\vareps>0$ there exists a path $(v^{\vareps}_t)$ from $v_0$ to $v_1$ which is a $C^{\infty}$ deformation of the segment $\big((1-t)v_0+tv_1\big)$, satisfying the equation 
$\big(\ddot{v^{\vareps}_t}-\big|\partial\dot{v^{\vareps}_t}\big|_{\omega^{\vareps}_t}^2\big)(\omega^{\vareps}_t)^m=\vareps\omega_{bp}^m$, where $\omega^{\vareps}_t=\omega_{bp}+i\ddbar v^{\vareps}_t$. 
There exists $C>0$ such that for all $\vareps$, 
$\big|v^{\vareps}_t-\big((1-t)v_0+tv_1\big)\big|$, $|dv^{\vareps}_t|_{\omega_{bp}}$, $\big|\ddot{v^{\vareps}_t}\big|$, $\big|d\dot{v^{\vareps}_t}\big|_{\omega_{bp}}$, $\big|i\ddbar v^{\vareps}_t\big|_{\omega_{bp}}\leq C$ where $d$, $\partial$ and $\dbar$ are those of $\XD$ and $\dot{}$ stands for $\partial_t$. 
 \end{crl} 
\prf. Take $\vareps>0$ small, and for all $t\in[0,1]$ denote by $v^{\vareps}_t$ the function $\Phi_{\vareps}(\cdot,t,\cdot)$, with $\Phi_{\vareps}$ that of Theorem \ref{thm_MA_hmg} (with $\vareps$ instead of $r$); this makes sense, since every summand is $S^1$-invariant. 
Moreover, $v^{\vareps}_{\tau}=v_{\tau}$, $\tau=0,1$, since $(\Phi_{\vareps}-\Xi)|_{(\XD)\times\partial\Sigma}\equiv0$, and $(v^{\vareps})_{t\in[0,1]}\in \mathcal{E}\big((\XD)\times[0,1]\big)$. 
To assert that $(v^{\vareps}_t)$ is a \textit{path} from $v_0$ to $v_1$, we thus only have to check that $\omega^{\vareps}_t=\omega_{bp}+i\ddbar v^{\vareps}_t$ is quasi-isometric to $\omega$ for all $t\in[0,1]$, where $\partial$ and $\dbar$ are those of $\XD$; this simply follows from the fact that for all $t\in[0,1]$, $\omega^{\vareps}_t$ is the restriction of ${\pr_{\XD}}^*\omega_{bp}+i\ddbar\phi_{\vareps}$ ($\partial$ and $\dbar$ of $(\XD)\times\Sigma$) to the subbundle $\Lambda^{1,1}_{(\XD)\times\{t\}}$ of $\Lambda^{1,1}_{(\XD)\times\Sigma}$, 
and from the mutual bound between ${\pr_{\XD}}^*\omega_{bp}+i\ddbar\phi_{\vareps}$ and ${\pr_{\XD}}^*\omega_{bp}+\frac{i}{2}dw\wedge d\overline{w}$ required in \eqref{MA_htrg}.

We furthermore have from Theorem \ref{thm_MA_hmg} a bound on $i\ddbar\Phi_{\vareps}$, \textit{independent of }$\vareps$, $\partial$ and $\dbar$ being those of the product $(\XD)\times\Sigma$. 
This tells us that there is some $C$ such that for all small $\vareps>0$, 
 \begin{equation} \label{unif_C11_ctrol}
  \big|\ddot{v^{\vareps}_t}\big|, \big|d\dot{v^{\vareps}_t}\big|_{\omega_{bp}}, \big|i\ddbar v^{\vareps}_t\big|_{\omega_{bp}}\leq C
 \end{equation}
(the linear part $(1-t)v_0+tv_1$ of $(v^{\vareps}_t)$ is killed by $\partial_t^2$, and $dv_0$, $dv_1$, $i\ddbar v_0$ and $i\ddbar v_1$ are bounded).

Finally, expressing \eqref{MA_htrg} 
on $\XD$ and forgetting the pull-backs, we have: 
 \begin{equation*}
  \big(\ddot{v^{\vareps}_t}-\big|\partial\dot{v^{\vareps}_t}\big|_{\omega^{\vareps}_t}^2\big)(\omega^{\vareps}_t)^m\wedge\frac{i}{2}dw\wedge d\overline{w}=
                            4c\vareps\omega_{bp}^m\wedge\frac{i}{2}dw\wedge d\overline{w};
 \end{equation*}
since we focus on small $\vareps$, we can assume up to rescaling that $4c=1$, and hence $\big(\ddot{v^{\vareps}_t}-\big|\partial\dot{v^{\vareps}_t}\big|_{\omega^{\vareps}_t}^2\big)(\omega^{\vareps}_t)^m=\vareps\omega_{bp}^m$. \hfill $\square$

\begin{df} For $\vareps>0$, a path as in Corollary \ref{crl_eps_gdscs} is called an $\vareps$-geodesic between $v_0$ and $v_1$.
\end{df}

  \subsection{The continuity method} \label{cont_met}

Observe that equation \eqref{MA_hmg} can be rewritten as 
 \begin{equation*} 
  \big(({\pr_{\XD}}^*\omega_{bp}+\frac{i}{2}dw\wedge d\overline{w}) +i\ddbar(\Phi+t(1-t))\big)^{m+1}\equiv 0,
 \end{equation*}
since $\tfrac{i}{2}dw\wedge d\overline{w}=dt\wedge ds=-i\ddbar\big(t(1-t)\big)$ ($w$ is a local holomorphic coordinate on $\Sigma$ such that $dw=dt+ids$). 
This rewriting takes into account that $\hat{\omega}:={\pr_{\XD}}^*\omega_{bp}+\tfrac{i}{2}dw \wedge d\overline{w}$ is a Kähler form on $(\XD)\times\Sigma$, whereas ${\pr_{\XD}}^*\omega_{bp}$ is degenerate in the $\Sigma$ direction. Having said this, we finally study the equation 
\begin{equation*}
  (\omega+i\ddbar\phi)^{m+1}\equiv 0
 \end{equation*}
on $(\XD)\times \Sigma$ with $\omega$ now standing for $\hat{\omega}$, which as a product of two reference metrics will be our reference product metric on $(\XD)\times\Sigma$. 
To generalize the definitions of Hölder spaces of functions or tensors on $(\XD)\times\Sigma$ in an easy way we use quasi-coordinates, by replacing the polydiscs we used on $\XD$ by their product with (half-)balls of coordinates of homogeneous diameter forming an atlas of $\Sigma$. In this part $C^{k,\alpha}$, $C^{\infty}$, $\Gamma^{k,\alpha}(\Lambda^{1,1})$ and so on will thus denote such spaces on $(\XD)\times\Sigma$, unless otherwise specified. 
One last remark is the $S^1$-invariance of our new $\omega$, as well as this of $\phi$, if this latter stands for some $\Phi+t(1-t)$.

Now let us give ourselves an $S^1$-invariant function $\theta:[0,1]\times(\XD)\times\Sigma\rightarrow\R$ strictly increasing in $r\in[0,1]$ at every point, such that $\theta(0,\cdot)\equiv0$, $\theta(1,\cdot)\equiv1$, bounded below by $cr$ for some positive constant $c$ and with nice derivatives, namely such that $\theta$ would be in a space denoted by $C^{\infty}\big([0,1]\times(\XD)\times\Sigma\big)$. 
The \textit{continuity method} we propose consists in resolving for $r\in(0,1]$ the family of equations
 \begin{equation*} \tag{$E_r$}
  \left\{
   \begin{aligned}
         &\big(\omega+i\ddbar\phi\big)^{m+1}= \theta(r)\big(\omega+i\ddbar\phi_1\big)^{m+1} & & \\
         &\phi|_{(\XD)\times\{\tau\}\times S^1}=v_{\tau},\,\tau=0,1                         & &\text{(boundary conditions)}\\
         &c\omega \leq \omega+i\ddbar\phi \leq c^{-1}\omega                                 & &\text{for some constant }c>0\\
         &\phi-\phi_1 \in C^{\infty}\big((\XD)\times\Sigma\big).                            & &
   \end{aligned} \right.
 \end{equation*}
where $\phi_1$ is itself the solution of $(E_1)$, meaning that $\omega+i\ddbar\phi_1$ is $C^{\infty}$-quasi-isometric to $\omega$ on $(\XD)\times\Sigma$ and $\phi_1|_{(\XD)\times\{\tau\}\times S^1}=v_{\tau}$, $\tau=0,1$.

The first step is ensuring that such a $\phi_1$ exists. Actually, an easy computation provides it as a $C^{\infty}\big((\XD)\times\Sigma\big)$ deformation of the segment $\Xi$ joining to $v_0$ and $v_1$. 
Namely, if $C>0$ is a constant, and if $\phi_1:=(1-t)v_0(z)+tv_1(z)-Ct(1-t)$ (notice it is $S^1$-invariant), one has:
 \begin{equation*}
  \omega+i\ddbar\phi_1=(1-t)\omega_{v_0}+t\omega_{v_1}+2\mathfrak{Re}\big(i\partial(v_1-v_0)\wedge d\overline{w}\big)
                       +(C+1)\frac{i}{2}dw\wedge d\overline{w},
 \end{equation*}
which clearly is $C^{\infty}$-quasi-isometric to $\omega$ when $C$ is big enough, since $dv_0$ and $dv_1$ are in $\Gamma^{\infty}(\Lambda^1,\XD)$ (and in particular bounded for a Poincaré type metric).

Having settled this question, our strategy is to show the following:
 \begin{prop}  \label{prop_r0=0}
  Let $r_0$ the infimum of the $r$ such that $(E_{r'})$ admits a unique solution for all $r'\in(r,1]$. Then $r_0=0$.
 \end{prop}
The proof is done in section \ref{ccl_cont_meth}, but requires a significant preparatory work, in particular in obtaining a priori estimates for solutions of equations $(E_r)$. For now, let us deal with the uniqueness of their solutions, as well as some $C^0$ estimates.

\subsection{Uniqueness and a priori $C^0$ estimates of intermediate solutions} \label{uniq_sol}

 \begin{prop}[Uniqueness and $C^0$ estimate] \label{prop_unqss_C0_estimate}
  For any $r\in(0,1]$, the solution $\phi$ of $(E_r)$ is unique if exists; in particular it is $S^1$-invariant. Moreover, $\phi_1\leq\phi\leq\phi_1+h$ for some bounded function $h\in C^{\infty}$ vanishing on $(\XD)\times\partial\Sigma$, and if $\phi'$ is the solution of $(E_{r'})$, $r'\in(0,1]$, then $r'\leq r$ implies $\phi\leq\phi'$, and reverse.
 \end{prop}
\prf. The idea underlying the technique used here consists in making apparent some functions sub/over-harmonic with respect to well-chosen metrics and which vanish on $(\XD)\times \partial\Sigma$ 
and then apply an appropriate maximum principle (Lemma \ref{lem_strong_max_principle}). This latter states, in a weak form:

 \begin{lem}  \label{lem_max_principle_w/_bdry}
 Let $v$ be a $C^2_{loc}$ function bounded above on $(\XD)\times\Sigma$, such that $\sup_{(\XD)\times\Sigma}v>\sup_{(\XD)\times\partial\Sigma}v$. 
Assume $(\XD)\times\Sigma$ is endowed with a Kähler metric $\omega'$ quasi-isometric to $\omega$. Then there exists a sequence $(x_j)_{j\geq0}$ of points of $(\XD)\times\mathring{\Sigma}$ such that
  \begin{align*}
   \lim_{j\rightarrow \infty}v(x_j)             =\sup_{(\XD)\times\Sigma}v,\,\,\,
   \lim_{j\rightarrow \infty}|dv(x_j)|          =0,\text{ and }\,\,\,
   \liminf_{j\rightarrow \infty} \Delta_{\omega'} v(x_j)\geq 0. 
  \end{align*}
 \end{lem}
\textit{Proof of Lemma }\ref{lem_max_principle_w/_bdry}. It is very similar to that of Wu's maximum principle \cite{wu} p.406, but adapted to the boundary context. We do it for $\omega$. 
For $\vareps>0$ set $v_{\vareps}=v-\vareps{\pr_{\XD}}^*\mathfrak{u}$, so that it goes to $-\infty$ near $D\times\Sigma$. Suppose that for all $\vareps>0$, $v_{\vareps}$ reaches its maximum at some $x_{\vareps}\in(\XD)\times\partial \Sigma$. 
It is then not hard to see that for all fixed $x\in(\XD)\times\Sigma$, $\liminf_{\vareps\rightarrow0}v(x_\vareps)\geq v(x)$, hence a contradiction with the assumption $\sup_{(\XD)\times\Sigma}v>\sup_{(\XD)\times\partial\Sigma}v$ since of course $\liminf_{\vareps\rightarrow0}v(x_{\vareps})\leq \sup_{(\XD)\times\partial\Sigma}v$.

Having said this, we know there is an $\vareps_0>0$ such that $v_{\vareps_0}$ raises its maximum at some $x_{\vareps_0}\in(\XD)\times\mathring{\Sigma}$. 
Applying the reasoning above to the $\vareps\in(0,\vareps_0)$ gives an $\vareps_1$ such that $v_{\vareps_1}$ raises its maximum at some $x_{\vareps_1}\in(\XD)\times\mathring{\Sigma}$, and so on. 
Set $x_j=x_{\vareps_j}$; a glance Wu's proof shows $(x_j)$ verifies the stated assertions. \hfill $\blacksquare$

~

The following will be useful to strengthen our maximum principle:
 \begin{lem}  \label{lem_alpha_pot}
  Assume $(\XD)\times\Sigma$ is endowed with a Kähler metric $\omega'$ $C^{\infty}$-quasi-isometric to $\omega$. Then there exists a $C^{\infty}\big((\XD)\times \Sigma\big)$ solution $\alpha$ to the Dirichlet problem
  \begin{equation*}
   \left\{
    \begin{aligned}
     &\Delta_{\omega'}\alpha =1\\
     &\alpha|_{(\XD)\times\partial\Sigma}=0.
    \end{aligned}
   \right.
  \end{equation*}
  Moreover, $0\leq\alpha\leq Ct(1-t)$ for some $C=C(\omega')>0$.
 \end{lem}
\textit{Proof of Lemma }\ref{lem_alpha_pot}. This follows from an exhaustion argument; namely, $\alpha$ is obtained as the $C^2$-limit on each compact subset of $(\XD)\times\Sigma$ of some subsequence of a sequence of $C^{2,\beta}$ solutions $(\alpha_p)$ ($\beta\in (0,1)$) of the analogous Dirichlet problem on an exhaustive sequence of compact subdomains $(V_p)$ of $(\XD)\times\Sigma$. 
From the uniform ellipticity of $\Delta_{\omega'}$ on (half-)balls of quasi-coordinates, it suffices to have a uniform $C^0$ control on the $\alpha_p$ to get a uniform $C^{2,\beta}$ control and then perform some extraction. 
The nonnegativity of the $\alpha_p$, hence that of $\alpha$, is clear. 
Moreover $\Delta_{\omega'}\big(t(1-t)\big)=2\tr^{\omega'}(idt\wedge ds)=|dt|^2_{\omega'}\geq c|dt|_{\omega}^2=c$ where $c=c(\omega')>0$ is such that $\omega'\geq c\omega$. 
Finally, $\Delta_{\omega'}\big(c^{-1}t(1-t)-\alpha_p\big)\geq0$ on $V_p$ and $\big(c^{-1}t(1-t)-\alpha_p\big)\geq0$ on $\partial V_p$ so $c^{-1}t(1-t)-\alpha_p\geq0$ i.e. $\alpha_p\leq c^{-1}t(1-t)$ for all $p$, so we are done for the sought $C^0$ estimate. 
This estimate passes to $\alpha$, which then is $C^{\infty}\big((\XD)\times \Sigma\big)$ still by uniform ellipticity of $\Delta_{\omega'}$ on (half-)balls of quasi-coordinates and since $1\in C^{\infty}\big((\XD)\times \Sigma\big)$. \hfill $\blacksquare$

 \begin{rmk}  \label{rmk_isom1}
  Similar arguments give isomorphisms $\Delta_{\omega'}:C^{k+2,\beta}_0\rightarrow C^{k,\beta}$ for every $(k,\beta)\in \N\times(0,1)$ for any $\omega'$ quasi-isometric to $\omega$, the 0 index meaning "vanishing on $(\XD)\times\partial\Sigma$".
 \end{rmk}

Combining the last two lemmas, one gets:
 \begin{lem}  \label{lem_strong_max_principle}
  Let $v$ be a $C^2_{loc}$ function bounded above on $(\XD)\times\Sigma$, nonpositive on $(\XD)\times\partial\Sigma$. Assume $(\XD)\times\Sigma$ is endowed with a Kähler metric $\omega'$ $C^{\infty}$-quasi-isometric to $\omega$ and that $\Delta_{\omega'}v\leq0$. Then $v\leq0$.
 \end{lem}
\textit{Proof of Lemma }\ref{lem_strong_max_principle}. Suppose there exists a point $x\in(\XD)\times\Sigma$ such that $v(x)>0$. 
Then for some $\vareps>0$ small enough $(v-\vareps\alpha)(x)>0$ hence $\sup_{(\XD)\times\Sigma}(v-\vareps\alpha)>0\geq\sup_{(\XD)\times\partial\Sigma}$. 
Take a sequence $(x_j)$ as in Lemma \ref{lem_max_principle_w/_bdry} for $v-\vareps\alpha$; in particular, $\liminf_{j\rightarrow\infty}\Delta_{\omega'}(v-\vareps\alpha)(x_j)\geq0$, whereas this is equal to $-\vareps+\liminf_{j\rightarrow\infty}\Delta_{\omega'}v(x_j)\leq-\vareps$, hence a contradiction. \hfill $\blacksquare$

~

Let us come back to the proof of Proposition \ref{prop_unqss_C0_estimate}. Denote by $\psi$ the difference $\phi-\phi_1$, so that $\psi\in C^{\infty}_0$.
We claim that $\psi$ is over-harmonic with respect to $\omega'=\omega+i\ddbar\phi$ (which is $C^{\infty}$-quasi-isometric to $\omega$); this can be seen at any point $x$ by taking coordinates $(z_1,\dots,z_{m+1})$ such that
$\omega'=\sum_j idz_j\wedge d\overline{z_j}$ and $i\ddbar\psi=\sum_j \lambda_j idz_j\wedge d\overline{z_j}$ at $x$. 
From logarithm concavity we write $\big(\frac{1}{\theta(r)}\big)^{\tfrac{1}{m+1}}= \prod_{j=1}^{m+1}(1-\lambda_j)^\frac{1}{m+1} \leq \tfrac{\sum_{j=1}^{m+1}(1-\lambda_j)}{m+1} 
= 1+\frac{\tfrac{1}{2}\Delta'\psi}{m+1}$ i.e. $\Delta'\psi\geq2(m+1)\big(\theta(r)^{-1/(m+1)}-1\big)\geq 0$ where $\Delta'$ is the Laplacian associated to $\omega'$. 
From the latter proposition above, this gives: $\psi\geq0$, i.e. $\phi\geq\phi_1$.
Using the same techniques we show that:
 \begin{itemize}  \setlength{\itemsep}{0pt}
  \item[$\bullet$] If $\phi'$ denotes a solution of $(E_{r'})$, $r\leq r'\leq1$, then $\Delta'(\phi-\phi')\geq0$ so $\phi\geq\phi'$. Reverse inequality comes from symmetry. This provides the uniqueness, and hence the $S^1$-invariance, statements.
  \item[$\bullet$] Keep the notation $\psi=\phi-\phi_1$. If $h$ denotes the $C^{\infty}_0$ function such that $\Delta_1h=2(m+1)$ given by Lemma \ref{lem_alpha_pot} --- 
then $\psi\leq h$, i.e. $\phi\leq \phi_1+h$. This comes from the inequality $\omega=\omega_1+i\ddbar\psi\geq0$; taking its trace with respect to $\omega_1$ provides: $m+1-\tfrac{1}{2}\Delta_{\omega_1}\psi\geq0$. \hfill $\square$
 \end{itemize}

\subsection{Second order estimates} \label{higher_order_est}
Let us denote by $f$ the function $\tfrac{\omega_1^{m+1}}{\omega^{m+1}}$ so that $f\in C^{\infty}\big((\XD)\times\Sigma\big)$, $f\geq c$ for some positive constant $c>0$; notice that $(E_r)$, $r\in(0,1]$, sums up as $\big(\omega+i\ddbar\phi\big)^{m+1}=\theta(r)f\omega^{m+1}$, in addition to mutual boundedness and boundary conditions.
We still have some freedom on the definition of $\theta$; for instance we can take $\theta(r)=r\big((1-\chi(r))cf^{-1}+\chi(r)\big)$ where $\chi$ is an increasing smooth function on $[0,1]$ equal to $0$ (resp. 1) in a neighbourhood of $0$ (resp. $1$) and $c=\inf_{(\XD)\times\Sigma}f >0$. 
This way,  $\theta(r)f=cr$ when $r$ is close to 0.

On the other hand, since $\psi=\phi-\phi_1$ is the function that has a chance to be bounded (in general, $\phi_1$ is not if it is constructed from a segment joining unbounded potentials, so neither is $\phi$), it is convenient to look at our equations in the form: $\big(\omega_1+i\ddbar\psi\big)^{m+1}=\theta(r)\omega_1^{m+1}$. 
Let us call this latter $(E'_r)$ after adding to it mutual boundedness ($c^{-1}\omega_1\leq \omega_1\leq c^{-1}\omega_1$ for some $c>0$) and boundary ($\psi|_{(\XD)\times\Sigma}\equiv0$) conditions. Then:
 \begin{prop} \label{prop_lapl_bdry}
  There exists a constant $C$ independent of $r\in(0,1]$ such that for the solution $\psi$ of any $(E'_r)$, 
   \begin{equation}  \label{ineq_lapl_bdry}
     \sup_{(\XD)\times\Sigma}\big(m+1-\tfrac{1}{2}\Delta_1 \psi\big)\leq C\big(1+\sup_{(\XD)\times\partial\Sigma}(m+1-\tfrac{1}{2}\Delta_1 \psi)\big).
   \end{equation}
 \end{prop}
\prf. It uses an inequality due to Yau \cite{yau1}, whose proof is purely local, and writes within our setting:
 \begin{lem}[Yau] \label{lem_yau1}
  If $\inf_{j\neq l}\riem^{\omega_1}_{j\jbar l\bar{l}}$ means the infimum on $(\XD)\times\Sigma$ of the quantities $\big(\riem^{\omega_1}(\tfrac{\partial }{\partial z_j}, \tfrac{\partial }{\partial \overline{z_l}})\tfrac{\partial }{\partial z_j},\tfrac{\partial }{\partial \overline{z_l}} \big)_{\omega_1}$, $j\neq l$ where the $\tfrac{\partial }{\partial z_k}$ are taken $\omega_1$-orthonormal at the point of computation and $\Delta'$ is the Laplacian associated to $\omega_1+i\ddbar\psi$ then:
  \begin{align*} 
  -\Delta'\big(e^{-\kappa\psi}&(m+1-\tfrac{1}{2}\Delta_1 \psi)\big)\geq \\
                            &e^{-\kappa\psi}\big(\Delta_{1} \log\big(c+\chi(r)(f-c)\big)-(m+1)^2\inf_{j\neq l}\riem^{\omega_1}_{j\jbar l\bar{l}}\big)\\
                            & - \kappa e^{-\kappa\psi}(m+1)\big(m+1-\tfrac{1}{2}\Delta_1 \psi\big)                                                   \\ 
                            &+\big(\kappa+\inf_{j\neq l}\riem^{\omega_1}_{j\jbar l\bar{l}}\big)e^{-\kappa\psi}\big(m+1-\tfrac{1}{2}\Delta_{1} \psi\big)^{1+1/m}\cdot   
                              \big(r(c+\chi(r)(f-c))\big)^{-1} 
  \end{align*}
  where the only constraint on the constant $\kappa$ is: $\kappa+\inf_{j\neq l}\riem^{\omega_1}_{j\jbar l\bar{l}}>1$. In particular $\kappa$ can be chosen independently of $r$.
 \end{lem}
Fix the constant $\kappa$ of the lemma once for all. 
Now it is easy to find $K_0$ independent of $r$ such that 
\begin{align*}
1\leq & e^{-\kappa\psi}\big(\Delta_{1} \log\big(c+\chi(r)(f-c)\big)
        -(m+1)^2\inf_{j\neq l}\riem^{\omega_1}_{j\jbar l\bar{l}}\big)\\
      & -\kappa (m+1)K +\big(\kappa+\inf_{j\neq l}\riem^{\omega_1}_{j\jbar l\bar{l}}\big)K^{1+1/m}\cdot               \big(r(c+\chi(r)(f-c))\big)^{-1}
\end{align*}
as soon as $K\geq K_0$. Now, either $e^{-\kappa\psi}(m+1)\big(m+1-\tfrac{1}{2}\Delta_1 \psi\big)$ is $\leq K_0+1$ on $(\XD)\times\Sigma$ and we are done ($\psi$ is bounded), or its supremum is $>K_0$.
In this latter case suppose the supremum is not reached along $(\XD)\times\partial\Sigma$, and use Lemma \ref{lem_max_principle_w/_bdry} to get a sequence of points $(x_j)$ such that $e^{-\kappa\psi}(m+1)\big(m+1-\tfrac{1}{2}\Delta_1 \psi\big)(x_j)$ tends to our supremum, and $\Delta\big(e^{-\kappa\psi}(m+1)(m+1-\tfrac{1}{2}\Delta_1 \psi)\big)(x_j)$ to a nonnegative quantity. With our definition of $K_0$, this contradicts the formula of Lemma \ref{lem_yau1}, hence the result, since $\psi$ is bounded independently of $r$. \hfill $\square$

~

Now we can control the right-hand-side term of the inequality (\ref{ineq_lapl_bdry}) with the help of first order terms in $\psi$:
 \begin{prop} \label{prop_ctrl_lapl_diff}
 There exists a constant $C$ independent of $r\in(0,1]$ such that for the solution $\psi$ of any $(E'_r)$, 
   \begin{equation*}
   \sup_{(\XD)\times\partial \Sigma}\big(m+1-\tfrac{1}{2}\Delta_{\omega_1} \psi\big) \leq C\big(1+\sup_{(\XD)\times \Sigma}|d\psi|_{\omega_1}^2\big).
  \end{equation*}
 \end{prop}
\textit{Sketch of proof.} The proof of this Proposition is rather technical, but follows closely Chen's \cite{chen1}, Theorem 1, so we are only saying a few words about what needs to be adapted in our non-compact set-up.

First, Chen's proof works considering any point $p$ of the boundary with a half-ball $B$ of coordinates such that his reference metric is bounded above by twice of the euclidian metric, and below by one half of it on $B$. 
Moreover the radius does not depend on $p$, and the $m$ first coordinates parametrize the base whereas the last one, $z$ say parametrizes $\Sigma$; more precisely, $\Sigma$ is given by $\{\mathfrak{Re}(z)\geq0\}$ in $B$.

Of course we cannot proceed like this with our kind of metrics (the injectivity radius goes to 0), but we already know that having uniform estimates on the pullbacks by some quasi-coordinate system (the $\Phi_{\delta}$) provides global bounds. 
So that we replace Chen's coordinate half-balls by quasi-coordinate half-balls, namely we fix a ball of radius $\delta>0$ in $\C^m\times \{\mathfrak{Re}(z)\geq0\}$ and consider a family $(\pi_p)_{p\in(\XD)\times\partial\Sigma}$ of holomorphic immersions $B\rightarrow(\XD)\times\Sigma$ such that for all $p\in(\XD)\times\partial\Sigma$, $\pi_p$ sends $0$ to $p$, $B\cap(\C^m\times\{0\})$ in $(\XD)\times\partial\Sigma$ and $\tfrac{1}{2}\omega_{euc}\leq {\pi_p}^*\omega_1\leq 2\omega_{euc}$. 
This way, we can apply Chen's techniques to ${\pi_p}^*\psi$, and get analogous results, in particular the fact that the normal-tangential (resp. tangential-tangential) second derivatives at $p$ are controlled by the $L^{\infty}$ norm (resp. the squared $L^{\infty}$ norm) of its differential, control which does not depend on $p$.

One subtlety though; to prove the nonnegativity of Chen's barrier function $\nu$ when $\delta$ is small enough, instead of using positive lower bounds on $\Delta'\nu$ (or ${\pi_p}^*(\Delta'\nu)$), we directly use the definition of this function, and the fact that for some constant $C$ independent of $p$, if $x$ stands for $\mathfrak{Re}(z)$, we have $0\leq x\leq C'h$. 
Indeed, we can take $\pi_p(*,z)=\big(\star,c(t+i(s-s_p))\big)$ or $\pi_p(*,z)=\big(\star,c((1-t)-i(s-s_p))\big)$, depending on which component of $(\XD)\times\partial\Sigma$ $p$ is, with $c>0$ small independent of $U$, so that we are done if we know that $t(1-t)\leq Ch$ on $\XD$ for some $C>0$. 
But such a constant exists since for $C$ big enough, $\Delta_{1}\big(Ch-t(1-t))=2C(m+1)-\Delta_{1}\big(t(1-t)\big)\geq0$ on $\XD$, and from Lemma \ref{lem_strong_max_principle}. We refer to \cite{chen1}, p.204-208, for the details. \hfill $\square$

~

Let us conclude this section with a definitive control on $\Delta_1\psi$:
 \begin{prop}   \label{prop_C11_estimate}
  There exists a constant $C$ independent of $r\in(0,1]$ such that for the solution $\psi$ of any $(E'_r)$, 
  \begin{equation*}
   \sup_{(\XD)\times \Sigma}|d\psi|_{\omega_1} \leq C.
  \end{equation*}
  In particular in view of Proposition \ref{prop_ctrl_lapl_diff}, $\sup_{\XD}\big|i\ddbar\psi\big|_{\omega_1}$ is bounded above by a constant independent of the parameter $r$.
 \end{prop}
\textit{Sketch of proof.} Here again we can adapt Chen's argument, namely his blowing-up analysis --- \cite[\S 3.2]{chen1} --- so we will not repeat it entirely here, but rather underline a few necessary changes in the proof. 

So we suppose there exists a sequence $(r_j)$ such that $\vareps_j^{-1}:=\sup_{(\XD)\times \Sigma}|d\psi_{r_j}|_{\omega_1}$ goes to $+\infty$, and we look at a sequence $(p_j)$ of points of $(\XD)\times\Sigma$ such that $|d\psi_{r_j}(p_j)|_{\omega_1}\geq\vareps_j^{-1}-1$ for all $j$.
Because in general we cannot extract from $(p_j)$ a sequence converging in $(\XD)\times\Sigma$, we \textit{follow} these points, and define objects on (half-)balls around them. 
Here nonetheless, we have to differentiate two cases: up to an extraction, $w_j:=\pr_{\Sigma}(p_j)$ converges to a point $w$ of $\Sigma$, and:
 \begin{enumerate} \setlength{\itemsep}{0pt}
  \item if $w\in\partial\Sigma$, we take $\delta>0$ small enough and give ourselves a half-disc $D_{\delta}$ of coordinate with nonnegative real part, whose radius is $\delta$, centered in $w$ and with $D_{\delta}\cap\{\mathfrak{Re}=0\}$ sent parallel to $\partial\Sigma$; 
  \item if $w\in\mathring{\Sigma}$, we take $\delta>0$ small enough and give ourselves a disc $D_{\delta}$ of coordinate in $\Sigma$,  whose radius is $\delta$ and centered in $w$.
 \end{enumerate}
In both cases, forgetting the extraction, $w_j$ is in the considered neighbourhood of $w$, and even in that of half-radius. 
Moreover we take a ball $B_{\delta}'$ of quasi-coordinate of radius $\delta$ centered in $z_j=\pr_{\XD}(p_j)$ in $\XD$, and then we have immersions
 \begin{align*}
   \pi_j : B_{\delta} &\longrightarrow  (\XD)\times \Sigma\\
           0          &\longmapsto      (z_j,w_j),
 \end{align*}
at our disposal, where $B_{\delta}$ denotes the (half-)ball of radius $\delta$ of $\C^{m+1}$ included in $B'_{\delta}\times D_{\delta}$. 
Our construction of $\omega_1$ allows us furthermore to assume ${\pi_j}^*\omega_1$ is trivial at $(0,0)$ and that its derivatives are bounded in $B_{\delta}$ independently of $j$. 
We then set, for $j$ big enough and $(z,w)\in B_{\delta/\vareps_j}$,
 \begin{equation*}
  \tilde{\psi}_j(z,w)={\pi_j}^*\psi_{r_j}\big(\vareps_j(z,w)\big),
 \end{equation*}
which defines on every compact a sequence of functions we are going to study. Similarly, for those $j$, $(z,w)$, we set 
 \begin{equation*}
  \tilde{h}_j(z,w)={\pi_j}^*h\big(\vareps_j(z,w)\big),
 \end{equation*}
and finally we set $\tilde{\omega}_j(z,w)={\pi_j}^*\omega_1\big(\vareps_j(z,w)\big)$; the previous remark ensures us that these $\tilde{\omega}_j$ converge in $C^{\infty}$ on every compact, and that we can assume $\delta$ small enough to always have $\tfrac{1}{2}\omega_{euc}\leq\tilde{\omega}_j\leq 2\omega_{euc}$.

Now this rescaling implies for all big enough $j$ that $\big|d\tilde{\psi}_j(z,w)\big|_{\tilde{\omega}_j}\leq 1$ wherever it makes sense, $\big|d\tilde{\psi}_j(0,0)\big|_{\tilde{\omega}_j}\geq 1-\vareps_j$, and $\big|\Delta_{\tilde{\omega}_j}\tilde{\psi}_j(z,w)\big|\leq C$ where $C$ is that of Proposition \ref{prop_ctrl_lapl_diff}. 
Moreover the inequalities $0\leq\psi\leq h\leq \|h\|_{C^0}$ propagate and give $0\leq\tilde{\psi}_j\leq \tilde{h}_j\leq \|h\|_{C^0}$. 
We deduce from those and from standard Schauder estimates that $\big(\tilde{\psi}_j\big)$ is $C^{1,\alpha}$ bounded on every compact as soon as it makes sense ($\alpha\in(0,1)$), hence two diagonal extractions give us a subsequence we still call $(\tilde{\psi}_j)$ which converges $C^{1,\beta}$ on every compact to some function $\tilde{\psi}$ which belongs to $C^{1,\beta}_{loc}\big(\C^m\times\{\mathfrak{Re}\geq0\}\big)$ in the first case raised above, and to $C^{1,\beta}_{loc}\big(\C^{m+1}\big)$ in the second one ($\beta\in(0,\alpha)$). In addition $\tilde{\psi}$ is bounded by $\|h\|_{C^0}$ on its whole domain, and the inequalities on the $\big|\tilde{\psi}_j(0,0)\big|_{\tilde{\omega}_j}$ tell us, passing to the limit: $\big|\tilde{\psi}(0,0)\big|_{euc}=1$.

However in the first case, it is easy to see that $\tilde{h}_j(z,w)$ tends to 0 when $j$ goes to infinity for every fixed $(z,w)$ from the very definition of the $\tilde{h}_j$. 
This implies $\tilde{\psi}\equiv0$, which contradicts $\big|\tilde{\psi}(0,0)\big|_{euc}=1$.

In the second case, using the nonnegativity of the $\omega+i\ddbar\psi_{r_j}$, we can show that on every complex line $\Pi$ passing through $0\in\C^{m+1}$, $\Delta_{\Pi}\psi\leq 0$ in the sense of distributions, hence $\tilde{\psi}$ is constant on every such $\Pi$, hence constant (it is bounded), which contradicts again $\big|\tilde{\psi}(0,0)\big|_{euc}=1$. \hfill $\square$

 \subsection{$C^2$ and $C^{2,\eta}$ estimates}  \label{C2_and_C2alpha_est}
 
We have proved a uniform (independent of the parameter $r$) estimate for the differential and the complex Hessian of our potentials $\psi$; 
notice that from $\psi|_{(\XD)\times\partial\Sigma}\equiv0$, this gives a uniform complete $C^2$ estimate of $\psi$ along $(\XD)\times\partial\Sigma$. We now give such a $C^2$ estimate on $(\XD)\times\Sigma$, which however is no more uniform, at least when $r$ goes to $0$:
 \begin{prop}  \label{prop_C2_estimate}
  Assume $\psi$ is a 
solution of some $(E'_r)$, $r\in(0,1]$. 
Then there exists some constant $C$ independent of $r$ such that $\big\|(\nabla_{\omega_1})^2\psi\big\|_{C^0}\leq\tfrac{C}{r}$.
 \end{prop}
\textit{Proof.} Here we adapt B\l ocki's proof of his Theorem 3.2 in \cite{blocki}. 
This proof uses the compactness of the underlying manifold in a crucial way, namely in working at a point where some function attains its maximum. 
Instead of making up for this lack of compactness by using, for instance, our maximum principle (Lemma \ref{lem_max_principle_w/_bdry}), we are seeing what happens when following a sequence of points such that the function in question tends to its supremum along this sequence.

To begin with, fix $r\in(0,1]$, take $\psi$ as in the statement and define a function $B$ by
 \begin{equation*}
  B: x\longmapsto \sup_{\substack{Y\in T_x ((\XD)\times\Sigma)\\ |Y|_{\omega_1}=1}}(\nabla_Yd\psi)(Y) 
 \end{equation*}
where $\nabla$ stands for $\nabla_{\omega_1}$. Notice that $B(x)$ is nothing but the biggest eigenvalue of $(\nabla_{\omega_1})^2\psi$ at $x$, so that we have to produce the desired estimate on $B$ (up to some first order term and a factor 2, this last object and $i\ddbar\psi$ have the same trace); a bound above will even be enough. 
Then define $A=B+|d\psi|_{\omega_1}$, and set $M=\sup_{(\XD)\times\Sigma}A$. 
Since $A$ is already controlled on $(\XD)\times\partial\Sigma$, we can assume that $M>\sup_{(\XD)\times\partial\Sigma}A$, and even that there exists some positive $\delta$ such that
 \begin{equation*}
  M=\sup_{(\XD)\times[\delta,1-\delta]\times S^1} A.
 \end{equation*}
Indeed, if such a $\delta$ did not exist, the $C^3$ bound we have assumed on $\psi$ would schematically provide that we can reach $M$ following a sequence of points whose projection on $[0,1]$ would tend to $0$ or $1$ and give $M=\sup_{(\XD)\times\partial\Sigma}A$.

So we have balls of quasi-coordinates $B_{\delta}\xrightarrow{\pi_j}B_j\subset (\XD)\times\Sigma$ of radius $\delta$ centered at points $O_j$ such that for all $j$:
 \begin{itemize} \setlength{\itemsep}{0pt}
  \item[a)] ${\pi_j}^*A(0)=A(O_j)\geq M-\tfrac{1}{2^j}$, and $\inf_{(\XD)\times\Sigma}A\leq{\pi_j}^*A\leq M$ ;
  \item[b)] $\tfrac{1}{2}\omega_{euc}\leq{\pi_j}^*\omega_1\leq 2\omega_{euc}$, ${\pi_j}^*\omega_1=\omega_{euc}$ at $0$, and ${\pi_j}^*\big(i\ddbar\psi\big)$ is diagonal at $0$ ;
  \item[c)] $\big({\pi_j}^*\omega_1+i\ddbar\pi_j^*\psi\big)^{m+1}={\pi_j}^*\theta(r)\cdot{\pi_j}^*(\omega_1)^{m+1}$ ;
  \item[d)]  there exists $Y_j$ of norm $1$ at $O_j$ such that ${\pi_j}^*A(0)={\pi_j}^*(\nabla_{Y_j}d\psi)(Y_j)+{\pi_j}^*(|d_{Y_j}\psi|)$.
 \end{itemize}
This way if we denote the pullbacks with hats, we have on $B_{\delta}$ for all $j$:
 \begin{itemize} \setlength{\itemsep}{0pt}
  \item[a)] $\hat{A}_j(0)\geq M-\tfrac{1}{2^j}$,and $\inf_{(\XD)\times\Sigma}A\leq\hat{A}_j\leq M$ ;
  \item[b)] $\tfrac{1}{2}\omega_{euc}\leq\hat{\omega}_j\leq 2\omega_{euc}$, $\hat{\omega}_j=\omega_{euc}$ at $0$, and $i\ddbar\hat{\psi}_j$ is diagonal at $0$ ;
  \item[c)] $\big(\hat{\omega}_j+i\ddbar\hat{\psi}_j\big)^{m+1}=\hat{\theta}_j(r)\cdot(\hat{\omega}_j)^{m+1}$ ;
  \item[d)] $\big|\hat{Y}_j\big|_{\hat{\omega}_j}=1$ (at $0$) and $\hat{A}_j(0)=\big(\nabla_{\hat{Y}_j}d\hat{\psi}_j\big)(\hat{Y}_j)+\big|d_{\hat{Y}_j}\hat{\psi}_j\big|$.
 \end{itemize}
The idea now is to let $j$ go to $\infty$ and bring the problem to the situation in which B\l ocki's proof works. 
Nonetheless we cannot assume so far that the $\hat{A}_j$ are regular, and this is why we start by some local regularizations. 
For this reason we extend the $\hat{Y}_j$ to the whole $B_{\delta}$ as constant vector fields, and for all $j$ we consider:
 \begin{equation*}
  \hat{A}'_j:=\frac{1}{|\hat{Y}_j|^2_{\hat{\omega}_j}}\big(\nabla_{\hat{Y}_j}d\hat{\psi}_j\big)(\hat{Y}_j)+\big|d\hat{\psi}_j\big|_{\hat{\omega}_j},
 \end{equation*}
so that $\hat{A}'_j\leq \hat{A}_j \leq \hat{A}_j(0)+\tfrac{1}{2^j}=\hat{A}'_j(0)+\tfrac{1}{2^j}$. 
Moreover $\hat{A}'_j$ is $C^{2,\eta}$, and bounded in $C^{2,\eta}(B_{\delta})$ independently of $j$ (thanks to similar $C^{4,\eta}$ controls on the $\hat{\psi}_j$).
On the other hand we have such $C^{k,\eta}$ controls on the $\hat{\omega}_j$ ($k=3$), $\hat{\theta}_j(r)$ ($k=4$, plus a lower bound $cr$ for this latter).
We can then simultaneously extract from our sequences weakly $C^{k,\eta}$ converging sequences, hence up to another extraction strongly $C^k$ converging sequences, with convergence to $C^{k,\eta}$ objects (and convergence in $S^{2m+1}$ for $\big(\hat{Y}_j\big)$). 
Let us simply drop the index to denote the limit; the relations above give, by passing to the limit:
  \begin{itemize} \setlength{\itemsep}{0pt}
  \item[a)] $\hat{A}'(0)= M$,and $\leq\hat{A}'\leq M$, with $\hat{A}'=\frac{1}{|\hat{Y}|^2_{\hat{\omega}}}\big(\nabla_{\hat{Y}}d\hat{\psi}\big)(\hat{Y})+\big|d\hat{\psi}\big|_{\hat{\omega}}$ ;
  \item[b)] $\tfrac{1}{2}\omega_{euc}\leq\hat{\omega}\leq 2\omega_{euc}$, $\hat{\omega}=\omega_{euc}$ at $0$, and $i\ddbar\hat{\psi}$ is diagonal at $0$ ;
  \item[c)] $\big(\hat{\omega}+i\ddbar\hat{\psi}\big)^{m+1}=\hat{\theta}(r)\cdot(\hat{\omega})^{m+1}$, $cr\leq\hat{\theta}(r)$ and control on the derivatives of $\hat{\theta}(r)$ up to order $k-1$ independent of $r$, and the same for $\tilde{\omega}$ ;
  \item[d)] $\big|\hat{Y}\big|_{\hat{\omega}}=1$ (at $0$) and $\hat{A}(0)=\big(\nabla_{\hat{Y}}d\hat{\psi}_j\big)(\hat{Y})+\big|d_{\hat{Y}_j}\hat{\psi}\big|$.
 \end{itemize}

Now we can use normal coordinates at 0 and apply B\l ocki's proof, since we have enough regularity on our objects, to get at $0$
 \begin{equation*}
  \Delta_{\hat{\omega}+i\ddbar\hat{\psi}}\hat{A}'\leq   
  -K\Big(\frac{\big(\nabla_{\hat{Y}}d\hat{\psi}\big)(\hat{Y})}{|\hat{Y}|^2_{\hat{\omega}}}-K'\Big)^2+C_r,
 \end{equation*}
with $C_r$ depending only on $r$ (essentially, $C_r\leq\tfrac{C}{r^2}$ with $C$ depending only on $\omega_1$ and its derivatives up to order 3, $|d\psi|_{\omega_1}$, $\Delta_1\psi$, and hence does not depend on $r$, and neither does $K$ nor $K'$).
Now $\hat{A}$ reaches its maximum at 0, so the left-hand side of the latter inequality is nonnegative, hence an upper bound on $\hat{A}'$ by some $\tfrac{C}{r}$ with $C$ independent of $r$ and $\eta$. 
This gives the desired control on $M$. \hfill $\square$

~

Now using the techniques of \cite{ckns}, and working as usual in (half-)balls of quasi-coordinates instead of (half-)balls of coordinates, one can show:
 \begin{prop}  \label{prop_C2beta_estimate}
  There exists some $\beta\in(0,1)$ and some constant $C$ such that $\|\psi\|_{C^{2,\beta}}\leq C$ if $\psi$ is the solution of some $(E'_r)$, $r\in(0,1]$; more precisely, such $\beta$ and $C$ can be taken independent of $r$ if it stays away from 0.
 \end{prop}
\textit{Proof.} Cover $(\XD)\times\Sigma$ of (half-)balls $B^{(+)}\subset \C^m\times\C^{(+)}$ (where $\C^+=\{\mathfrak{Re}\geq0\}$) of quasi-coordinates $(z_1,\dots,z_m,z)$ of radius $\delta>0$ independent of $r$ such that:
 \begin{itemize} \setlength{\itemsep}{0pt}
  \item[$\bullet$] the collection of  (half-)balls of radius $\delta/2$ still covers $(\XD)\times\Sigma$;
  \item[$\bullet$] any point in $(\XD)\times\partial\Sigma$ is the center of a half-ball;
  \item[$\bullet$] the part $T=(\C^m\times\{0\})\cap B^+$ of a half-ball corresponds to $(\XD)\times\partial\Sigma$ i.e. if $\pi$ is one of the immersions associated to $B^+$ then $T=B^+\cap\pi^{-1}\big((\XD)\times\partial\Sigma\big)$;
  \item[$\bullet$] on any (half-)ball $\tfrac{1}{2}\omega_{euc}\leq\pi^*\omega\leq 2\omega_{euc}$, $\pi^*\theta(r)\geq cr$ and the derivatives of $\pi^*\omega$ and $\pi^*\theta(r)$ are bounded; all these controls are independent of $\pi$ and $r$;
  \item[$\bullet$] according to Proposition \ref{prop_C2_estimate}, we have bounds on the $\pi^*\psi$ up to order $(4,\eta)$ which are independent of $\pi$. Moreover those bounds are independent of $r$ on $|\pi^*\psi|$, $|d\pi^*\psi|$ and 
$\big|i\ddbar\pi^*\psi\big|$; they remain so on $|\nabla^2\pi^*\psi|$ as long as $r$ stays away from 0;  
 \end{itemize}
To get $C^{2,\beta}$ estimates on the balls, we write the pull-back of $(E'_r)$ as $F[\pi^*\psi]=0$ where 
 \begin{equation*}
  F[u]=\log\left[\det\Big((\pi^*\omega)_{j\bar{k}}+\frac{\partial^2u}{\partial z_j\partial \overline{z_k}}\Big)\right]-\log\big(\pi^*\theta(r)\big),u\in C^{2}_{loc}(B);
 \end{equation*}
this way $\pi^*\psi$ and $F$ satisfies the hypotheses of Theorem 17.14 in \cite{gt}, with in particular the ellipticity of $F$ coming from
 \begin{equation*}
  \sum_{j,k=1}^{2m+2}F^{jk}\xi_j\xi_k=|\xi|_{(\pi^*\omega')^{-1}},\xi\in\R^{2m+2}=\C^{m+1}
 \end{equation*}
with $(\pi^*\omega')^{-1}$ the (1,1)-form whose matrix in the coordinates of $B$  is the inverse of that of $\pi^*\omega'=\pi^*\big(\omega_1+i\ddbar\psi\big)$, the estimates on $i\ddbar\psi$ ensuring us about the existence of some $c>0$ independent of $r$ and $\pi$ such that:
  $cg_{euc}\leq (\pi^*g_{\phi})^{-1} \leq c^{-1}r^{-1}g_{euc}$.
The theorem gives us an estimate on the $|\nabla^2(\pi^*\psi)|_{C^{0,\beta}(\tfrac{1}{2}B)}$ with $\beta$ depending only on $\delta$, $\lambda$ et $\Lambda$ such that $\lambda\omega_{euc}\leq\pi^*\omega_1+i\ddbar(\pi^*\psi)\leq\Lambda\omega_{euc}$ and $|\nabla^2(\pi^*\psi)|_{C^0(B)}$, so that $\beta$ can be taken independent of $r$ if it stays away from $0$.

The case of (half-)balls is a bit more delicate; nonetheless, let us say some words about it. 
We want to apply Theorem 9.15 of \cite{gt}, and for this we need an estimate on the modulus of continuity of $\nabla(\pi^*\psi)$ around points of the boundary.
Applying techniques of \cite{ckns}, in particular those of \S 2.2, one gets
 \begin{lem}
 There exists a constant $C$ depending only on $|\pi^*\psi|_{C^2(B^+)}$, $\lambda$, $\Lambda$, $\pi^*\theta(r)$ --- so in particular $C$ can be taken independent of $r$ if it stays away from $0$, and independent of $\pi$  --- such that for all $z_0\in \tfrac{2}{3}T:=T\cap\tfrac{2}{3}B^+$ on has 
  \begin{equation*}
   \big|\nabla^2\pi^*\psi(z_0)-\nabla^2\pi^*\psi(z)\big|\leq\frac{C}{1+\big|\log|z-z_0|\big|}
  \end{equation*}
 for all $z\in B^+$ such that $|z-z_0|<\delta/3$.
 \end{lem}
Now differentiate the pullback of $(E'_r)$ with respect to some tangential operator $\mathcal{D}$ equal to $\pm\tfrac{\partial}{\partial x_j}$ or $\pm\tfrac{\partial}{\partial y_j}$, $1\leq j\leq m$ to get
 \begin{equation}  \label{diff_MA}
  \Delta_{\pi^*\omega'}(\mathcal{D}\pi^*\psi)=-\mathcal{D}\log \big(\pi^*\theta(r)\det (\pi^*g_1)_{j\bar{k}}\big)+\sum_{j,k=1}^{m+1}(\pi^*\omega')^{j\bar{k}}\mathcal{D}(\pi^*\omega_1)_{j\bar{k}}.
 \end{equation}
and apply Theorem 9.15 of \cite{gt} with $L=-\Delta_{\pi^*\omega'}$, $u=\mathcal{D}\pi^*\psi$ and $p>\tfrac{2m+2}{1-\beta}$ fixed. 
This gives us an estimation $\big|\mathcal{D}(\pi^*\psi)\big|_{L^{p,2}\big(\tfrac{2}{3}B^+\big)}\leq C$ with $C$ only depending on a lower bound on $r$. 
It is converted to a $C^{1,\eta}\big(\tfrac{2}{3}B^+\big)$ estimate on the $\mathcal{D}(\pi^*\psi)$ thanks to our choice of $p$, thus so far we control the $\tfrac{\partial^2\pi^*\psi}{\partial z_j\partial \overline{z_k}}$ on $C^{0,\eta}\big(\tfrac{2}{3}B^+\big)$, $1\leq j,k\leq m$. 
A similar control on $\tfrac{\partial^2\pi^*\psi}{\partial z\partial \overline{z}}$ comes from the very equation $(E'_r)$: develop the determinant with respect to the last column and express $\tfrac{\partial^2\pi^*\psi}{\partial z\partial \overline{z}}$ as a function of all the other terms. \hfill $\square$

 \subsection{Proofs of Proposition \ref{prop_r0=0} and Theorem \ref{thm_MA_hmg}}  \label{ccl_cont_meth}
  \subsubsection{Proof of Proposition \ref{prop_r0=0}}

Since equations $(E_r)$ and $(E'_r)$ are equivalent under the translation $\psi\mapsto\phi_1+\psi$, we can take $r_0$ as the infimum of the $r$ such that $(E'_{r'})$ admits a solution for all $r'\in(r,1]$. 
We first show that $r_0<1$, which is somehow the easy part, and then that $r_0$ cannot be positive, which uses the estimates we proved in sections \ref{uniq_sol}, \ref{higher_order_est} and \ref{C2_and_C2alpha_est}. 
Notice that uniqueness has already been proved in Proposition \ref{prop_unqss_C0_estimate}.

  \paragraph{Equation $(E'_r)$ admits (regular) solutions for $r$ close to 1.} 
  The vital remark here is the following: if $P$ denotes the operator
 \begin{align*}
  P:C^{\infty}_0 &\longrightarrow \Gamma^{\infty}(K_{(\XD)\times\Sigma}) \\
    \psi         &\longmapsto \big(\omega_1+i\ddbar\psi\big)^{m+1},
 \end{align*}
and $\psi$ is strictly $\omega_1$-pluri-subharmonic (i.e. $\omega_1+i\ddbar\psi>0$), then up to a $-\tfrac{1}{2}$ factor, the linearization of $P$ at $\psi$ is the Laplacian of $\omega_1+i\ddbar\psi$ multiplied by its volume from, that is:
 \begin{equation*}
  d_{\psi}P(\chi)=-\frac{1}{2}(\Delta_{\omega_1+i\ddbar\psi}\chi)\cdot\big(\omega_1+i\ddbar\psi\big)^{m+1},
 \end{equation*}
and this remains true when restricting $P$ to $C^{k+2,\beta}_0$ to $\Gamma^{k,\beta}(K)$, $(k,\beta)\in\N\times(0,1)$.
In particular $d_0P=-\tfrac{1}{2}(\Delta_{\omega_1}\cdot)\omega_1^{m+1}$ which is an isomorphism from $C^{4,\beta}_0$ to $\Gamma^{2,\beta}(K)$ (since $\tfrac{1}{2}\Delta_{\omega_1}$ is an isomorphism from $C^{4,\beta}_0$ to $C^{2,\beta}$, see remark \ref{rmk_isom1}). 
Take any $\gamma \in(0,1)$. 
Because of the latter isomorphism, and since $\theta(r)\in C^{2,\gamma}$, we know from the implicit functions theorem that $(E'_r)$ admits $C^{4,\gamma}_0$ solutions for $r$ close to 1, and away from 0 if necessary, say $r\in J$; the only point to be checked is that for such solutions $\psi_r$, $\omega'=\omega_1+i\ddbar\psi$ are equivalent to $\omega$ ("uniformly equivalent" is not necessary). 
Notice that $J$, $r\mapsto\psi_r$ is continuous for the $C^{4,\gamma}$ norm, and so is the function mapping those $r$ to the smallest eigenvalue of $\omega'$. 
Because $\theta(r)$ never vanishes, neither does this eigenvalue, which remains positive, as well as the other eigenvalues. 
So far there is no evidence for the existence of some $c>0$ such that $\omega'>c\omega_1$ \textit{globally on} $(\XD)\times\Sigma$ for all $r\in J$, but we can assume that $\|\psi_r\|_{C^{4,\gamma}}$, and in particular $\|\psi_r\|_{C^{2}}$ remains bounded for those $r$.
This tells us that there exists some $C>0$ such that for all $r\in J$, $\omega'\leq C\omega_1$; since $\theta(r)=\det^{\omega_1}(\omega')$ is positively and uniformly bounded below on $(\XD)\times\Sigma\times J$, it turns out that such a $c$ exists.

This is now a standard bootstrap argument to show that those solutions are $C^{\infty}$. 
Fix $r$ and choose some quasi-coordinate system like in the proof of Proposition \ref{prop_C2beta_estimate}; select a (half-)ball, with coordinate $(z_1,\dots,z_{m+1})$, and denote as usual by $\pi$ the associated immersion and $\mathcal{D}$ some first order differential operator, namely one of the $\partial_{x_j}$ or $\partial_{y_j}$, $j\in\{1,\dots,m+1\}$.
Differentiate the pulled back Monge-Ampère equation $(E'_r)$ with respect to $\mathcal{D}$; this writes $\Delta_{\pi^*\omega'}(\mathcal{D}\pi^*\psi)=f$, with $f$ as in \eqref{diff_MA} hence bounded up to order $(2,\gamma)$ independently of $\pi$.
Now $\Delta_{\pi^*\omega'}$ is an elliptic operator with $C^{2,\gamma}$ coefficients, and both its ellipticity (lower and upper bounds on its principal symbol) and the $C^{2,\gamma}$ bounds on the coefficients are independent of $\pi$. 
Standard Schauder estimates thus tell us that $\mathcal{D}\pi^*\psi$ is $C^{4,\gamma}$ on say the (half-)ball of half radius, and provide $C^{2,\gamma}$ on those smaller balls independent of $\pi$. 
Collecting all those regularity statements and estimates for all the $\mathcal{D}$ and $\pi$ in game, we get that $\psi\in C^{5,\gamma}\big((\XD)\times \Sigma\big)$: we have improved regularity by one order.
Going back to a (half-)ball of quasi-coordinate, the differentiate Monge-Ampère equation writes with an elliptic $C^{3,\gamma}$ operator and a $C^{3,\gamma}$ right-hand-side, with ellipticity and bounds independent on the immersion. 
We have this way $C^{5,\gamma}$ regularity and bounds on the $\mathcal{D}\pi^*\psi$ independent of $\pi$: $\psi\in C^{6,\gamma}\big((\XD)\times \Sigma\big)$. 
Going on this induction it is clear that $\psi\in C^{\infty}\big((\XD)\times \Sigma\big)$, for all $r\in J$.

  \paragraph{Equation $(E'_r)$ admits (smooth) solutions for all $r\in(0,1]$.}
Denote by $r_0$ the infimum of the $r\in(0,1]$ such that $(E'_{r'})$ for all admits a solution in $C^{\infty}_0\big((\XD)\times \Sigma\big)$ for all $r'\in(r,1]$.
We already know that $r_0<1$; let us suppose it is $>0$. 
Choose some sequence $(r_j)_{j\geq1}$ of elements of $(r_0,1]$ tending to $r_0$.
By Proposition \ref{prop_C2beta_estimate} we have some $\beta\in(0,1)$ and some constant $C$ such that $\|\psi_{r_j}\|_{C^{2,\beta}}\leq C$ for all $j\geq1$. Playing the same game as above, it is easy to provide a uniform $C^{4,\beta}$ bound on $(\psi_{r_j})$. 
Two diagonal extractions give us a $C^{4,\gamma}_{loc}$ converging subsequence with some $\gamma\in(0,\beta)$ to some function $\psi$; moreover the uniform $C^{4,\beta}$ bound on the whole $(\XD)\times\Sigma$ provides a uniform $C^{4,\gamma}$ bound which pass to the limit (use quasi-coordinates), hence $\psi\in C^{4,\gamma}$, and even $C^{4,\gamma}_0$.
By local $C^2$ convergence, $\omega_1+i\ddbar\psi\geq0$ and $\big(\omega_1+i\ddbar\psi\big)^{m+1}=\theta(r_0)\omega_1^{m+1}$.
Since $\psi\in C^2\big((\XD)\times\Sigma\big)$ we know from above that $\omega_1+i\ddbar\psi$ is mutually bounded with $\omega_1$. 
Then the bootstrap argument applies and we get that $\psi\in C^{\infty}\big((\XD)\times\Sigma\big)$.

To conclude, apply the implicit function theorem with $\omega_1+i\ddbar\psi$ replacing $\omega_1$; from this we know that there exists $C^{4,\gamma}$ solutions to $(E'_r)$ with $r$ in some neighbourhood $J$ of $r_0$. 
Shrinking $J$ if necessary, $0\notin J$, and it turns out as above that those solutions are in $C^{\infty}\big((\XD)\times\Sigma\big)$, which contradicts the definition of $r_0$, since $(r_0-\vareps,r_0]\subset J$ as soon as $\vareps>0$ is small enough.
Proposition \ref{prop_r0=0} is proved.

  \subsubsection{Proof of Theorem \ref{thm_MA_hmg}}

Theorem \ref{thm_MA_hmg} almost follows from Proposition \ref{prop_r0=0}, except for the uniform bounds on the $\Phi_{\vareps}$, which come from Proposition \ref{prop_C11_estimate}, and the statement about the limit obtained when letting $\vareps$ go to 0. 
This latter is understood in the theory developed in \cite{be-ta} and is an application of the monotonicity theorem in this paper; even if it is stated for a decreasing sequence of pluri-sub-harmonic functions, we can apply it to our sequence $(\Phi_{\vareps})$ which increases when $\vareps$ goes to 0. 
Take indeed an exhaustive sequence $(K_j)$ of compact subsets of $(\XD)\times\Sigma$, and a decreasing sequence of $(\vareps_j)$ going to 0 such that for every $j$, $m_j:=\sup_{K_j}\big|\Phi_{\vareps_j}-\Phi_{\vareps_{j+1}}\big|\leq\tfrac{1}{2^j}$. 
Then on every compact subset, $\big(\Phi_{\vareps_j}+\sum_{k\geq j-1}m_k\big)$ decreases from a certain rank, to the same limit as the $C^0_{loc}$-limit of the $(\Phi_{\vareps})$, and this limit satisfies \eqref{MA_hmg} in the sense of currents by the monotonicity theorem. 

This ends the proof of Theorem \ref{thm_MA_hmg}, and the present part.

\section{Calabi-Yau theorem on $\XD$ and negative Ricci forms}  \label{approx_CY}
 \subsection{Statement and motivation}  \label{approx_CY_motiv}

   In order to state properly Theorem \ref{thm_approx_CY}, which is a generalization of the celebrated Calabi-Yau theorem, we first need to introduce \textit{weighted Hölder spaces}, in which the decay of the functions 
is taken into account near the divisor. 
 \begin{df}
  Let $(k,\alpha)\in\N\times[0,1)$, $\gamma\in\R$. We set
  \begin{equation}  \label{df_wtd_holder}
   C^{k,\alpha}_{\gamma}=\big\{f\in C^{k,\alpha}_{loc}(\XD)|\,\rho^{\gamma}f\in C^{k,\alpha}(\XD)\big\}=\rho^{-\gamma}C^{k,\alpha}(\XD),
  \end{equation}
  where $C^{k,\alpha}(\XD)$ is that of section \ref{class_ric_form}. We endow this space with the obvious norm, denoted by $\|\cdot\|_{C^{k,\alpha}_{\gamma}}$.

  We also set $C^{\infty}_{\gamma}=\bigcap_{k\in\N,\alpha\in(0,1)}C^{k,\alpha}_{\gamma}$.
 \end{df}

Let us comment briefly this definition. The right hand side inequality in \eqref{df_wtd_holder} and \eqref{df_wtd_holder_forms} comes from the control on the derivatives of $\rho$, especially $|d\rho|_g$ is comparable to $\rho$ near $D$, and that $\nabla_{g}^k\rho=O(\rho)$ for any $k\geq1$. 
Notice that we can also compute norms using quasi-coordinates. 
For instance, if $U$ is a polydisc $(c\Delta)^k\times (\Delta)^{m-k}$ around a neighbourhood of a point of $D$ such that $D\cap U=\{(0,\dots,0)\}\times (\Delta)^{m-k}$ covered by a union $\bigcup_{\delta\in(0,1)^k}\Phi_{\delta}\big((\tfrac{1}{2}\Delta)^k\times (\Delta)^{m-k}\big)$ as in \S\ref{model} and if $f$ is $C^{k,\alpha}_{loc}$, with support in $U$, then
 \begin{equation*}
  \|f\|_{C^{k,\alpha}_{\gamma}}\sim\sup_{\delta\in(0,1)^k} \frac{1}{\big((1-\delta_1)\cdots({1-\delta_k})\big)^{\gamma}}
                                                          \|{\Phi_{\delta}}^{*}f\|_{C^{k,\alpha}(\mathcal{P}_k)}
 \end{equation*}
(where $\mathcal{P}_k=(\tfrac{1}{2}\Delta)^k\times (\Delta)^{m-k}$) because ${\Phi_{\delta}}^{*}\rho$ is uniformly mutually bounded with $\tfrac{1}{(1-\delta_1)\cdots(1-\delta_k)}$ on $\mathcal{P}_k$ for $\delta\in(0,1)^k$, as we already saw it in \S\ref{model}.

~

We can now state the following "logarithmic" version of the Calabi-Yau theorem (see for instance \cite[ch.5]{joy} for a review on the Calabi conjecture and its resolution by Yau):
 \begin{thm}  \label{thm_approx_CY}
  Let $\omega'\in\mom$, $\nu>0$ and $f\in C^{\infty}_{\nu}(\XD)$ such that $\int_{\XD} e^f\vol^{\omega'}=\vl$. 
Then there exists $\varphi\in C^{\infty}(\XD)$ such that $\big(\omega'+i\ddbar\varphi\big)^m=e^f(\omega')^m$. 
More precisely, $\varphi$ is for all $k\geq 0$ a $C^k_{loc}$-limit of $(\varphi_{\vareps})_{0<\vareps\leq 1}$ when $\vareps$ goes to 0, where $\big(\omega'+i\ddbar\varphi_{\vareps}\big)^m=e^{f+\vareps\varphi_{\vareps}}(\omega')^m$ for all $\vareps>0$. 
Moreover there are $C^k$-bounds independent of $\vareps$ on those $\varphi_{\vareps}$, and there exists $c>0$ such that $\varphi_{\vareps}\in C^{\infty}_{c\vareps}(\XD)$ when $\vareps$ is small enough.
 \end{thm}

This approach of $\vareps$-perturbed Monge-Ampère equation is quite close to that of \cite{tian-yau2} and \cite[part 4]{hein}.

We postpone the proof to part \ref{prf_approx_CY} below. 
The existence of the family $(\varphi_{\vareps})_{\vareps>0}$ with elements in $C^{\infty}(\XD)$ is not new, and follows from \cite{tian-yau1}; actually, they do it with $\Omega=2\pi(K[D])$ with $K[D]$ assumed ample and $\vareps=1$, but what really matters here for $\omega'$ is being of Poincaré type, and that $\vareps>0$. 
It also follows from this work that $|\varphi_{\vareps}|_{C^0}\leq \vareps^{-1}|f|_{C^0}$ for all $\vareps\in(0,1]$.
However, new  are the uniform $C^k$ bounds, and that the $\varphi_{\vareps}$ lie in positively weighted Hölder spaces. 

An interesting observation is the following: for $\vareps>0$, $\varrho_{\omega'+i\ddbar\varphi_{\vareps}}=\varrho_{\omega'}-i\ddbar f-\vareps i\ddbar\varphi_{\vareps}$, which tends \textit{at any order uniformly on $\XD$} to $\varrho_{\omega'}-i\ddbar f$. 
In other words, suppose that $\varrho_{\omega'}-i\ddbar f$ is "interesting" in some sense; then we can realize it as the Ricci form of a metric differing from $\omega'$ by some fast decaying potential, up to an arbitrary small error term in the $\Gamma^{\infty}(\Lambda^{1,1},\XD)$ topology. 
More concretely, our theorem allows us to construct metrics with Ricci form strictly negative in the Poincaré sense:
 \begin{thm}  \label{thm_negative_ricci}
  Assume $K[D]$ is ample on $X$. Then there exists $\varpi\in\mom$ such that $\varrho_{\varpi}\leq-c\varpi$ for some positive constant $c$.
 \end{thm}
The proof is rather long, so the next section is devoted to it.

 \subsection{Proof of Theorem \ref{thm_negative_ricci}}

Before starting, we shall mention that we proceed by induction on the highest codimension in $X$ of the crossings of $D$. We shall also introduce more functional spaces, as our weighted Hölder spaces defined so far fail to contain the functions appearing in the upcoming proof. 
 \begin{df}
  Let $g$ be a metric $C^{\infty}$-quasi-isometric to the model $\omega$ of section \ref{model} and let $(k,\alpha)\in\N\times[0,1)$, $\gamma\in\R^{+}$. 
Given $v_1,\dots,v_n$ such that $v_j\equiv 1$ in a neighbourhood of the connected component $\mathcal{D}_j$ of $D$ and $v_j\equiv 0$ in the neighbourhood of $\mathcal{D}_l$ if $l\neq j$ for $j=1,\dots,n$ (so that $D=\bigsqcup_{j=1}^n\mathcal{D}_j$), we set
   \begin{equation*}
  E^{k,\alpha}_{\gamma}(g)=\Big\{f\in C^{k,\alpha}_{\gamma}\oplus\bigoplus_{j=1}^n \R v_j| \int_{\XD}f\vol^g=0 \Big\}.
 \end{equation*}
If $\gamma>0$, 
we set $\|f\|_{E^{k,\alpha}_{\gamma}(g)}=\|h\|_{C^{k,\alpha}_{\gamma}}+\sum_j|a_j|$ (we get each $a_j$ back as the limit of $f$ near $\mathcal{D}_j$).
 \end{df}

Those spaces are indeed relevant in the weighted $\ddbar$-lemma we are going to use in the proof of Theorem \ref{thm_negative_ricci}, as the spaces where lie the $\ddbar$-potentials of real closed (1,1)-forms which are $O(\rho^{-\delta})$ at any order for some $\delta>0$, as described in the weighted $\ddbar$-lemma (Proposition \ref{prop_weighted_ddbar_lemma}) stated and proved in section \ref{wtd_ddbar_lemma} below. 

  \subsubsection{The smooth divisor case}

As aforementioned, we start when the codimension of the crossings equals 1 in $X$, meaning actually there are no proper crossings, but instead that $D$ is smooth.
Choose some smooth negative $\varrho_0\in -2\pi c_1(K[D])$, and remember $\omega_0$ was a smooth Kähler form on $X$, such that $\Omega=[\omega_0]_{dR}$. 
The adjunction formula says, if $D=\sum_{j=1}^N D_j$ is the decomposition of $D$ into irreducible \textit{disjoint} components, for all $j$: 
 \begin{align*}
  K[D]|_{D_j}&=(K\otimes[D_1]\otimes\cdots\otimes[D_N])|_{D_j}\\
             &=(K\otimes[D_j])|_{D_j}\otimes\big([D_1]\otimes\cdots\otimes\widehat{[D_j]}\otimes\cdots\otimes[D_N]\big)|_{D_j}\\
             &\cong K_{D_j}\otimes1 =K_{D_j},
 \end{align*}
so that $\varrho_0|_{D_j}$ (meaning "the closed form induced in $\Lambda^{1,1}_{D_j}$") is in $-2\pi c_1(K_{D_j})$. 
Now for all $j$, the Calabi-Yau theorem for smooth Kähler compact manifolds applies on $D_j$ which is smooth and compact, and provides some potential $\psi_j\in C^{\infty}(D_j)$ such that $\varrho_0|_{D_j}=\varrho_{\omega_0|_{D_j}+i\ddbar\psi_j}$. 
Denote by $p_j$ the projection on $D_j$, defined in a tubular neighbourhood $\mathcal{N}_j$ of $D_j$, and by $\chi_j$ a smooth function equal to 1 in a small neighbourhood of $D_j$, with support in $\mathcal{N}_j$. 
This way, $\varphi:=\sum_{j=1}^n\chi_j{p_j}^*\psi_j$ is well defined and smooth on $X$; moreover, $\omega_0+i\ddbar\varphi$ induces $\omega_0|_{D_j}+i\ddbar\psi_j$ on every $D_j$ as soon as the $\mathcal{N}_k$ are disjoint.

The point is that this closed real (1,1)-form $\omega_0+i\ddbar\varphi$ has no reason in general to be positive; nevertheless, the lack of positivity is essentially in the direction normal to $D$, so that it can be corrected, in our Poincaré metrics setting, by "$\log\log$ potentials". 
More explicitly, let $\chi_0:\R\rightarrow[0,1]$ such that $\chi_0\equiv0$ on $(-\infty,0]$ and $\chi_0\equiv1$ on $[1,+\infty)$. 
Remember that $u_j=\log(\lambda+\rho_j)$ for some $\lambda\geq 0$, and that there we can assume $\rho_j$ to be constant on the $\mathcal{N}_k$ for all $1\leq j\neq k\leq N$; take also $A_1,\dots,A_N>0$. 
In those conditions, we claim that 
 \begin{equation*}
  \omega':=\omega_0-\sum_{j=1}^N A_ji\ddbar u_j+\sum_{j=1}^Ni\ddbar\big(\chi_0(u_j^{1/2}-K){p_j}^*\psi_j\big)\in\mom
 \end{equation*}
when $\lambda$ and $K$ are big enough. 
It even turns out that $\omega'=\omega_0|_{D_j}+i\ddbar\psi_j+\tfrac{A_jidz\wedge d\overline{z}}{|z|^2\log^2(|z|^2)}+O^{\infty}(\rho_j^{-1})$ in any neighbourhood of any point of $D_j$ in which $D_j$ is given by $\{z=0\}$, this last assertion being independent of $\lambda$ and $K$. 
These asymptotics being showed in the same way than those of Proposition \ref{prop_asymptotics_omega}, the only point to be checked is that $\omega'>0$ on $\XD$ for a suitable choice of $\lambda$ and $K$. 
Since the $u_j$ are constant near the $D_k$, $k\neq j$, we can assume that $N=1$ to show this positivity, and we drop the $j$ indexes. 
First fix $\lambda\geq0$ big enough so that $\omega_0-Ai\ddbar u>0$ on $\XD$. 
Then take $\vareps\in(0,\tfrac{1}{4})$ small enough so that $\omega_0|_D+i\ddbar\psi\geq 4\vareps\omega_0|_D$, that is $i\ddbar\psi\geq (4\vareps-1)\omega_0|_D$. 
If one takes a collection of open sets of coordinates which in $X$  are neighbourhoods of open sets covering $D$ and in which $D$ is given by $z=0$, we can assume those neighbourhoods small enough so that $i\ddbar(p^*\psi)\geq -(1-3\vareps)\omega_0-Cidz\wedge d\overline{z}$. 
Since $\chi_0$ takes its value in $[0,1]$, we will have similarly that $\chi_0(u^{1/2}-K)i\ddbar(p^*\psi)\geq -(1-3\vareps)\omega_0-Cidz\wedge d\overline{z}$ on our open sets, whose union is denoted by $V$. 
Then we take $K$ big enough so that $V_K:=\{u\geq K^2\}\subset V$, and this way on $X\backslash V_K$, $\chi_1(u^{1/2}-K)$ is 0 and $\omega'$ equals $\omega_0-Ai\ddbar u$, hence is $>0$. 
Thus, it suffices to show $\omega'>0$ on $V_K$, and our lower bound on $\chi_1(u^{1/2}-K)i\ddbar(p^*\psi)$ goes in this sense.

Indeed, since we have besides:
 \begin{equation*}
  i\ddbar\chi_1(u^{1/2}-K)=\chi_1''(u^{1/2}-K)\frac{i\partial u\wedge\dbar u}{4u}+\chi_1'(u^{1/2}-K)\Big(\frac{i\ddbar u}{2u^{1/2}}-\frac{3i\partial u\wedge\dbar u}{4u^{3/2}}\Big)
 \end{equation*}
we can again, up to increasing $K$ once more, assume we have $\big|(p^*\psi) i\ddbar\chi_1(u^{1/2}-K)\big|\leq \vareps(\omega_0-Ai\ddbar u)$ on $V_K$, remembering that $du$ and $i\ddbar u$ are bounded for Poincaré type metrics. 
Similarly, we can assume that
 \begin{equation*}
  \big|i\big(\partial\chi_1(u^{1/2}-K)\wedge\dbar(p^*\psi) +\partial(p^*\psi)\wedge\dbar\chi_1(u^{1/2}-K)\big)\big|\leq\vareps(\omega_0-Ai\ddbar u)
 \end{equation*}
on $V_K$. 
Finally, on $V_K$, or rather in its intersection with any of our open sets of coordinates assumed small enough so that $-Ai\ddbar u\geq \tfrac{Aidz\wedge d\overline{z}}{2|z|^2\log^2(|z|^2)}-\vareps\omega_0\geq2Cidz\wedge d\overline{z}-\vareps\omega_0$ up to increasing $K$ once again, we have the minoration
 \begin{align*}
  \omega_1 =  &\omega_0-Ai\ddbar u +\chi_1(u^{1/2}-K)i\ddbar(p^*\psi)\\
              &+i\big(\partial\chi_1(u^{1/2}-K)\wedge\dbar(p^*\psi) +\partial(p^*\psi)\wedge\dbar\chi_1(u^{1/2}-K)\big)\\
              &+(p^*\psi) i\ddbar\chi_1(u^{1/2}-K)\\
          \geq&(1-2\vareps)(\omega_0-Ai\ddbar u)+\chi_1(u^{1/2}-K)i\ddbar(p^*\psi)\\
          \geq&\vareps\omega_0-Cidz\wedge d\overline{z}-(1-2\vareps)Ai\ddbar u \\
          \geq&2\vareps^2\omega_0+(1-4\vareps)Cidz\wedge d\overline{z} \,\,\,\text{ car }-Ai\ddbar u\geq-\vareps\omega_0+2Cidz\wedge d\overline{z},
 \end{align*}
which is positive. 

Having dealt with that point, thanks the asymptotics of $\omega'$, it is easy to compute asymptotically its Ricci form; schematically, it writes $\varrho_{\omega_0|_{D_j}+i\ddbar\psi_j}-\tfrac{2idz\wedge d\overline{z}}{|z|^2\log^2(|z|^2)}+O^{\infty}(\rho_j^{-1})$ near each $D_j$. 
These asymptotics are exactly those of $\varrho_0-2\sum_{j=1}^Ni\ddbar u_j$, which we can suppose $\leq -c\omega$ on $\XD$ for some $c>0$ for the same reasons than $\omega$ is Kähler of Poincaré type. 
In a nutshell, $\varrho_{\omega'}+\big(\varrho_0-2\sum_{j=1}^Ni\ddbar u_j\big)\in \Gamma^{\infty}_1(\Lambda^{1,1})$, and this form lives in the zero cohomology $L^2$ class. 
Applying the weighted $\ddbar$ lemma (Proposition \ref{prop_weighted_ddbar_lemma}), we thus can write $\varrho_{\omega'}+\big(\varrho_0-2\sum_{j=1}^Ni\ddbar u_j\big)=i\ddbar f$ for some $f\in E^{\infty}_{1}(\omega')$ and $c\in\R$ so that $\int_{\XD} e^{f+c}\vol^{\omega'}=\vl$ ; we will not use $\int_{\XD}f\vol^{\omega'}=0$, so we can assume $c=0$. 
This function $f$ has no reason to tend to 0 near $D$; nonetheless, we can correct it in a compact subset of $\XD$ so that it does, and so that $\rho_{\omega'}+i\ddbar f'\leq -c'\omega$ and $\int_{\XD}e^{f'}\vol^{\omega'}=\vl$.

To do so, set $a_j:=\lim_{D_j}f$ for all $j=1,\dots,N$. 
Start by assuming $a_1>0$ (if $a_1=0$, consider $a_2$, and if $a_1<0$, the technique is the same). 
According to the beginning of the proof, when $K$ is big enough, and $\mu\geq 1$, then
 \begin{equation*}
  \varrho_{\omega'}-i\ddbar f -a_1i\ddbar \Big(-\chi_0\big(u_j^{1/2}-K\big)+\chi_0\big(u_j^{1/2}-\mu K\big)+\chi_0\big(u_j^{1/2}-(\mu+1)K\big)\Big)\leq -c_1\omega
 \end{equation*}
for some $c_1>0$ independent of $\mu$. 
Set $f_{1,\mu}=f+a_j\Big(-\chi_0\big(u_j^{1/2}-K\big)+\chi_0\big(u_j^{1/2}-\mu K\big)+\chi_0\big(u_j^{1/2}-(\mu+1)K\big)\Big)$, so that $f_{1,\mu}\in E^{\infty}_{1}(\omega')$ up to its mean, $f_{1,\mu}$ goes to 0 near $D_1$ and to $a_j$ near $D_j$, $j\geq 2$. 
Moreover, $\mu\mapsto\int_{\XD}e^{f_{1,\mu}}\vol^{\omega'}$ is continuous on $[1,+\infty)$, is strictly greater than $\vl$ for $\mu=1$ and its limit is strictly less than $\vl$ when $\mu$ goes to $+\infty$ (if $a_j<0$, the inequalities are inversed). 
Hence there exists some $\mu_1$ so that $\int_{\XD}e^{f_{1,\mu_1}}\vol^{\omega'}=\vl$. 
Repeat this construction near $D_2, \dots,D_N$, to get a function $f'\in C^{\infty}_{1}(\XD)$ such that $\int_{\XD}e^{f'}\vol^{\omega'}=\vl$ and $\varrho_{\omega'}-i\ddbar f'\leq -c'\omega$. 
Denote by $\eta$ the difference $i\ddbar f-i\ddbar f'$, so that 
 \begin{equation} \label{formula_eta}
  \eta=i\ddbar\Big[\sum_{j=1}^Na_j\Big(-\chi_0\big(u_j^{1/2}-K\big)+\chi_0\big(u_j^{1/2}-\mu_j K\big)+\chi_0\big(u_j^{1/2}-(\mu_j+1)K\big)\Big)\Big]
 \end{equation}
with well-chosen $\mu_j$ (this will be useful below), and notice it has compact support on $\XD$.
 
Now apply the first part of Theorem \ref{thm_approx_CY} to $\omega'$ and $f'$; then $\varrho_{\omega'+i\ddbar\varphi}=\varrho_{\omega'}-i\ddbar f'\leq -c'\omega$. Set $\varpi=\omega'+i\ddbar\varphi$ to conclude.
Notice that we could have applied the second part with $\vareps>0$ small enough, and still get $\varrho_{\omega'+i\ddbar\varphi_{\vareps}}=\varrho_{\omega'}-i\ddbar f'-\vareps i\ddbar\varphi_{\vareps}\leq -c'_{\vareps}\omega$ with $c'_{\vareps}>0$. 
Notice moreover that we can take arbitrary positive $A_j$, in particular we can take them \textit{equal}.

  \subsubsection{Proof of Theorem \ref{thm_negative_ricci}: the general case}

Assume now there exist some codimension 2 crossings, and that it is the highest possible codimension. 
In what precedes, we first solve Calabi problem on the divisor, and then construct a potential on $\XD$ from the data of potentials on the divisor. 
We are following here the same process, now we know approximately how to solve Calabi problem when the divisor is smooth. 
For the sake of simplicity, assume that $D=D_1+D_2$, and that the decomposition of $D'=D_1\cap D_2$ into irreducible smooth components writes $\sum_{j=1}^{N'}D'_j$, and observe that $D_1\backslash D'$, $D_2\backslash D'$ are endowed with Poincaré type Kähler metrics, namely $\omega|_{D_1\backslash D'}$ and $\omega|_{D_2\backslash D'}$.
Once again, the adjunction formula applies nicely to give
 \begin{equation*}
  K_X[D]|_{D_1}\cong K_X[D_1]|_{D_1}\otimes [D_2]|_{D_1}\cong K_{D_1}\otimes [D']_{D_1},
 \end{equation*}
that is the intrinsic $K_{D_1}[D']$ on $D_1$.  
One step further we have 
 \begin{equation*}
  K_X[D]|_{D'_j}=(K_X[D]|_{D_1})|_{D'_j}\cong(K_{D_1}[D'])|_{D'_j}\cong K_{D'_j}
 \end{equation*}
for all $j=1,\dots, N'$. Thus $\varrho_0|_{D'_j}\in -2\pi c_1(K_{D'_j})$ as soon as $\varrho_0$ is smooth in $-2\pi c_1(K[D])$. Take such a $\varrho_0$, such that $\varrho_0<0$ on $X$.

Set as usual $u_j=\log(\lambda+\rho_j)$ on $X\backslash D_j$, $j=1,2$ ($\lambda\geq 0$ adjustable), and notice that when $k\neq j$, $u_j|_{D_k}$ plays the role of $\sum_{l=1}^{N'} u'_{k,l}$ on $D_k$ where $u'_{k,l}$ would be defined on $D_k$ as a function with "$\log\log$" behaviour near $D'_l$. 

As for the smooth divisor case, if $\psi_j\in C^{\infty}(D'_j)$ is such that $\omega_0|_{D'_j}+i\ddbar\psi_j$ has Ricci form $\varrho_0|_{D'_j}$ (Calabi-Yau theorem for smooth manifolds), we can extend it as a smooth function $\tilde{\psi}$ on $X$ so that 
 \begin{align*}
  \omega_1 &:=\omega_{0}|_{D_1}+i\ddbar (\tilde{\psi}|_{D_1})-i\ddbar (u_2|_{D_1\backslash D_2}) \\
 \text{ and } &\omega_2:=\omega_{0}|_{D_2}+i\ddbar (\tilde{\psi}|_{D_2})-i\ddbar (u_1|_{D_2\backslash D_1})
 \end{align*}
are Poincaré type metrics, respectively on $D_1\backslash D_2$ and $D_2\backslash D_1$, with respective asymptotics $\omega_0|_{D'_j}+i\ddbar\psi_j+\tfrac{idw\wedge d\overline{w}}{|w|^2\log^2(|w|^2)}+O^{\infty}(\rho_2^{-1})$ near $D'_j$ given by $\{w=0\}$ in $D_1$, 
and $\omega_0|_{D'_j}+i\ddbar\psi_j+\tfrac{idz\wedge d\overline{z}}{|z|^2\log^2(|z|^2)}+O^{\infty}(\rho_1^{-1})$ near $D'_j$ given by $\{z=0\}$ in $D_2$, and moreover such that 
their Ricci forms have respective asymptotics $\varrho_0|_{D'_j}-\tfrac{2idw\wedge d\overline{w}}{|w|^2\log^2(|w|^2)}+O^{\infty}(\rho_2^{-1})$ near $D'_j$ in $D_1$ and 
$\varrho_0|_{D'_j}-\tfrac{2idz\wedge d\overline{z}}{|z|^2\log^2(|z|^2)}+O^{\infty}(\rho_1^{-1})$ near $D'_j$ in $D_2$. 

Now applying the construction of the previous paragraph, we find $\varphi_1\in \varphi_1\in C_{\gamma}^{\infty}(D_1\backslash D')$ and $\varphi_2\in C_{\gamma}^{\infty}(D_2\backslash D')$
such that 
 \begin{align*}
  \varrho_{\omega_1+i\ddbar\varphi_1}=&(\varrho-i\ddbar u_1-i\ddbar u_2)|_{D_1\backslash D_2}+\vareps i\ddbar\varphi_1 +\eta_2|_{D_1}\\
 \text{ and } \varrho_{\omega_2+i\ddbar\varphi_2}&=(\varrho-i\ddbar u_1-i\ddbar u_2)|_{D_2\backslash D_1}+\vareps i\ddbar\varphi_2 +\eta_1|_{D_2}
 \end{align*}
with $\vareps$ arbitrarily small (and $\gamma=\gamma(\vareps)$) and $\eta_1$ with compact support in $\XD_2$ arbitrarily small in $C^{\infty}(\XD_2)$, $\eta_2$ with compact support in $\XD_1$ arbitrarily small in $C^{\infty}(\XD_1)$, constructed as $\eta$ in \eqref{formula_eta} (notice that the formulas respectively make sense on the whole $\XD_2$ and $\XD_1$). 
This only changes the asymptotics near the $D'_j$ by putting an exponent $-\gamma$ instead of an exponent $-1$ in the $O^{\infty}$. 

Now consider a function $\varphi$ on $\XD$ such that $\varphi$ is the sum of a function in $C^{\infty}_{\gamma}(\XD)$ and a smooth function on $X$, such that $\varphi|_{D_1\backslash D_2}=\varphi_1$, $\varphi|_{D_2\backslash D_1}=\varphi_2$, and $\omega_0+i\ddbar\varphi-i\ddbar u_1-i\ddbar u_2>0$ on $\XD$. 
We observe then that its Ricci form differs from $\rho_0+2(i\ddbar u_1+i\ddbar u_2)+\eta_1+\eta_2-\vareps i\ddbar\varphi\leq -c\omega$ by some $O^{\infty}(\rho^{-\gamma})$. 

Following the same process than in what precedes (weighted $\ddbar$-lemma, correction of constants near $D$, approximate Calabi-Yau theorem,), there exists some $\varphi'$ in $C^{\infty}_{\gamma'}(\XD)$, $\gamma'>0$, such that $\omega_0+i\ddbar(\varphi+\varphi')-i\ddbar u_1-i\ddbar u_2>0$ on $\XD$ and its Ricci form is arbitrarily close to $\rho_0+2(i\ddbar u_1+i\ddbar u_2)+\eta_1+\eta_2-\vareps i\ddbar\varphi$, and 
can in particular be taken $-c\omega$ for some $c>0$. 
This rules out the simplest codimension 2 case. 

The proofs of the cases where there are more $D_j$ (with possibly some disjoint from the others) or where the codimensions of the crossings are higher are just careful repetitions of the techniques used here. \hfill $\square$

  \subsection{The weighted $\ddbar$-lemma}   \label{wtd_ddbar_lemma}
We precise that everything in this section is independent of the ampleness of $K[D]$.
We formalize what is a real (1,1)-form which is a $O(\rho^{-\gamma})$ at any order (or at order $(k,\alpha)$) by the following:
 \begin{df}
  Let $(k,\alpha)\in\N\times[0,1)$, $\gamma\in\R$. We set:
   \begin{equation*}  \label{df_wtd_holder_forms}
     \Gamma^{k,\alpha}_{\gamma}\big(\Lambda^{(1,1)}\big)=\rho^{-\gamma}\Gamma^{k,\alpha}\big(\Lambda^{(1,1)}\big),
   \end{equation*}
  and endow the it with the obvious norm. We also set $\Gamma^{\infty}_{\gamma}\big(\Lambda^{1,1}\big)=\rho^{-\gamma}\Gamma^{\infty}\big(\Lambda^{1,1}\big)$.
\end{df}

The result we used in the proof of Theorem \ref{thm_negative_ricci} writes:
 \begin{prop}[Weighted $\ddbar$-lemma]  \label{prop_weighted_ddbar_lemma}
  Let $(k,\alpha)\in\N\times(0,1)$, $\eta\in \Gamma^{k,\alpha}_{\beta}\big(\Lambda^{(1,1)}\big)$ an $L^2$ exact 2-form, $\beta>0$, and $\varphi$ the $\ddbar$-potential of $\eta$ with zero mean w.r.t. some Kähler metric of Poincaré type $\omega'$.
 Then $\varphi$ is in fact in $E^{k+2,\alpha}_{\beta}(\omega')$, and there exists a constant $C=C(\beta,k,\alpha,\omega')$ such that $\|\varphi\|_{E^{k+2,\alpha}_{\beta}}\leq C\|\eta\|_{\Gamma^{k,\alpha}_{\beta}(\Lambda^{1,1})}$.
 \end{prop}
\prf. We decompose it with the help of three intermediate lemmas: 
 \begin{lem} \label{lem_pert_lapl}
  Let $g$ be a Kähler metric $C^{\infty}$-quasi-isometric to the model $\omega$. 
There exists a constant $c=c(g)>0$ such that for any $(k,\alpha)\in\N\times(0,1)$, $\vareps\in(0,1]$ and $\gamma\in[0,c\vareps)$ the $\vareps$-perturbed Laplacian $\Delta_g+\vareps:C^{k+2,\alpha}_{\gamma}\rightarrow C^{k,\alpha}_{\gamma}$ is an isomorphism.
 \end{lem}
\textit{Proof of Lemma \ref{lem_pert_lapl}}. For $\gamma$ and $\vareps>0$, one has to check that the conjugate operator $\mathcal{L}_{\gamma,\vareps}=\rho^{\gamma}(\Delta_g+\vareps)(\rho^{-\gamma}\cdot)=\Delta_g+\big(\vareps-\gamma\tfrac{\Delta_g\rho}{\rho}-2\gamma(\gamma+1)\big|\tfrac{d\rho}{\rho}\big|_g^2\big)+2\gamma\big(\cdot,\tfrac{d\rho}{\rho}\big)_g$ is an isomorphism from $C^{k+2,\alpha}$ to $C^{k,\alpha}$. 
Following \cite{tian-yau1}, p.589, this is true when $\sup_{\XD}\big\{\gamma\tfrac{\Delta_g\rho}{\rho}+2\gamma(\gamma+1)\big|\tfrac{d\rho}{\rho}\big|_g^2\big\}<\vareps$ (the first order term in $\mathcal{L}_{\gamma,\vareps}$ does not matter). 
Taking $A=\sup_{\XD}\tfrac{\Delta_g\rho}{\rho}$, $B=\sup_{\XD}\big|\tfrac{d\rho}{\rho}\big|_g$, this latter inequality is easy to check for all $\vareps\in(0,1]$ and $\gamma\in[0,c\vareps)$ where $c=\tfrac{4B}{2(A+B)\sqrt{4B+(A+B)^2}}$. \hfill $\blacksquare$

 \begin{lem} \label{lem_thm_weighted_ddbar_lemma1}
  Under the assumptions of Proposition \ref{prop_weighted_ddbar_lemma}, $\varphi\in C^{k+2,\alpha}$, and there exists a constant $C$ such that $\|\varphi\|_{C^{k+2,\alpha}}\leq C\|\eta\|_{\Gamma^{k,\alpha}_{\beta}}$.
 \end{lem}
\textit{Proof of Lemma \ref{lem_thm_weighted_ddbar_lemma1}}. It uses a perturbed Moser's iteration scheme, with parameter $\vareps$. Namely, define $\varphi_{\vareps}$ as the solution of $\Delta_{\omega'}\varphi_{\vareps}+\vareps\varphi_{\vareps}=-2\tr^{\omega'}(\eta)$, given by Lemma \ref{lem_pert_lapl}. Once noticed that $\int_{\XD}\varphi_{\vareps}\vol^{\omega'}=0$ (integrate the equation satisfied by $\varphi_{\vareps}$), just copy word by word the proof of Proposition \ref{prop_CY_C0_est} below, replacing $1-e^{f+\vareps\varphi_{\vareps}}$ by $-2\tr^{\omega'}(\eta)-\vareps\varphi_{\vareps}$ (again, the $\vareps$ are not a problem, and merely play in our favor; for instance, $\int_{\XD}\big|d\varphi_{\vareps}\big|^2\vol^{\omega'}=\int_{\XD}\varphi_{\vareps}\Delta\varphi_{\vareps}\vol^{\omega'}=\int_{\XD}\varphi_{\vareps}(-2\tr^{\omega'}(\eta)-\vareps\varphi_{\vareps})\vol^{\omega'}\leq  \|2\tr^{\omega'}(\eta)\|_{L^2}\|\varphi_{\vareps}\|_{L^2}$), 
$T_{\vareps}$ by $(\omega')^{m-1}$, and noticing constants $C$ and $C'$ of Proposition \ref{prop_p_to_espilon_p} , which depend on $f$, can be replaced by constants independent of $\eta$ times $\|\eta\|_{\Gamma^{k,\alpha}_{\beta}(\Lambda^{1,1})}$. 
Then notice $\varphi$ is a $C^0_{loc}$-limit of $(\varphi_{\vareps})_{\vareps>0}$ with $\vareps$ going to 0. 
\hfill $\blacksquare$

 \begin{rmk}
  This could appear a bit short, but we have preferred to develop the computations of such a Moser's iteration scheme in the slightly more difficult case that is Theorem \ref{thm_approx_CY}. 
Notice that both proofs use the Sobolev embedding (Lemma \ref{lem_sob_inj}) stated in paragraph \ref{prf_approx_CY_C0_est}.
 \end{rmk}

Next, we come to the technical core of the proof: 
 \begin{lem}  \label{lem_thm_weighted_ddbar_lemma2}
  Under the assumptions of Proposition \ref{prop_weighted_ddbar_lemma}, set $\beta'=\min\{2,\beta\}$. Then $\varphi=\psi+\sum_k a_kv_k$ with $\big\|\sum_j\big|\log|\sigma_j|\big|^{\beta'}\psi\big\|_{C^0}\leq C\|\eta\|_{C^{0,\alpha}_{\beta}}$.
 \end{lem}
\textit{Proof of Lemma \ref{lem_thm_weighted_ddbar_lemma2}}. We limit ourselves to the case of codimension at most 2 of the crossings for the sake of simplicity. 
Notice that assuming this proposition, the $a_k$ are automatically controlled as in the statement of Proposition \ref{prop_weighted_ddbar_lemma}, because for all $k$, $a_k=\lim_{x\rightarrow \mathcal{D}_k} (\varphi-\psi)(x)$, and we already control $\|\varphi\|_{C^0}$.
We choose a connected component $\mathcal{D}_k$ of $D$, which we can split into smooth irreducible components $D_1,\dots,D_j$. We first work around $D_1$, which we cover with polydiscs of coordinates $\{|z_1|,\dots,|z_m|\leq \tfrac{1}{e}\}$ where $D_1$ is given by $z_1=0$, and in case of a crossing the other component is given by $\{z_2=0\}$. 
Now if we choose one of these polydiscs, $\mathcal{P}$ say, two situations can occur:
 \begin{enumerate}
  \item \textit{There is no crossing in} $\mathcal{P}$. 
Write $z_1=re^{i\theta}=e^{-e^t+i\theta}$, $z'=(z_2,\dots,z_m)$, and equip each punctured disc $\{0<|z_1|\leq\tfrac{1}{e}\}$ with the standard cusp metric $dt^2+e^{-2t}d\theta^2$. 
The equation $\big(i\ddbar\varphi\big)_{1\bar{1}}=\eta_{1\bar{1}}:=\tfrac{f}{|z_1|^2\log^2(|z_1|^2)}$ rewrites $\big((\partial_t^2-\partial_t)+e^{2t}\partial_{\theta}^2\big)\varphi=f$, with $|f(z_1,z')|\leq C\|\eta\|_{\Gamma^0_{\beta}}e^{-\beta t}$ where $C$ depends only on our polydisc.
Now decompose $\varphi$ into $\varphi_0+\varphi_{\perp}$, with $\varphi_0$ invariant with respect to $\theta$, and $\varphi_{\perp}$ orthogonal to the constants on each $S^1$. 
In the same way, decompose $f$ into $f_0+f_{\perp}$ ; $|f_0(t,z')|$, $|f_{\perp}(t,z')|\leq C\|\eta\|_{\Gamma^0_{\beta}}e^{-\beta t}$ still hold with a possibly bigger $C$ still depending only on $\mathcal{P}$. 
Then $\varphi_0$ verifies $(\partial_t^2-\partial_t)\varphi_0=f_0(t,z')$, and we solve this writing: 
   \begin{equation}\label{integral_formula_1}
    \varphi_0(t,z')=a_{\mathcal{P}}(z')+\int_t^{+\infty}e^{t'}dt'\int_{t'}^{+\infty}e^{-t''}f_0(t'',z')dt''
   \end{equation}
(with the fact that $\varphi_0$ is $L^2$ for $e^{-t}dt\vol^D$ to get rid of an additional term $\chi(z')e^t$). 
Notice that for each $(t,z')$, the double integral is in absolute value less than $\tfrac{1}{\beta(1+\beta)}\|v(\cdot,z')e^{\beta\cdot}\|_{\Gamma^0}e^{-\beta t}$, and in particular $\varphi_0(t,z')$ tends to $a_{\mathcal{P}}(z')$ exponentially fast when $t$ goes to infinity. Moreover we can write 
$a_{\mathcal{P}}(z')=\varphi_0(0,z')-\int_0^{+\infty}e^{t'}dt'\int_{t'}^{+\infty}e^{-t''}f_0(t'',z')dt''$ for all $z'$, which gives a $C^0$ bound on $a_{\mathcal{P}}$ depending only on the polydisc and $\|\eta\|_{C^{0,\alpha}_{\beta}}$.

We still have to deal with $\varphi_{\perp}$. Since it is orthogonal to the constants on every circle, we can write $|\varphi_{\perp}(t,\theta,z')|\leq \pi^2 \sup_{\theta'\in[0,2\pi]}\big|\tfrac{\partial^2 \varphi_{\perp}}{\partial\theta^2}(t,\theta',z')\big|$ for every $(t,\theta,z')$, and this can be rewritten as:
 \begin{equation*}
  |\varphi_{\perp}(t,\theta,z')|\leq \pi^2 e^{-2t}\sup_{\theta'\in[0,2\pi]}\Big|\Big(-\frac{\partial^2 \varphi_{\perp}}{\partial t^2}+\frac{\partial \varphi_{\perp}}{\partial t}+f_{\perp}\Big)(t,\theta',z')\Big|
 \end{equation*}
for every $(t,\theta,z')$ thanks to the equation verified by $\varphi_{\perp}$. The $\sup$ in the latter right hand side is smaller than $C\|\eta\|_{C^{0,\alpha}_{\beta}}$ with $C$ depending only on the polydisc, thanks to the $C^2$ estimate on $\varphi$ and hence on $\varphi_{\perp}$ we got from Proposition \ref{lem_thm_weighted_ddbar_lemma1}. 

We can sum this up saying that on $\mathcal{P}$, $|\varphi-a_{\mathcal{P}}|\leq C_{\mathcal{P}}\|\eta\|_{\Gamma^{0,\alpha}_{\beta}}e^{-\beta't}$, with $|a_{\mathcal{P}}|\leq C_{\mathcal{P}}\|\eta\|_{\Gamma^{0,\alpha}_{\beta}}$ (and $a_{\mathcal{P}}$ independent of $z_1$ and continuous).

 \item \textit{There is a crossing in} $\mathcal{P}$. Write again $z_1=e^{-e^t+i\theta}$, and $z''=(z_3,\dots,z_m)$. Nothing impedes us to lead the same analysis as above in $\mathcal{P}$ but outside of the $(m-1)$-dimensional polydisc $\{z_2=0\}$ (we can still write down the integrals and take suprema) to see once more that 
$|\varphi-a_{\mathcal{P}}|\leq C_{\mathcal{P}}\|\eta\|_{\Gamma^{0,\alpha}_{\beta}}e^{-\beta't}$, with $|a_{\mathcal{P}}(z_2,z'')|\leq C_{\mathcal{P}}\|\eta\|_{\Gamma^0_{\beta}}$, $z_2\neq0$, and $a_{\mathcal{P}}$ continuous outside of $\{z_2=0\}$.
 \end{enumerate}
We want to improve our analysis of the crossing case, but so far we can collect the information we got when working around $D_1$. Notice that the $a_{\mathcal{Q}}$ patch together; indeed, if $x\in D_1$ has coordinates ${z'}^{(1)}$ in $\mathcal{Q}_1$ and ${z'}^{(2)}$ in $\mathcal{Q}_2$ (and up to increasing the number of polydiscs we can assume $\mathcal{Q}_1\cap\mathcal{Q}_2$ has nonempty interior), we have that
 \begin{equation*}
  a_{\mathcal{Q}_1}\big({z'}^{(1)}\big)=\lim\limits_{\substack{y\rightarrow x\\ y\in\mathcal{Q}_1\backslash D_1}}\varphi(y)
                                       =\lim\limits_{\substack{y\rightarrow x\\ y\in(\mathcal{Q}_1\cap\mathcal{Q}_2)\backslash D_1}}\varphi(y)
                                       =\lim\limits_{\substack{y\rightarrow x\\ y\in\mathcal{Q}_2\backslash D_1}}\varphi(y)
                                       =a_{\mathcal{Q}_2}\big({z'}^{(2)}\big).
 \end{equation*}
Let us denote this function induced on $D_1$ by $a$. Thus $a$ is continuous and bounded on $D_1\backslash \bigcup_{l=2}^j D_l$. It is furthermore pluriharmonic. Indeed, choose a test $(m-2,m-2)$-form $\chi$ on $D_1\backslash \bigcup_{l=2}^j D_l$, or more precisely on some polydisc in $D_1$ from which we remove $\bigcup_{l=2}^j D_l$. 
There is no loss in generality in assume that this polydisc is $\mathcal{Q}\cap D_1$ for one the $\mathcal{Q}$ considered above. It the sense of currents, $\langle i\partial_{z'}\overline{\partial}_{z'}\psi,\chi\rangle= \int_{D_1\cap\mathcal{Q}} \psi i\partial_{z'}\overline{\partial}_{z'}\chi$. Now since we are working on $\mathcal{Q}$, we can write this as $\int_{D_1\cap\mathcal{Q}} \psi_{\mathcal{Q}} i\partial_{z'}\overline{\partial}_{z'}\chi$, and this rewrites 
$\int_{\{\log(-\log|z_1|)=t\}}\varphi(z_1,z')i\partial_{z'}\overline{\partial}_{z'}\chi+O(e^{\beta't})$. This latter integral is equal to $\int_{\{\log(-\log|z_1|)=t\}}\eta'\wedge\chi$ by Stokes' theorem, where $\eta'=\sum_{p,q\geq2}\eta_{p\bar{q}}dz_pd\overline{z_q}$. Clearly this is $O(e^{-\beta t})$, so finally  $\langle i\partial_{z'}\overline{\partial}_{z'}\psi,\chi\rangle=0$.

Finally from the fact that for a Poincaré metric, the domain of the Laplacian $L^2\rightarrow L^2$ is $H^2$, we get that $a$ is $H^2$ on $D_1$, and an integration by parts says it is constant.
Let us go back to the polydisc $\mathcal{P}$ with a crossing, and write $z_2=e^{-e^s+i\zeta}$. 
We now know there is a constant $c_1$ such that $|\varphi(z_1,z_2,z'')-c_1|\leq C_1\|\eta\|_{C^{0,\alpha}_{\beta}}e^{-\beta' t}$. 
By symmetry, there is a constant $c_2$ such that $|\varphi(z_1,z_2,z'')-c_2|\leq C_2\|\eta\|_{C^{0,\alpha}_{\beta}}e^{-\beta' s}$. Since $t$ and $s$ are arbitrarily big, this forces $c_1$ and $c_2$ to be equal. 
In the general case, this shows that the constant induced by $\varphi$ is the same on all the irreducible components of a common connected component $\mathcal{D}_k$ of $D$, and we can write, if this constant is denoted by $a_k$, that there exists a constant $C$ such that on a fixed neighbourhood of $\mathcal{D}_k$, $|\varphi-a_k|\leq C\|\eta\|_{\Gamma^{0,\alpha}_{\beta'}}\big(\sum_{l=1}^{j}\big|\log|\sigma_l|\big|\big)^{-\beta'}$. 
The lemma is proved. \hfill $\blacksquare$

~

\noindent\textit{End of the proof of Proposition \ref{prop_weighted_ddbar_lemma}}. 
We are now performing a last improvement to our controls to get the result of Proposition \ref{prop_weighted_ddbar_lemma}. 
First, assume $\beta\leq 2$, so that $\beta'=\beta$. Then, Schauder estimates on balls of quasi-coordinates with mixed weights give us a constant $C$ such that 
 \begin{equation*}
  \Big\|\sum_j\big|\log|\sigma_j|\big|^{\beta'}\Big(\varphi-\sum_{k=1}a_k v_k\Big)\Big\|_{C^{2,\alpha}}\leq C\|\eta\|_{C^{0,\alpha}_{\beta}}.
 \end{equation*}
To prove the theorem (at order $(2,\alpha)$), one needs to change the sum of weights $\big|\log|\sigma_j|\big|^{\beta'}$ into a product of these weights. 
Notice that this is automatic away from the crossings of $D$, so that we can focus on what happens near a crossing of the component $\mathcal{D}_k$, which is assumed to have codimension 2. 
As in the proof of Lemma \ref{lem_thm_weighted_ddbar_lemma2}, we consider a polydisc $\mathcal{P}=\{|z_1|\leq\tfrac{1}{e}\}\times\cdots\times\{|z_m|\leq\tfrac{1}{e}\}$ around a crossing given by $\{z_1=0\}\cup\{z_2=0\}$, and we set 
$z_1=e^{-e^t+i\theta}$, $z_2=e^{-e^s+i\zeta}$ and $z''=(z_3,\dots,z_m)$. Take a punctured disc $\{z_1\}\times\{0<|z_2|\leq\tfrac{1}{e}\}\times\{z''\}$, on which we can write, according to the proof of Lemma \ref{lem_thm_weighted_ddbar_lemma2}
 \begin{equation*}
  \varphi(z_1,s,\zeta,z'')=a_k+\int_s^{+\infty}e^{s'}ds'\int_{s'}^{+\infty}e^{-s''}f_0(z_1,s'',z'')ds''+\varphi_{\perp_2}(z_1,s,\zeta,z'')
 \end{equation*}
where $f=|z_2|^2\log^2(|z_2|^2)\eta\big(\tfrac{\partial}{\partial z_2}, \tfrac{\partial}{\partial \overline{z_2}}\big)$ and $f_0$ is its $\zeta$-invariant part (so that one has $|f_0|,|f|\leq C\|\eta\|_{C^0}e^{-\beta(t+s)}$, $C$ depending only on $\mathcal{P}$), and $\varphi_{\perp_2}$ orthogonal to the constants on each circle $\{s=\text{constant}\}$ verifying 
 \begin{equation*}
  \frac{\partial^2 \varphi_{\perp_2}}{\partial s^2}-\frac{\partial \varphi_{\perp_2}}{\partial s}+e^{2s}\frac{\partial^2 \varphi_{\perp_2}}{\partial \zeta^2}=v_{\perp_2}.
 \end{equation*}
so that $\tfrac{\partial^2 \varphi_{\perp_2}}{\partial \zeta^2}=e^{-2s}\big(f_{\perp_2}-\tfrac{\partial^2 \varphi_{\perp_2}}{\partial s^2}+\tfrac{\partial \varphi_{\perp_2}}{\partial s}\big)$, which is, among others, smaller than $C\|\eta\|_{C^{0,\alpha}_{\beta}}e^{-\beta'(t+s)}$ thanks to the $C^2$ control we have on $\big|\log|z_1|\big|^{\beta'}(\varphi-a_k)$, hence the $C^0$ control on $e^{\beta't}\big(\frac{\partial^2 \varphi_{\perp_2}}{\partial s^2}-\frac{\partial \varphi_{\perp_2}}{\partial s}\big)$. 
Then conclude using $\big|\varphi_{\perp_2}(z_1,s,\zeta,z'')\big|\leq \pi^2\sup_{\zeta'\in[0,2\pi]}\big|\tfrac{\partial^2 \varphi_{\perp_2}}{\partial \zeta^2}(z_1,s,\zeta',z'')\big|$.

When $\beta$ is $>2$, observe that on the $S^1$-invariant parts we have the desired control. 
Now, for the orthogonal part, after using Schauder estimates we have a $C^{2,\alpha}_{\beta'}$ control on $(\varphi-a_k)$. 
Applying the technique above first give a $C^{0}$ control on $\prod_j\big|\log|\sigma_j|\big|^2\big(\sum_j\big|\log|\sigma_j|\big|\big)^{\min\{\beta-2,2\}}(\varphi-a_k)$, hence a $C^{2,\alpha}$ control, and applying it once more gives a $C^0_{\min\{\beta,4\}}$, hence a $C^{2,\alpha}_{\min\{\beta,4\}}$ control. 
Just repeat the argument as many times as necessary to get a $C^0_{\beta}$ control, and conclude with weighted Schauder estimates. \hfill $\square$

 \section{Proof of the Calabi-Yau theorem on $\XD$}  \label{prf_approx_CY}

As underlined in our comment following the statement of Theorem \ref{thm_approx_CY}, which is new here decomposes into to parts: uniform unweighted bounds on the approximate solutions, to which we devote the next section, and the fact that approximate solutions lie in weighted Hölder spaces, to which we devote section \ref{prf_approx_CY_wtd_spaces}. 
Moreover, an easy observation on the statement of Theorem \ref{thm_approx_CY} is that it is enough to have the announced uniform $C^k$ bounds on the $\varphi_{\vareps}$ to get $\varphi$ with a diagonal extraction, so we are only looking for such bounds, and not for possible weighted $C^k$ bounds. 

 \subsection{Uniform bounds} \label{prf_approx_CY_unif_bnds}

Before starting, notice that when $X$ is a Riemann surface, the theorem follows at once from Lemma \ref{lem_pert_lapl} and Proposition \ref{prop_weighted_ddbar_lemma}. 
In this section and in the following, we are thus assuming that $m\geq 2$ (nonetheless Lemma \ref{lem_sob_inj} also holds if $m=1$).

  \subsubsection{Order zero estimate}  \label{prf_approx_CY_C0_est}
To get a $C^0$ estimate, we follow a Moser's iteration scheme. 
Nonetheless, it will be more convenient to work also on the normalized potentials $\psi_{\vareps}:=\varphi_{\vareps}-a_{\vareps}$ with $a_{\vareps}:=\tfrac{1}{\vl}\int_{\XD}\varphi_{\vareps}\vol^{\omega'}$. 
In what follows, all the $L^p$ norms are taken with respect to $\vol^{\omega'}$ unless another measure is specified. 
Similarly, $g$ refers to the Riemannian metric $\omega'(\cdot,J\cdot)$, $\nabla$ to its Levi-Civita connection, and Hölder norms of functions and tensors are computed with respect to $\omega'$.
 \begin{prop}  \label{prop_initial_est}
 In the conditions of Theorem \ref{thm_approx_CY}, there exist constants $C$ and $A$ only depending on $\omega'$ and $f$ such that for all $\vareps\in(0,1]$, $\|\psi_{\vareps}\|_{L^2}\leq C$ and $|a_{\vareps}|\leq A$. 
In particular, $\|\varphi_{\vareps}\|_{L^2}\leq C'=C+A\vl$.
 \end{prop}
\prf. We start with the $L^2$ estimate. Fix $\vareps\in(0,1]$, and set $\omega'_{\vareps}=\omega'+i\ddbar\varphi_{\vareps}$ and $T_{\vareps}:=(\omega')^{m-1}+(\omega')^{m-2}\wedge(\omega'_{\vareps})+\cdots+(\omega'_{\vareps})^{m-1}$. 
Notice that $T_{\vareps}$ is closed, greater than or equal to $(\omega')^{m-1}$ in the sense of real $(m-1,m-1)$-forms and that $i\ddbar\varphi_{\vareps}\wedge T_{\vareps}=(\omega'_{\vareps})^m-(\omega')^m=(e^{f+\vareps\varphi_{\vareps}}-1)(\omega')^m$. 
Now, since 
 \begin{equation*}
  i\ddbar\big(\varphi_{\vareps}^2T_{\vareps}\big)=2i\partial\varphi_{\vareps}\wedge\dbar\varphi_{\vareps}\wedge T_{\vareps}+2\varphi_{\vareps}i\ddbar\varphi_{\vareps}\wedge T_{\vareps},
 \end{equation*}
and $\int_{\XD}i\ddbar\big(\varphi_{\vareps}T_{\vareps}\big)=0$ (Gaffney-Stokes), we have:
 \begin{equation*}
  \int_{\XD}i\partial\varphi_{\vareps}\wedge\dbar\varphi_{\vareps}\wedge T_{\vareps}+\int_{\XD}\varphi_{\vareps}(e^{f+\vareps\varphi_{\vareps}}-1)(\omega')^m=0.
 \end{equation*}
Noticing that $i\partial\varphi_{\vareps}\wedge\dbar\varphi_{\vareps}\wedge T_{\vareps}\geq i\partial\varphi_{\vareps}\wedge\dbar\varphi_{\vareps}\wedge(\omega')^{m-1}$ and that $(e^{f+\vareps\varphi_{\vareps}}-1)\varphi_{\vareps}=e^f(e^{\vareps\varphi_{\vareps}}-1)\varphi_{\vareps}+(e^f-1)\varphi_{\vareps}\geq (e^f-1)\varphi_{\vareps}$ (because $e^{\vareps\varphi_{\vareps}}-1$ has the sign of $\varphi_{\vareps}$), we get that:
 \begin{equation*}
  \int_{\XD}i\partial\varphi_{\vareps}\wedge\dbar\varphi_{\vareps}\wedge(\omega')^{m-1}\leq \int_{\XD}\varphi_{\vareps}(1-e^f) (\omega')^m.
 \end{equation*}
Since $\int_{\XD}(1-e^f) (\omega')^m=0$ and $\partial\varphi_{\vareps}=\partial\psi_{\vareps}$, this rewrites $\int_{\XD}i\partial\psi_{\vareps}\wedge\dbar\psi_{\vareps}\wedge(\omega')^{m-1}\leq \int_{\XD}(1-e^f)\psi_{\vareps} (\omega')^m$. 
The left-hand-side term of the latter inequality is $\|d\psi_{\vareps}\|^2_{L^2_{\omega'}}$. 
From unweighted Poincaré inequality (\ref{unwtd_poinc_ineq}) ($\psi_{\vareps}$ has zero mean) and Cauchy-Schwarz inequality, it then follows that $\|\psi_{\vareps}\|_{L^2}\leq C_P\|1-e^f\|_{L^2}:=C$, which is independent of $\vareps$.

Let us now estimate $a_{\vareps}$, and to begin with, see how to get an upper bound on it. 
Integrating both parts of the equation $\big(\omega'+i\ddbar\varphi_{\vareps}\big)^m=e^{f+\vareps\varphi_{\vareps}}(\omega')^m$ yields $\int_{\XD}e^{f+\vareps\varphi_{\vareps}}\vol^{\omega'}=\vl$.
Hence Jensen inequality says $\int_{\XD}\vareps\varphi_{\vareps}e^{f}\vol^{\omega'}\leq 0$ i.e. $\int_{\XD}\varphi_{\vareps}e^{f}\vol^{\omega'}\leq 0$. 
Now $a_{\vareps}\vl=\int_{\XD}e^f(\varphi_{\vareps}-\psi_{\vareps})\vol^{\omega'}\leq -\int_{\XD}e^f\psi_{\vareps}\vol^{\omega'}$ so by Cauchy-Schwarz,
 \begin{equation*}
  a_{\vareps}\leq \frac{\|\psi_{\vareps}\|_{L^2}\|e^f\|_{L^2}}{\vl}
             \leq \frac{C\|e^f\|_{L^2}}{\vl},
 \end{equation*}
which does not depend on $\vareps$. 

On the other hand, in order to get a lower bound on $a_{\vareps}$, let us define $b_{\vareps}:=\tfrac{1}{\vl}\int_{\XD}\varphi_{\vareps}\vol^{\omega'_{\vareps}}$, i.e. $b_{\vareps}$ is the mean of $\varphi_{\vareps}$ with respect to the metric it defines. 
This way, $\vl=\int_{\XD}e^f\vol^{\omega'}=\int_{\XD}e^{-\vareps\varphi_{\vareps}}\vol^{\omega'_{\vareps}}\geq \int_{\XD}(1-\vareps\varphi_{\vareps})\vol^{\omega'_{\vareps}}=(1-\vareps b_{\vareps})\vl$, so $b_{\vareps}\geq0$. 
Therefore, $a_{\vareps}\geq a_{\vareps}-b_{\vareps}=\tfrac{1}{\vl}\int_{\XD}\varphi_{\vareps}(\vol^{\omega'}-\vol^{\omega'_{\vareps}}) =\tfrac{1}{\vl}\int_{\XD}\big(1-e^{f+\vareps\varphi_{\vareps}}\big)\psi_{\vareps}\vol^{\omega'}\geq -\tfrac{\|1-e^{f+\vareps\varphi_{\vareps}}\|_{L^2}\|\psi_{\vareps}\|_{L^2}}{\vl}$. 
To conclude, repeat $\|\psi_{\vareps}\|_{L^2}\leq C$, and use $|\vareps\varphi_{\vareps}|\leq |f|_{C^0}$ to get $a_{\vareps}\geq -\tfrac{2C\|f\|_{C^0}e^{2\|f\|_{C^0}}}{\vl^{1/2}}$, which is again independent of $\vareps$. \hfill $\square$

~

Next proposition is central in our upcoming iteration scheme:
 \begin{prop}  \label{prop_ibp_ineq}
  Under the assumptions of Theorem \ref{thm_approx_CY}, for all $\vareps\in(0,1]$ and all $p\geq 2$ we have:
   \begin{equation}   \label{ibp_ineq}
    \int_{\XD}i\partial\big(|\varphi_{\vareps}|^{p/2}\big)\wedge\dbar\big(|\varphi_{\vareps}|^{p/2}\big)\wedge(\omega')^{m-1}
                   \leq \frac{p^2}{4(p-1)}\int_{\XD}|\varphi_{\vareps}|^{p-2}\varphi_{\vareps}(1-e^f) (\omega')^m.
   \end{equation}
 \end{prop}
\prf. Fix $p\geq 2$ and $\vareps\in(0,1]$. 
Again from the inequalities $T_{\vareps}\geq (\omega')^{m-1}$ and $(1-e^{f+\vareps\varphi_{\vareps}})|\varphi_{\vareps}|^{p-2}\varphi_{\vareps}= (1-e^f)|\varphi_{\vareps}|^{p-2}\varphi_{\vareps}+e^{f}(1-e^{\vareps\varphi_{\vareps}})|\varphi_{\vareps}|^{p-2}\varphi_{\vareps}\leq (1-e^f)|\varphi_{\vareps}|^{p-2}\varphi_{\vareps}$, the proposition is proved if we show the identity
 \begin{equation*}
  \int_{\XD}i\partial\big(|\varphi_{\vareps}|^{p/2}\big)\wedge\dbar\big(|\varphi_{\vareps}|^{p/2}\big)\wedge T_{\vareps}
                   = \frac{p^2}{4(p-1)}\int_{\XD}|\varphi_{\vareps}|^{p-2}\varphi_{\vareps}(1-e^{f+\vareps\varphi_{\vareps}}) (\omega')^m.
 \end{equation*}
But this simply follows from the direct computation 
 \begin{equation*}
  i\ddbar\big(|\varphi_{\vareps}|^{p-1}\varphi_{\vareps}T_{\vareps}\big)=
                    p|\varphi_{\vareps}|^{p-2}\varphi_{\vareps}i\ddbar\varphi_{\vareps}\wedge T_{\vareps}
                   +p(p-1)|\varphi_{\vareps}|^{p-2}i\partial\varphi_{\vareps}\wedge\dbar\varphi_{\vareps}\wedge T_{\vareps},
 \end{equation*}
the identities $i\ddbar\varphi_{\vareps}\wedge T_{\vareps}=(e^{f+\vareps\varphi_{\vareps}}-1)(\omega')^m$, $|\varphi_{\vareps}|^{p-2}i\partial\varphi_{\vareps}\wedge\dbar\varphi_{\vareps}= \tfrac{4}{p^2}i\partial\big(|\varphi_{\vareps}|^{p/2}\big)\wedge\dbar\big(|\varphi_{\vareps}|^{p/2}\big)$ and $\int_{\XD}i\ddbar\big(|\varphi_{\vareps}|^{p-1}\varphi_{\vareps}T_{\vareps}\big)=0$ by Gaffney-Stokes. \hfill $\square$

~

Before deriving inductive controls in our iteration scheme, we have to enlighten which Sobolev embedding we are going to use, and in particular between which (weighted) Sobolev spaces: 
 \begin{df}
  Let $q\in[1,+\infty)$. 
  We set:
   \begin{equation*}
    L^{q,0}_{0}=\Big\{v\in L^q_{loc}|\,\int_{\XD}|v|^q\rho\vol^{\omega'}<+\infty\Big\}=L^q_{\rho\vol^{\omega'}}
   \end{equation*}
and if $k\geq\N$, we call $L^{q,k}_{0}$ the space of functions $v\in L^{q,k}_{loc}$ such that $\big|\nabla^j v\big|\in L^{q,0}_{0}$ for $j=0,\dots,k$.
 \end{df}

 \begin{lem}  \label{lem_sob_inj}
  Let $q_2\geq q_1$. 
  Then one has the continuous injection $L^{q_1,1}_{0}\hookrightarrow L^{q_2,0}_{0}$ as soon as $\tfrac{1}{q_1}\leq \tfrac{1}{q_2}+\tfrac{1}{2m}$.
 \end{lem}
\prf. We only need to look at what happens near the divisor, and even near crossings, since the smooth divisor case is ruled in \cite{biq}, Lemma 4.4. 
We assume for simplicity that the crossings have codimension two in $X$. 
Let us consider a small polydisc $U$ around a point in such a crossing, and let us cover $U\backslash D$ by a union (with the notations of \ref{asymp_prop})
 \begin{equation*}
  \bigcup_{k,l\geq 0} \Phi_{\delta_{k,l}}(\mathcal{P}), \quad \mathcal{P}=\big(\tfrac{3}{4}\Delta\big)^2\times \Delta^{m-2},
 \end{equation*}
where the $\delta_{k,l}=(\delta^1_k,\delta^2_l)\in(0,1)^2$ can be chosen so that 
$1-\delta_k^1\sim\tfrac{1}{2^k}$ (resp. $1-\delta_l^1\sim\tfrac{1}{2^l}$) when $k$ (resp. $l$) goes to $\infty$ (we choose them in a similar way as the $\delta_k$ of the proof of Lemma \ref{lem_sobolev_estimate}). In this way, ${\Phi_{\delta_{k,l}}}^*\rho_1$ (resp. ${\Phi_{\delta_{k,l}}}^*\rho_2$) is mutually bounded with $2^k$ (resp. $2^l$), i.e. ${\Phi_{\delta_{k,l}}}^*\rho$ is mutually bounded with $2^{k+l}$.
We can also assume that the metric on $U\backslash D$ is the product $g_U$ of two standard cusp metrics by a euclidian metric, so that all the ${\Phi_{\delta_{k,l}}}^*g_U$ give the same metric on $\mathcal{P}$. 
We also have $|w|_{L^p(U\backslash D)}^p\sim \sum_{k,l\geq 0}\frac{1}{2^{k+l}}\big\|{\Phi_{\delta_{k,l}}}^*w\big\|_{L^p}^{p}$ for any $p\geq1$.
Now we know that there exists $C>0$ such that for all $w\in L^{q_1,1}(\mathcal{P})$, $\|w\|_{L^{q_2,0}(\mathcal{P})}\leq C\|w\|_{L^{q_1,1}(\mathcal{P})}$.

Take $v\in L^{q_1,1}_{0}$; then
 \begin{align*}
  \hspace{0.3cm}|v|_{L^{q_2,0}_{0}(U\backslash D)}^{q_2}
                      &\sim \sum_{k,l\geq 0}\frac{1}{2^{k+l}}\big\|{\Phi_{\delta_{k,l}}}^*(v\rho^{1/q_2})\big\|_{L^{q_2}}^{q_2}
                       \sim \sum_{k,l\geq 0}\frac{1}{2^{k+l}}2^{k+l}\big\|{\Phi_{\delta_{k,l}}}^*v\big\|_{L^{q_2}}^{q_2}\\
                      &\leq C\sum_{k,l\geq 0}\big\|{\Phi_{\delta_{k,l}}}^*v\big\|_{L^{q_1,1}}^{q_2}
                       \leq C\Big(\sum_{k,l\geq 0}\big\|{\Phi_{\delta_{k,l}}}^*v\big\|_{L^{q_1,1}}^{q_1}\Big)^{\tfrac{q_2}{q_1}}
                       \text{ since } q_2\geq q_1\\
                      &\sim C\Big(\sum_{k,l\geq 0}\frac{1}{2^{k+l}}\big\|{\Phi_{\delta_{k,l}}}^*(v \rho^{1/q_1}) \big\|_{L^{q_1,1}}^{q_2}
                                \Big)^{\tfrac{q_2}{q_1}}
                       \sim C |v|_{L^{q_1,1}_{0}(U\backslash D)}^{q_2}.   \hspace{2cm} \square
 \end{align*}

Let us come back to our iteration scheme. Set $\epsilon=\min\big\{\tfrac{3}{2}, \tfrac{m}{m-1}, 1+\nu\big\}>1$ (beware $\epsilon$ is not related to $\vareps$), so that we have a continuous embedding $L^{2,1}_0\hookrightarrow L^{2\epsilon,0}_0$, of norm $C_{Sob}$ say, according to the latter lemma. 
Let $d\mu$ be the measure $\rho^{1-\epsilon}\vol^{\omega'}$. We have the following inductive control formula:
 \begin{prop}  \label{prop_p_to_espilon_p}
  Under the assumptions of Theorem \ref{thm_approx_CY}, there exist two constants $C$ and $C'$ such that for all $\vareps\in(0,1]$ and all $p\geq 2$,
   \begin{equation*}
    \|\varphi_{\vareps}\|^p_{L^{p\epsilon}_{d\mu}}\leq 
                                C\|\varphi_{\vareps}\|^p_{L^{p}_{d\mu}}+C'p\|\varphi_{\vareps}\|^{p-1}_{L^{p\epsilon}_{d\mu}}.
   \end{equation*}
 \end{prop}
\prf. We are going to use the inequality of Proposition \ref{prop_ibp_ineq}, but first, if $B$ denotes $\sup_{\XD}\big|\tfrac{d\rho}{\rho}\big|_{\omega'}$, an easy computation yields ($p$ and $\vareps$ are fixed):
 \begin{equation*}
  \int_{\XD}\big|d(\rho^{-1/2}|\varphi_{\vareps}|^{p/2})\big|_{\omega'}^2\rho\vol^{\omega'}\leq 2\int_{\XD}\big|d|\varphi_{\vareps}|^{p/2}\big|_{\omega'}^2\vol^{\omega'}+\frac{1}{2}B\int_{\XD}|\varphi_{\vareps}|^p\vol^{\omega'}
 \end{equation*}
so that $\big\|\rho^{-1/2}|\varphi_{\vareps}|^{p/2}\big\|_{L^{2,1}_0}\leq 2\int_{\XD}\big|d|\varphi_{\vareps}|^{p/2}_{\omega'}\big|^2\vol^{\omega'}+(\tfrac{1}{2}B+1)\int_{\XD}|\varphi_{\vareps}|^p\vol^{\omega}$. 
Applying Poincaré inequality (with a mean term) to $|\varphi_{\vareps}|^{p/2}$ we get:
 \begin{equation*}
  \big\|\rho^{-1/2}|\varphi_{\vareps}|^{p/2}\big\|_{L^{2,1}_0}\leq C\int_{\XD}\big|d|\varphi_{\vareps}|^{p/2}\big|_{\omega'}^2\vol^{\omega'}+C'\Big(\int_{\XD}|\varphi_{\vareps}|^{p/2}\vol^{\omega'}\Big)^2
 \end{equation*}
with $C=2+\big(\tfrac{1}{2}B+1\big)C_P$ and $C'=\big(\tfrac{1}{2}B+1\big)\vl^{-1}$. 
Now, $\big(\int_{\XD}|\varphi_{\vareps}|^{p/2}\vol^{\omega'}\big)^2\leq \big(\int_{\XD}|\varphi_{\vareps}|^p\rho^{1-\epsilon}\vol^{\omega'}\big)\big(\int_{\XD}\rho^{\epsilon-1}\vol^{\omega'}\big)$, and this latter integral is finite since $\epsilon-1\leq\tfrac{1}{2}<1$. 
We also know from (\ref{ibp_ineq}) that $\int_{\XD}\big|d|\varphi_{\vareps}|^{p/2}\big|^2_{\omega'}\vol^{\omega'}\leq \tfrac{cp^2}{4(p-1)}\int_{\XD} |1-e^f| |\varphi_{\vareps}|^{p-1}\vol^{\omega'}$, and by Hölder this smaller than  $\tfrac{cp^2}{4(p-1)}\big(\int_{\XD}|\varphi_{\vareps}|^{p}\rho^{1-\epsilon}\vol^{\omega'}\big)^{(p-1)/p}$ $\big(\int_{\XD}(|1-e^f|\rho^{\epsilon-1})^p\rho^{1-\epsilon}\vol^{\omega'}\big)^{1/p}$, and the last factor is always less or equal to 
$\big(\int_{\XD}\rho^{1-\epsilon}\vol^{\omega'}\big)^{1/p}\|1-e^f\|_{C^0_{\nu}}\leq C$ for some $C$ depending only on $\omega'$, $\rho$ and $\epsilon$ (the "parameters"). 
 To sum all this up, we say there are constants $C$ and $C'$ only depending on the parameters so that for all $p>2$ (and for $p=2$ with similar arguments) and all $\vareps\in(0,1]$,
 \begin{equation*}
  \big\|\rho^{-1/2}|\varphi_{\vareps}|^{p/2}\big\|_{L^{2,1}_0}\leq
                                     C\|\varphi_{\vareps}\|_{L^p_{d\mu}}^p+C'p\|\varphi_{\vareps}\|_{L^p_{d\mu}}^{p-1}.
 \end{equation*}
Applying the Sobolev embedding stated in Lemma \ref{lem_sob_inj} to $\rho^{-1/2}|\varphi_{\vareps}|^{p/2}$, we exactly get that $\|\varphi_{\vareps}\|_{L^{p\epsilon}_{d\mu}}^p$ is less than or equal to $C_{Sob}^2$ times the left-hand side of the latter inequality, so finally up to renaming the constants there are $C$ and $C'$ only depending on the parameters such that for all $p\geq 2$ and $\vareps\in(0,1]$,
 \begin{equation*}
  \hspace{4cm} \|\varphi_{\vareps}\|_{L^{p\epsilon}_{d\mu}}^p \leq      
                      C\|\varphi_{\vareps}\|_{L^p_{d\mu}}^p+C'p\|\varphi_{\vareps}\|_{L^p_{d\mu}}^{p-1}. \hspace{4cm} \square
 \end{equation*}

~

Since under the conditions of Theorem \ref{thm_approx_CY} we have an initial estimate on $\|\varphi_{\vareps}\|_{L^2}$ and hence on $\|\varphi_{\vareps}\|_{L^2_{d\mu}}$ ($\vol^{\omega'}$ dominates $d\mu$) independent of $\vareps$, it is now an easy exercise to show that there exists two positive constants $Q$ and $C_1$ depending only on the parameters such that
 \begin{equation*}
  \|\varphi_{\vareps}\|_{L^p_{d\mu}}\leq Q\big(C_1p\big)^{-m/p}
 \end{equation*}
for all $p\geq 2$ and $\vareps\in(0,1]$. Letting $p$ go to $\infty$, we have thus proved:
 \begin{prop}[Uniform $C^0$ estimate] \label{prop_CY_C0_est}
  Under the assumptions of Theorem \ref{thm_approx_CY}, there exists $Q=Q(\nu, \omega', \|f\|_{C^0_{\nu}})$ such that for all $\vareps\in(0,1]$, $\|\varphi_{\vareps}\|_{C^0}\leq Q$. 
 \end{prop}

  \subsubsection{Second order estimate} 

We are now looking for second order estimates, which as usual when dealing with Monge-Ampère equations derive from the $C^0$ estimate. 
If we denote by $\Delta$ (resp. $\Delta_\vareps$) the Laplacian of $\omega'$ (resp. $\omega'_{\vareps}=\omega+i\ddbar\varphi_{\vareps}$), then Joyce's computation \cite{joy}, p.111 (replace $f$ by $f+\vareps\varphi_{\vareps}$) for Aubin-Yau formula writes: 
 \begin{equation}  \label{yau_formula}
  \begin{aligned} 
   \Delta_{\vareps}(\Delta \varphi_{\vareps})
        =&-2\Delta (f+\vareps\varphi_{\vareps}) +4g^{\alpha\bar{\lambda}}g_{\vareps}^{\mu\bar{\beta}}g_{\vareps}^{\gamma\bar{\nu}}\nabla_{\alpha\bar{\beta}\gamma}\varphi_{\vareps}
          \nabla_{\bar{\lambda}\mu\bar{\nu}}\varphi_{\vareps}\\
         &+4g_{\vareps}^{\alpha\bar{\beta}}g^{\gamma\bar{\delta}}\big({(\riem^{\omega'})^{\bar{\epsilon}}}_{\bar{\delta}\gamma\bar{\beta}}
         \nabla_{\alpha\bar{\epsilon}}\varphi_{\vareps}                                  
         -{(\riem^{\omega'})^{\bar{\epsilon}}}_{\bar{\beta}\alpha\bar{\delta}}\nabla_{\gamma\bar{\epsilon}}\varphi_{\vareps}\big)
  \end{aligned}
 \end{equation}
(factors 2 and 4 are due to the fact that Joyce works with half-Laplacians and $dd^c$ instead of $i\ddbar$). 
Thus Aubin-Yau inequality becomes:
 \begin{prop}  \label{prop_aubin_yau_ineq}
  Under the assumptions of Theorem \ref{thm_approx_CY}, let $\vareps\in(0,1]$ and set $F_{\eta}=\log(2m-\Delta\varphi_{\vareps})-\kappa\varphi_{\vareps}$ where $\kappa$ is some real number to be fixed later. 
Then 
 \begin{equation*}
  \Delta_{\vareps}F_{\vareps}\leq \frac{\Delta f+2\vareps m}{2m-\Delta\varphi_{\vareps}}
                         -\vareps+\kappa\big(2m-(g_{\vareps})^{\alpha\bar{\beta}}g_{\alpha\bar{\beta}}\big)
                         +C(g_{\vareps})^{\alpha\bar{\beta}}g_{\alpha\bar{\beta}}
 \end{equation*}
where $C$ is some constant depending only on $\big\|\riem^{\omega'}\big\|_{C^0}$.
 \end{prop}

 \begin{crl}[Uniform second order estimate]   \label{crl_scd_order_est}
  Under the assumptions of Theorem \ref{thm_approx_CY}, there exits some constant $Q_1=Q_1(\nu, \omega', \|f\|_{C^0_{\nu}},\|f\|_{C^2} )$ such that for all $\vareps\in(0,1]$, $2m-\Delta\varphi_{\vareps}\leq Q_1$. 
This provides in particular two constants $Q_2$ and $c>0$ such that for all $\vareps\in(0,1]$, $\big\|i\ddbar\varphi_{\vareps}\big\|_{C^0}\leq Q_2$ and $c\omega'\leq\omega'_{\vareps}\leq c^{-1}\omega'$.
 \end{crl}
\prf. Choose $\kappa=C+1$ in Proposition \ref{prop_aubin_yau_ineq}. 
Remember that $2m-\Delta\varphi_{\vareps}\geq 2me^{-2\|f\|_{C^0}/m}>0$ (look at the eigenvalues of $i\ddbar\varphi_{\vareps}$ with respect to $\omega'$).
It follows that at any point
 \begin{equation*}
  (g_{\vareps})^{\alpha\bar{\beta}}g_{\alpha\bar{\beta}}=(\kappa-C)(g_{\vareps})^{\alpha\bar{\beta}}g_{\alpha\bar{\beta}}
                                            \leq -\Delta_{\vareps}F_{\vareps}+\frac{\Delta f+2\vareps m}{2m-\Delta\varphi_{\vareps}}-\vareps+2m\kappa
                                            \leq -\Delta_{\vareps}F_{\vareps}+C'
 \end{equation*}
where $C'=2m\kappa+\tfrac{1}{2m}e^{2\|f\|_{C^0}/m}(\|\Delta f\|_{C^0}+2m)$, which is independent of $\vareps$. 
Apply Wu's maximum principle (\cite{wu}, Lemma 3.1) to $\Delta_{\vareps}$ (for fixed $\vareps\in(0,1]$, we know that $\omega_{\vareps}'$ is quasi-isometric to $\omega'$ from \cite{tian-yau1}) and $F_{\vareps}$ and make possibly some extraction to get a sequence $(x_j)$ of points of $\XD$ such that $\lim_{j\rightarrow\infty} F_{\vareps}(x_j)=\sup_{\XD}F_{\vareps}$ and $\lim_{j\rightarrow\infty} \Delta_{\vareps}F_{\vareps}(x_j)\geq 0$. 
Hence up to a shift on the indexes, $(g_{\vareps})^{\alpha\bar{\beta}}g_{\alpha\bar{\beta}}(x_j)\leq C'+1$ for all $j$.
Moreover, play with the eigenvalues to see that $2m-\Delta\varphi_{\vareps}\leq 2\big((g_{\vareps})^{\alpha\bar{\beta}}g_{\alpha\bar{\beta}}\big)^{m-1}e^{f+\vareps\varphi_{\vareps}}$ at any point, so for all $j$ this gives 
$2m-\Delta\varphi_{\vareps}(x_j)\leq 2(C'+1)^{m-1}e^{2\|f\|_{C^0}}$. 
Plug this into the definition of $F_{\vareps}$ to evaluate the $F_{\vareps}(x_j)$, and let $j$ go to $\infty$; this yields $\sup_{\XD} F_{\vareps}\leq 2\|f\|_{C^0}+(m-1)\log(C'+1)+\log2+\kappa\|\varphi_{\vareps}\|_{C^0}$.
Finally, again by definition of $F_{\vareps}$, this tells us that
 \begin{equation*}
  2m-\Delta\varphi_{\vareps}\leq 2(C'+1)^{m-1}\exp(2\|f\|_{C^0}+2\kappa\|\varphi_{\vareps}\|_{C^0}),
 \end{equation*}
which can be made independent of $\vareps$ by noticing we have a uniform bound on $\|\varphi_{\vareps}\|_{C^0}$ by Proposition \ref{prop_CY_C0_est}.  \hfill $\square$

  \subsubsection{Third and higher orders estimates}

We shall now prove:
 \begin{prop}[Uniform third order estimate]  \label{prop_third_order_est}
  Under the assumptions of Theorem \ref{thm_approx_CY}, there exists a constant $Q_3$ such that for all $\vareps\in(0,1]$, $\big\|\nabla i\ddbar\varphi_{\vareps}\big\|_{C^0}\leq Q_3$.
 \end{prop}
 \prf. The starting point is again due to a hard but local computation by Yau \cite{yau1} (see again \cite{aub}). 
Define, for $\vareps\in(0,1]$, $S_{\vareps}$ such that $4S_{\vareps}^2=\big|\nabla i\ddbar\varphi_{\vareps}\big|_{\omega'_{\vareps}}^2$, so that $S^2=(g_{\vareps}')^{\alpha\bar{\lambda}}(g_{\vareps}')^{\mu\bar{\beta}}(g_{\vareps}')^{\gamma\bar{\nu}}\nabla_{\alpha\bar{\beta}\gamma}\varphi_{\vareps} \nabla_{\bar{\lambda}\mu\bar{\nu}}\varphi_{\vareps}$ in local holomorphic coordinates. 
Yau's computation writes: 
 \begin{align*}
  -\Delta_{\vareps}(&S_{\vareps}^2)=\big|\nabla_{\bar{\alpha}\beta\bar{\gamma}\delta}\varphi_{\vareps}-(g_{\vareps})^{\lambda\bar{\mu}}\nabla_{\alpha\bar{\lambda}\gamma}\varphi_{\vareps}\nabla_{\beta\bar{\mu}\delta}\varphi_{\vareps}\big|^2_{\omega'_{\vareps}}\\
               &+\big|\nabla_{\alpha\beta\bar{\gamma}\delta}\varphi_{\vareps}-(g_{\vareps})^{\lambda\bar{\mu}}\nabla_{\alpha\bar{\gamma}\lambda}\varphi_{\vareps}\nabla_{\beta\bar{\mu}\delta}\varphi_{\vareps}
                     -g^{\lambda\bar{\mu}}\nabla_{\alpha\bar{\mu}\delta}\varphi_{\vareps}\nabla_{\lambda\bar{\gamma}\beta}\varphi_{\vareps}\big|^2_{\omega'_{\vareps}}\\
               &+P^{4,2,1}\big((g_{\vareps})^{\alpha\bar{\beta}}, \nabla_{\alpha\bar{\beta}\gamma}\varphi_{\vareps}, \nabla_{\alpha\bar{\beta}}(f+\vareps\varphi_\vareps)\big)            
                 +Q^{4,2,1}\big((g_{\vareps})^{\alpha\bar{\beta}}, \nabla_{\alpha\bar{\beta}\gamma}\varphi_{\vareps},{\riem^a}_{bcd}\big)\\
               &+P^{3,1,1}\big((g_{\vareps})^{\alpha\bar{\beta}}, \nabla_{\alpha\bar{\beta}\gamma}\varphi_{\vareps}, \nabla_{\bar{\alpha}\beta\bar{\gamma}}(f+\vareps\varphi_\vareps)\big)
                 +Q^{3,1,1}\big((g_{\vareps})^{\alpha\bar{\beta}}, \nabla_{\alpha\bar{\beta}\gamma}\varphi_{\vareps},\nabla_e{\riem^a}_{bcd}\big)
 \end{align*}
where the  $P^{j,k,l}$ et $Q^{j,k,l}$ are polynomials with constant universal coefficients in the entries of three matrices, exponents $j$, $k$ and $l$ indicating the degrees of the coefficients of those matrices. 

In view of Corollary \ref{crl_scd_order_est}, there exists a constant $C_1\geq 1$ depending only on the parameters such that for all $\vareps\in(0,1]$, $\Delta_{\vareps}(S_{\vareps}^2)\leq C_1(S_{\vareps}^2+S_{\vareps})$. 
On the other hand, we can use formula (\ref{yau_formula}) to assert there exists constants $c>0$ and $C_2$ depending only on the parameters such that $\Delta_{\vareps}(\Delta\varphi_{\vareps})\leq cS^2-C$ for all $\vareps\in(0,1]$. 
Collect those two inequalities to write $\Delta_{\vareps}(S_{\vareps}^2-2cC_1C_2\Delta\varphi_{\vareps})\leq -C_1\big(S_{\vareps}-\tfrac{1}{2}\big)+C$, $C:=2cC_1C_2+\tfrac{1}{4}C_1$. 
Now choose a sequence of points $(x_j)$ of $\XD$ such that $\lim_{j\rightarrow\infty}(S_{\vareps}^2-2cC_1C_2\Delta\varphi_{\vareps})(x_j)=\sup_{\XD} (S_{\vareps}^2-2cC_1C_2\Delta\varphi_{\vareps})$ and $\lim_{j\rightarrow\infty}\Delta_{\vareps}(S_{\vareps}^2-2cC_1C_2\Delta\varphi_{\vareps})(x_j)\geq0$. 
Then up to a reindexation $\Delta_{\vareps}(S_{\vareps}^2-2cC_1C_2\Delta\varphi_{\vareps})(x_j)\geq-C_1$ so that $\big(S_{\vareps}(x_j)-\tfrac{1}{2}\big)^2\leq 2cC_2+\tfrac{5}{4}$, or $S_{\vareps}(x_j)\leq C_3:= \big(\tfrac{1}{2}+(2cC_2+\tfrac{5}{4})^2\big)$ for all $j$. 
Letting $j$ go to $\infty$, this tells us that $\sup_{\XD}(S_{\vareps}^2-2cC_1C_2\Delta\varphi_{\vareps})\leq C_3+2cC_1C_2\|\Delta\varphi_{\vareps}\|_{C^0}$, hence $\|S_{\vareps}\|_{C^0}\leq ( C_3+4cC_1C_2\|\Delta\varphi_{\vareps}\|_{C^0})^{1/2}$, which can be made independent of $\vareps$ with help of Corollary \ref{crl_scd_order_est}. \hfill $\square$

~

As an immediate consequence, let us state: 
 \begin{crl}
  Let $\alpha\in(0,1)$. Under the assumptions of Theorem \ref{thm_approx_CY}, there exists a constant $Q_{\alpha}$ such that $\|\omega'_{\vareps}\|_{\Gamma^{0,\alpha}(\Lambda^{1,1})}\leq Q_{\alpha}$ for all $\vareps\in(0,1]$.
 \end{crl}

~

Finally, the usual bootstrap argument allows us to conclude. 
Indeed, fix $\alpha\in(0,1)$; in formula \eqref{yau_formula}, the operator $\Delta_{\vareps}$ is uniformly elliptic on a quasi-coordinate system and its coefficients are controlled in $C^{0,\alpha}$, these controls being independent of $\vareps$. 
We have a uniform $C^0$ control on the right-hand side terms and on $\Delta\varphi_{\vareps}$, independent of $\vareps$. 
Thus, using quasi-coordinates, Schauder estimates give a uniform $C^{1,\alpha}$ on $\Delta\varphi_{\vareps}$ independent of $\vareps$. 
Since we have a $C^0$ estimate on $\varphi_{\vareps}$ which does not depend on $\vareps$, we get a $C^{3,\alpha}$ control on $\varphi_{\vareps}$,  independent of $\vareps$. 
Plug this back into formula (\ref{yau_formula}); the operator $\Delta_{\vareps}$ has now its coefficients controlled in $C^{1,\alpha}$, and the right-hand side terms are controlled in $C^{0,\alpha}$, with controls independent of $\vareps$. 
We deduce from those a $C^{4,\alpha}$ control on $\varphi_{\vareps}$, again independent of $\vareps$. 
Going on this induction, we see that for all $k\geq 0$ there exists a $C^{k,\alpha}$ bound on the $\varphi_{\vareps}$ independent of $\vareps$.

 \subsection{The approximate solutions are in weighted spaces} \label{prf_approx_CY_wtd_spaces}

Once we know that one of the potentials $\varphi_{\vareps}$ of Theorem \ref{thm_approx_CY} is in some $C^0_{\gamma}$, $\gamma>0$, fast decay of its derivatives easily follows:
 \begin{prop}
  Under the assumptions of Theorem \ref{thm_approx_CY}, let $\vareps\in(0,1]$. Assume that $\varphi_{\vareps}\in C^0_{\gamma}$ for some $\gamma\in(0,\nu]$. Then $\varphi_{\vareps}\in C^{\infty}_{\gamma}$.
 \end{prop}
\prf. We start by proving that $\varphi_{\vareps}\in C^{1,\alpha}_{\gamma}$ where $\alpha\in(0,1)$ is fixed. 
The statement is local near $D$, so we are looking on what is happening there. 
We take a small polydisc $U=(c\Delta)^k\times \big(\tfrac{1}{2}\Delta\big)^{m-k}$ ($c>0$ small) around a point of a $k$-dimensional crossing, and in which the components of $D$ are given by respective vanishings of the first $k$ variables. 
Now set $\mathcal{P}= \big(\tfrac{1}{2}\Delta\big)^k\times\Delta^{m-k}$ and $\Phi_{\delta}:\mathcal{P}\rightarrow \Delta^m$ as in section \ref{model} for $\delta\in]0,1[^k$, so that $U\backslash D\subset \bigcup_{\delta\in]0,1[^k}\Phi_{\delta}\big(\tfrac{1}{2}\mathcal{P}\big)$. 
We need to estimate $\|\rho^{\gamma}\varphi_{\vareps}\|_{C^{1,\alpha}(U\backslash D)}$, which is comparable to
 \begin{equation*}
  \sup_{\delta\in(0,1)^k}\frac{1}{(1-\delta_1)^{\gamma}\cdots(1-\delta_k)^{\gamma}}\big\|{\Phi_{\delta}}^*\varphi_{\vareps}\big\|_{C^{1,\alpha}(\tfrac{1}{2}\mathcal{P})}.
 \end{equation*}
Now consider on $\mathcal{P}$ the second order operators $P_{\delta}:v\mapsto \tfrac{i\ddbar v\wedge {\Phi_{\delta}}^*[(\omega')^{m-1}+\cdots+(\omega'_{\vareps})^{m-1}]}{{\Phi_{\delta}}^*(\omega')^{m-1}}$ for $\delta\in(0,1)^k$. 
Since $\varphi_{\vareps}\in C^{\infty}(\XD)$, we have uniform ellipticity and uniform $C^l$ control for all $l\geq 0$ on the coefficients of the $P_{\delta}$, meaning for instance that there exists a constant $C$ such that for all $\delta\in(0,1)^k$,
 \begin{equation*}
  \|v\|_{C^{1,\alpha}(\tfrac{1}{2}\mathcal{P})}\leq C\big(\|P_{\delta}v\|_{C^0(\mathcal{P})}+\|v\|_{C^0(\mathcal{P})}\big)
 \end{equation*}
for all $v\in C^2(\mathcal{P})$ such that $\|P_{\delta}v\|_{C^0(\mathcal{P})}$ is finite. 
Now, notice that $P_{\delta}({\Phi_{\delta}}^*\varphi_{\vareps})=1-e^{{\Phi_{\delta}}^*f+\vareps{\Phi_{\delta}}^*\varphi_{\vareps}}$ for all $\delta \in(0,1)^k$. 
But $f\in C^0_{\nu}$ and $\varphi_{\vareps}\in  C^0_{\gamma}$, $\gamma\leq \nu$, so that $\|P_{\delta}({\Phi_{\delta}}^*\varphi_{\vareps})\|_{C^0(\mathcal{P})}\leq C(1-\delta_1)^{\gamma}\cdots(1-\delta_k)^{\gamma}$ for some $C$ independent of $\delta$. 
In the same way, $\|{\Phi_{\delta}}^*\varphi_{\vareps}\|_{C^0(\mathcal{P})}\leq C'(1-\delta_1)^{\gamma}\cdots(1-\delta_k)^{\gamma}$ for some $C'$ independent of $\delta$.
Thus $ \tfrac{1}{(1-\delta_1)^{\gamma}\cdots(1-\delta_k)^{\gamma}}\|{\Phi_{\delta}}^*\varphi_{\vareps}\|_{C^{1,\alpha}(\tfrac{1}{2}\mathcal{P})}$ is controlled independently of $\delta$ i.e. $\|\rho^{\gamma}\varphi_{\vareps}\|_{C^{1,\alpha}(U\backslash D)}$ is finite. 
Take enough of such $U$ to declare that $\|\rho^{\gamma}\varphi_{\vareps}\|_{C^{1,\alpha}(\XD)}$ is finite, that is $\varphi_{\vareps}\in C^{1,\alpha}(\XD)$. 
Reinject this in the previous argument (in the $\|P_{\delta}({\Phi_{\delta}}^*\varphi_{\vareps})\|_{C^{2l-1,\alpha}(\mathcal{P})}$ controls) to get step by step that $\varphi_{\vareps}\in C_{\gamma}^{2l+1,\alpha}(\XD)$ for all $l\geq1$. \hfill $\square$

~

With the latter proposition, in order to complete the proof of Theorem \ref{thm_approx_CY}, we finally have to show:
 \begin{prop}
  Under the assumptions of Theorem \ref{thm_approx_CY}, there exists a constant $c>0$ such that $\varphi_{\vareps}\in C_{c\vareps}^0(\XD)$ for all $\vareps\in(0,1]$. 
 \end{prop}
\prf. Take $\vareps\in(0,1]$. We start from the inequalities $\Delta\varphi_{\vareps}+2\vareps\varphi_{\vareps}\leq 2f$ and $\Delta_{\vareps}\varphi_{\vareps}+2\vareps\varphi_{\vareps}\geq 2f$. 
Take $\gamma\in\big(0,\min(\tfrac{1}{2},\nu)\big)$, and denote by $\mathcal{L}_{\gamma,\vareps}$ the operator $\rho^{\gamma}(\Delta+2\vareps)(\rho^{-\gamma}\cdot)$, so that $\mathcal{L}_{\gamma,\vareps}(\rho^{\gamma}\varphi_{\vareps})\leq2\rho^{\gamma}f\leq M$ for some real constant $M$. 
Denote by $\psi$ a $C^{\infty}(\XD)$ function such that $\mathcal{L}_{\gamma,\vareps}(\psi)=M$ outside a compact subdomain $K$ of $\XD$; such a $\psi$ exists according to the proof of Lemma \ref{lem_pert_lapl}, provided $\gamma\leq c\vareps$.
Now if $A$ is some big enough constant, $v:=\rho^{\gamma}\varphi_{\vareps}-\psi-A\leq 0$ on $\partial K$ and $\mathcal{L}_{\gamma,\vareps}(v)\leq 0$ on the complement $V$ of $K$ in $\XD$. 
We want to deduce from this that $v\leq 0$ on $V$, which would give an upper weighted estimate on $\varphi_{\vareps}$; for this we will use arguments similar to those of the proof of Lemma \ref{lem_control_of_potentials}.

Namely, take an exhaustive increasing sequence $(U_p)_{p\geq0}$ of relatively compact open subsets of $\XD$ containing $K$, and set $V_p=U_p\backslash K$ for all $p$, so that $V=\bigcup_p V_p$. 
Denote for all $p$ by $v_p$ the solution of the Dirichlet problem
 \begin{equation*}
  \left\{\begin{aligned}
          \mathcal{L}_{\gamma,\vareps}(v_p)&=\mathcal{L}_{\gamma,\vareps}(v) & &  \text{ on } V_p\\
                                       v_p &=v                               & &  \text{ on } \partial K\\
                                       v_p &=0                               & &  \text{ on } \partial U_p.
         \end{aligned}
  \right.
 \end{equation*}
Still following the proof of Lemma \ref{lem_control_of_potentials}, we know that it is enough, in order to conclude, to show that the $v_p$ are nonpositive, and that there exists on $\|v_p\|_{L^2(V_p)}$ a bound independent of $p$. 
Let us deal first with the nonpositivity; fix $p$. 
We already know that $v_p$ is nonpositive on the boundaries of its domain $V_p$; suppose it is positive somewhere in $V_p$, and denote by $x\in V_p$ a point such that $v_p(x)=\sup_{V_p} v>0$. 
At this point, $\Delta v_p(x)\geq 0$, whereas 
 \begin{equation*}
  0\geq \mathcal{L}_{\gamma,\vareps}(v_p)(x)=\Delta v_p(x)+\big(2\vareps-\gamma\tfrac{\Delta\rho(x)}{\rho(x)}-2\gamma(\gamma+1)\big|\tfrac{d\rho}{\rho}\big|_{\omega',x}^2\big)v_p(x),
 \end{equation*}
so provided that the parenthesis is $>0$, which is the case with our assumption on $\gamma$, $v_p(x)\leq \tfrac{-\Delta v_p(x)}{\big(2\vareps-\gamma\tfrac{\Delta\rho(x)}{\rho(x)}-2\gamma(\gamma+1)\big|\tfrac{d\rho}{\rho}\big|_{\omega',x}^2\big)}\leq 0$, a contradiction. 

There remains to control $\|v_p\|_{L^2(V_p)}$ independently of $p$. 
In order to do so, we decompose $v_p$ as the sum $\xi_p+\eta_p$, where $\xi_p|_{\partial V_p}\equiv 0$ and $\mathcal{L}_{\gamma,\vareps}(\eta_p)=0$, so we are done if we control $\|\xi_p\|_{L^2(V_p)}$ and $\|\eta_p\|_{L^2(V_p)}$ independently of $p$. 
The arguments above give $\eta_p\leq 0$ and $\inf_{V_p}\eta_p=\inf_{\partial V_p}\eta_p=\inf_{\partial K}\eta_p=\inf_{\partial K} v$, so $\|\eta_p\|_{L^2(V_p)}\leq \vl(V)^{1/2}\big|\inf_{\partial K} v\big|$. 

Finally, an integration by parts gives:
 \begin{equation*}
  \int_{V_p}\xi_p\mathcal{L}_{\gamma,\vareps}(\xi_p)\vol^{\omega'}
  =\int_{V_p}\big(2\vareps-\gamma^2\big|\tfrac{d\rho}{\rho}\big|_{\omega'}^2\big)\xi_p^2\vol^{\omega'} + \int_{V_p}|d\xi_p|^2\vol^{\omega'}.
 \end{equation*}
Up to reducing the constant $c$, $2\vareps-\gamma^2\big|\tfrac{d\rho}{\rho}\big|_{\omega'}^2\geq 0$ on $\XD$; moreover we have seen that $\int_{V_p}|d\xi_p|^2\vol^{\omega'}\geq \tfrac{\vl(K)}{C_P\vl}\int_{V_p}\xi_p^2\vol^{\omega'}$, so $\int_{V_p}\xi_p^2\vol^{\omega'}\leq \tfrac{C_P\vl}{\vl(K)}\int_{V_p}\xi_p\mathcal{L}_{\gamma,\vareps}(\xi_p)\vol^{\omega'}$. 
Now notice that $\int_{V_p}\mathcal{L}_{\gamma,\vareps}(\xi_p)^2\vol^{\omega'}=\int_{V_p}\mathcal{L}_{\gamma,\vareps}(v)^2\vol^{\omega'}=\int_{V_p}(\rho^{\gamma}(\Delta\varphi_{\vareps}+2\vareps\varphi_{\vareps})-M)^2\vol^{\omega'}\leq \int_{V\backslash D} (\rho^{\gamma}(\Delta\varphi_{\vareps}+2\vareps\varphi_{\vareps})-M)^2\vol^{\omega'}<+\infty$, since $\Delta\varphi_{\vareps}+2\vareps\varphi_{\vareps}$ is bounded and $\rho^{\gamma}$ is square integrable as $\gamma<\tfrac{1}{2}$, and conclude by using Cauchy-Schwarz inequality. 
We have proved that $\varphi_{\vareps}\leq C\rho^{-\gamma}$ on $V$, for any $\gamma\in[0,c\vareps)$ and for some $C$, possibly depending on $\vareps$ and $\gamma$.

The reverse inequality $\Delta_{\vareps}\varphi_{\vareps}+2\vareps\varphi_{\vareps}\geq 2f$ gives us near $D$ the weighted lower bound $\varphi_{\vareps}\geq -C\rho^{-\gamma}$ for any $\gamma\in[0,c\vareps)$, up to reducing $c$. 
This is done by working with $\omega'_{\vareps}$ instead of $\omega'$. 
Nonetheless, we can indeed take $c$ independent of $\vareps$ because $\omega'_{\vareps}$ and $\omega'$ are mutually bounded independently of $\vareps$ by Corollary \ref{crl_scd_order_est}. \hfill $\square$

~

The proof of Theorem \ref{thm_approx_CY} is complete, and this ends the present part.

\section[Uniqueness of constant scalar curvature metrics]{Uniqueness of constant scalar curvature metrics ($K[D]$ ample)}  \label{unqnss_cscK}

 \subsection{Statement of the result}

As an application of both our constructions of approximate geodesics and metrics with negative Ricci forms, we shall get to the following result:
 \begin{thm} \label{thm_unqnss_csck_met}
  Assume $K[D]$ is ample. If there exists $\omega'\in\mom$ such that its scalar curvature $\scal(\omega')$ is constant on $\XD$, then it is unique in $\mom$.
 \end{thm}

The proof is the purpose of the next section. 
For now, we state and show an auxiliary result, which will be useful at the end of this demonstration, and which also explains the uniqueness stated in the latter theorem. 
Indeed, it specifies that the group of automorphisms of $X$ and tangent to the divisor is discrete, and consequently why we get a proper uniqueness for a Poincaré type Kähler metric on $\XD$ with constant scalar curvature, and not only uniqueness up to the action of such automorphisms in the connected component of the identity. 
 \begin{lem}  \label{lem_holo_vf}
  Assume $K[D]$ is ample. Then the space of holomorphic vector fields which are $L^2$ with respect to some Poincaré type metric is reduced to 0.
 \end{lem}
\prf. Endow $\XD$ with Tian-Yau's Kähler-Einstein metric \cite{tian-yau1}, or more generally with any $\varpi$ in any $\mom$ such that $\varrho_{\varpi}\leq -c\varpi$ for some $c>0$; in any case, denote the metric by $\omega$. 
Let $Z$ be a $L^2_{\omega}$  holomorphic vector field, which, as such, is in $\Gamma^{\infty}_{loc}(T^{1,0})$. 
Since $\omega$ dominates any metric smooth through $D$, $Z$ is actually smooth on the whole $X$. 
The rest of the proof is the same as in the compact case; we indeed have $\Delta_{\omega}|Z|_{\omega}^2=\ric_{\omega}(Z,Z)-2|\nabla_{\omega} Z|^2\leq -c|Z|^2$, so that $\Delta_{\omega}|Z|_{\omega}^2 +\tfrac{c}{2}|Z|_{\omega}^2\leq 0$. 
Now the integration by parts 
 \begin{equation*}
  \int_{\XD} \big|d|Z|_{\omega}^2\big|^2\vol^{\omega}=\int_{\XD}|Z|_{\omega}^2\Delta_{\omega}|Z|_{\omega}^2\vol^{\omega}
                                                      \leq -c\|Z\|^4_{L^4_{\omega}}
 \end{equation*}
forces $|Z|^2_{\omega}$ to be constant, hence to vanish since then $\Delta_{\omega}|Z|_{\omega}^2=0$.  \hfill $\square$

 \subsection{Proof of the uniqueness theorem} \label{prf_unqnss_thm}
 
This proof follows really closely Chen's proof in \cite[\S 6]{chen1} for the compact case, and we give here an outline of it for the sake of completeness. 

Let us fix a few notations. 
We denote by $\omega$ a metric of $\mom$ such that $\varrho_{\omega}\leq -c\omega$, $c>0$, given by Theorem \ref{thm_negative_ricci}; we consider it as the base-point of $\mom$, so that Kähler potentials will be computed with respect to this $\omega$, i.e. $\varphi$ is a potential for $\omega'\in\mom$ if $\omega'=\omega+i\ddbar\varphi$. 
Take moreover two metrics $\omega_0$ and $\omega_1$ in $\mom$ with constant scalar curvature; we call $v_{\tau}$ the potential associated to $\omega_{\tau}$ such that $\int_{\XD}v_{\tau} \vol^{\omega_{\tau}}=0$, $\tau=0,1$. 
Finally we consider the $\vareps$-geodesic $(v_t^{\vareps})_{t\in[0,1]}$ from $v_0$ to $v_1$ for $\vareps>0$ small; thus we have $v_{\tau}^{\vareps}\equiv v_0$, $\tau=0,1$, and if one sets $f^{\vareps}_t:=\tfrac{\omega^m}{(\omega_{v_t^{\vareps}})^m}$, we have $\ddot{v_t^{\vareps}}-\big|\partial\dot{v_t^{\vareps}}\big|_{\omega_{v_t^{\vareps}}}^2=\vareps f^{\vareps}_t$ for all $t\in [0,1]:=I$. 

Remember that we have on $\ddot{v}_t^{\vareps}$ and $d\dot{v}_t^{\vareps}$, as well as on $\big|i\ddbar v_t^{\vareps}\big|_{\omega}$, uniform bounds on $(\XD)\times I$ independent of $\vareps$. Set $E^{\vareps}:t\mapsto \tilde{\energy}(v_t^{\vareps})$; according to Proposition \ref{prop_second_derivative_energy},
 \begin{equation} \label{second_derivative_energy_eps}
  \ddot{E}^{\vareps}(t)=\int_{\XD}|\mathcal{D}^{\vareps}_t\dot{v}^{\vareps}_t|_{\omega_{v^{\vareps}_t}}^2\vol^{\omega_{v^{\vareps}_t}}
                         - \int_{\XD} \vareps \scal_{v^{\vareps}_t}\vol^{\omega}+\vareps \overline{\scal}\vl
 \end{equation}
for all $t\in I$, where $\mathcal{D}^{\vareps}_t$ stands for $\nabla_{\omega_{v^{\vareps}_t}}^{-}d$.

  \subsubsection{A crucial inequality}

We work now with frozen $t$ on the second summand of the right-hand-side of (\ref{second_derivative_energy_eps}), which we can rewrite as $-\vareps\int_{\XD} f^{\vareps}_t \scal^{\vareps}_t \vol^{\vareps}_t$ with the obvious simplifications of notations. 
Since $\scal^{\vareps}_t=2\Lambda^{\vareps}_t \varrho^{\vareps}_t$ et $\varrho^{\vareps}_t=\varrho_{\omega}+i\ddbar\log(f^{\vareps}_t)$, it follows that 
 \begin{equation*}
  \int_{\XD} f^{\vareps}_t \scal^{\vareps}_t \vol^{\vareps}_t
               =2\int_{\XD}f^{\vareps}_t \big(\tr^{\vareps}_t(\varrho_\omega)-\Delta^{\vareps}_t\log(f^{\vareps}_t)\big) \vol^{\vareps}_t.
 \end{equation*}
Now $\log(f^{\vareps}_t)\in C^{\infty}(\XD)$, so integrating by parts yields $\int_{\XD}f^{\vareps}_t \Delta^{\vareps}_t\log(f^{\vareps}_t) \vol^{\vareps}_t=\int_{\XD} |d\log(f^{\vareps}_t)|_{\omega^{\vareps}_t}^2\vol^{\omega}$. 
Thus, (\ref{second_derivative_energy_eps}) rewrites, after noticing also $\vol^{\omega^{\vareps}_t}=\tfrac{\vol^{\omega}}{f^{\vareps}_t}$, integrating on $I$ and dividing by $\vareps$:
 \begin{equation*}
  \int_{(\XD)\times I}\frac{|\mathcal{D}^{\vareps}_t\dot{v}^{\vareps}_t|_{\omega^{\vareps}_t}^2}{\vareps f^{\vareps}_t}\vol^{\omega} dt
      -2\int_{(\XD)\times I} \big(\tr^{\vareps}_t(\varrho_\omega)-|d\log(f^{\vareps}_t)|_{\omega^{\vareps}_t}^2\big) \vol^{\omega}dt
      =-\overline{\scal}\vl
 \end{equation*}
because $\dot{E}^{\vareps}(0)=\dot{E}^{\vareps}(1)$, as our extremities are potentials of constant scalar curvature metrics. 
Use now the inequality $\varrho_{\omega}\leq -c\omega$ to get:
 \begin{align} \label{crucial_ineq}
   \int_{(\XD)\times I} \Bigg[\frac{|\mathcal{D}^{\vareps}_t\dot{v}^{\vareps}_t|_{\omega^{\vareps}_t}^2}{\vareps f^{\vareps}_t} 
                               +\big(2c\tr^{\vareps}_t(\omega)+2|d\log(f^{\vareps}_t)|_{\omega^{\vareps}_t}^2\big)\Bigg] \vol^{\omega}dt
                                      \leq -\overline{\scal}\vl:=C,
 \end{align}
for all $\vareps>0$. 
This inequality, or rather the three inequalities it contains (every summand in the bracket is nonnegative), are essential in obtaining the controls of the next paragraph.

   \subsubsection{$L^p$ bounds, weak limits, and conclusion} 

Now that we have inequality (\ref{crucial_ineq}) as well as uniform controls independent of $\vareps$ on some second-order derivatives of $(v^{\vareps}_t)$, we can use Chen's computation, and get some control on the following objects: $w^{\vareps}_t := \log(f^{\vareps}_t)$, $X_{t}^{\vareps} :=\sharp_{t}^{\vareps} \partial \dot{v}_{t}^{\vareps}$, $Y^{\vareps}_t:=e^{-w^{\vareps}_t}X^{\vareps}_t$, for all $t\in I$, $\vareps>0$. 
We can sum up those controls this way:
 \begin{lem}
  $(X^{\vareps}_{\cdot})_{\vareps>0}$ is bounded in $L^2\big(|\cdot|_{\omega},\vol^{\omega}dt\big)$, $\big(Y^{\vareps}_{\cdot}\big)_{\vareps>0}$ is bounded in $L^{\infty}\big(|\cdot|_{\omega},\vol^{\omega}dt\big)$, and $\big(\dbar Y^{\vareps}_{\cdot}\big)_{\vareps>0}$ is bounded in $L^{q}(|\cdot|_{\omega},\vol^{\omega}dt\big)$, $1<q<2$. 
  
  Moreover $(w^{\vareps}_{\cdot})_{\vareps>0}$ is bounded in $L^p\big(\vol^{\omega}dt\big)$ for all finite $p \geq 1$, and $\big(e^{-w^{\vareps}_{\cdot}}\big)_{\vareps>0}$ is bounded in $L^{\infty}$ and $\big(\dbar w^{\vareps}_{\cdot}\big)_{\vareps>0}$ is bounded $L^2\big(|\cdot|_{\omega},\vol^{\omega}dt\big)$. 
  
  Finally, $\big(e^{-w^{\vareps}_{\cdot}}\dbar X^{\vareps}_{\cdot}\big)_{\vareps>0}$ tends to 0 in the $L^{q}(|\cdot|_{\omega},\vol^{\omega}dt\big)$, $1<q<2$.
 \end{lem}
\prf. \cite{chen1}, p.225-229. \hfill $\blacksquare$

~

We extract subsequences converging in those respective $L^p$ spaces and denote the limits by replacing the $\vareps$ by $0$. 
This does not lead to some ambiguity, since for instance $\dbar Y^{0}_{\cdot}$ coincides with the weak limit of $\dbar Y^{\vareps}_{\cdot}$.

We want now to show that $\dbar X^0_{\cdot}=0$, from which we are not so far formally, since if everything was smooth we could write $\dbar X^0_{\cdot}=e^{w^0_{\cdot}}\big(\dbar Y^0_{\cdot} + \dbar w^0_{\cdot}\otimes Y^0_{\cdot}\big)=0$. 
To reach this, we have to make a detour by truncated versions of $X^0_{\cdot}$, namely the $X^{0,k}:=\big(\sum_{j=0}^k\frac{(w^0_{\cdot})^j}{j!}\big)Y^0_{\cdot}$ defined for $k\geq 0$.
This provides us that $\dbar X^0_{\cdot}=0$ in the sense of distributions, i.e. for every $\psi$ of the correct type, $\int_{(\XD)\times I}\big(X_t^0,\dbar\psi(t)\big)\vol^{\omega}dt =0$. 
From this we pass to the statement that on almost every slice $(\XD)\times\{t\}$, $\dbar X^0_t=0$; 
since the space of holomorphic $L^2_{\omega}$ vector fields is reduced to 0 by Proposition \ref{lem_holo_vf}, $X_t^0=0$ for those $t$.
Now in an open set of holomorphic coordinates, 
  $\partial\dot{v}_t^{\vareps}=\sum_{j,k=1}^m (g_t^{\vareps})_{j\bar{k}}(X_t^{\vareps})^{\bar{k}}dz_j.$
The right-hand-side term thus tends weakly to 0, as $g_t^{\vareps}$ is bounded independently of $\vareps$. 
Hence $d\dot{v}_t^{\vareps}$ tends weakly to 0 in $L^2_{\omega}$; on the other hand, for every $\vareps>0$ and at any point,
  $\dbar (v_1-v_0)=\dbar v_1^{\vareps}-\dbar v_0^{\vareps}=\int_0^1 \dbar \dot{v}_t^{\vareps} dt,$
hence for all $(2m-1)$-form $\psi$ with compact support in $\XD$,
 \begin{equation*} \label{dbar_v_0_v_1}
  \big(d(v_1-v_0),\psi\big)=\int_0^1 d\dot{v}_t^{\vareps}\wedge\psi \,dt
                            =\int_{(\XD)\times I}d\dot{v}_t^{\vareps}\wedge\psi \,dt
                            =\big\langle d\dot{v}_t^{\vareps}, \tilde{\psi}\big\rangle,
 \end{equation*}
where $\tilde{\psi}(\cdot,t)=\psi$ for all $t\in[0,1]$. Letting $\vareps$ go to 0, this says that $d(v_1-v_0)=0$ is the sense of distributions, hence vanishes since it is locally smooth. 
This implies $v_0=v_1$ up to a constant, hence $\omega_0=\omega_1$ (and $v_0=v_1$ by normalization). The reader is referred to \cite{chen1}, p.229-231, for the details. \hfill $\square$

\begin{footnotesize}

 \renewcommand{\refname}{References}

\end{footnotesize}

~

\small \textsc{École Normale Supérieure, UMR 8553}

\end{document}